\documentclass[english,a4paper,10pt]{amsart}
\usepackage{babel,amsmath,amssymb,amsbsy,amsfonts,latexsym,amsthm, amscd,amsxtra}
\usepackage{graphicx}
\usepackage[all]{xy}
\usepackage[T1]{fontenc}
\usepackage{graphics,latexsym}
\newtheorem{Teo}{Theorem}[section]
\newtheorem{Def}[Teo]{Definition}
\newtheorem{Cor}[Teo]{Corollary}
\newtheorem{Lem}[Teo]{Lemma}
\newtheorem{Prop}[Teo]{Proposition}

\newtheorem{Rema}[Teo]{Remark}
\newenvironment{Rem}{\begin{Rema} \begin{upshape}} {\end{upshape}\end{Rema}}
\newtheorem*{Pf}{Proof}
\newenvironment{Proof}{\begin{Pf} \begin{upshape}} {\end{upshape} \qed\end{Pf}}

\newcommand\beq[1]{ \begin{equation} \label{#1}}
\newcommand{\eeq}{ \end{equation} }
\newcommand{\beqno}{ \begin{equation*}}
\newcommand{\eeqno}{ \end{equation*}}
\newcommand\beqa[1]{ \begin{eqnarray} \label{#1}}
\newcommand{\eeqa}{ \end{eqnarray} }
\newcommand{\beqano}{ \begin{eqnarray*} }
\newcommand{\eeqano}{ \end{eqnarray*} }

\def\be{\begin{equation}}
\def\ee{\end{equation}}
\def\bea{\begin{eqnarray}}\def\eea{\end{eqnarray}}
\def\nn{\nonumber}

\newcommand{\T}{ {\mathbb T}   }
\newcommand{\N}{ {\mathbb N}   }
\newcommand{\R}{ {\mathbb R}   }
\newcommand{\Z}{ {\mathbb Z}   }

\newcommand{\A}{ {A}   }

\def\TTT{\T}%{\hbox{\msytw T}}
\def\RRR{\R}%{\hbox{\msytw R}}
\renewcommand \a {\alpha}
\newcommand \e {\varepsilon }
\newcommand \x {\xi}
\renewcommand \b  {\beta}
\renewcommand \r {\rho}
\renewcommand \d {\delta}

\newcommand \m {\mu}
\newcommand \n {\nu}
\newcommand \om {\omega}

\newcommand \f {\varphi}
\newcommand \F {\Phi}
\newcommand \g {\gamma}

\newcommand \s {\sigma}

\renewcommand \t {\tau}

\renewcommand \l {\lambda}
\renewcommand \L {\Lambda}
\newcommand \ph {\phi}

\newcommand \Le {L_{\eta}}
\newcommand \Fk {\phi_{\eta,k}}
\newcommand \Fc {\phi_{\eta,c(\Le)}}
\newcommand \Fac {\phi_{\eta,\a(c)}}

\newcommand \cA {{\mathcal A}}
\newcommand \cN {{\mathcal N}}
\newcommand \cS {{\mathcal S}}

\newcommand \cC {{\mathcal C}}
\newcommand \cB {{\mathcal B}}
\newcommand \cI {{\mathfrak I}}
\newcommand \cW {{\mathcal W}}
\newcommand \calM {\mathfrak{M}}

\newcommand \cM {{\mathcal M}}
\newcommand \cbA {\bar{{\mathcal A}}}
\newcommand \bd {\bar{\d}}
\newcommand \cL {{\mathcal L}}
\newcommand \cE {{\mathcal E}}

\newcommand \dpr {\partial}

\newcommand{\tu}{\tilde{u}}

\newcommand \rH {{\rm H}}
\newcommand \rT {{\rm T}}

\def\TT{{\mathcal T}}

\let\0=\noindent

\def\ie{\hbox{\it i.e.,\ }}
\def\eg{\hbox{\it e.g.,\ }}
\let\dpr=\partial

\let\==\equiv

\author{Alfonso Sorrentino}
\title[Lecture Notes on Mather's theory for Lagrangian systems]{Lecture Notes on Mather's theory for Lagrangian systems \footnote{Preliminary version: please send all comments and typos to the author.}}
\address{Department of Pure Mathematics and Mathematical Statistics, University of Cambridge, Wilberforce Road, Cambridge CB3 0WB, United Kingdom} 
\email{as998@dpmms.cam.ac.uk}
\date{\today}
\begin{document}
\selectlanguage{english}
\maketitle
\tableofcontents

%%%%%%%%%%%%%%%%%%%%%%%%%%%%%%%%%%%%%%%%%%%%%%%%%%%%%%%%%%%%%%%%%%%%%%%%%%%%
%%%%%%%%%%%%%%%%%%%%%%%% CHAPTER 0 - INTRODUCTION %%%%%%%%%%%%%%%%%%%%%%%%%%%%%%%%%%%
%%%%%%%%%%%%%%%%%%%%%%%%%%%%%%%%%%%%%%%%%%%%%%%%%%%%%%%%%%%%%%%%%%%%%%%%%%%%

\section{Introduction}\label{intro}

The celebrated Kolmogorov-Arnol'd-Moser (or KAM) theorem finally settled the old question concerning the existence of {\it quasi-periodic}  motions for {\it nearly-integrable} Hamiltonian systems, \ie Hamiltonian systems that are slight perturbation of an {\it integrable} one.  In the integrable case, in fact, the whole phase space is foliated by invariant {\it Lagrangian} submanifolds that are diffeomorphic to tori, and on which the dynamics is conjugated to a rigid rotation. These tori are generally called {\it KAM tori} (see also Definition \ref{kamtorus}). On the other hand, it is natural to ask  what happens to such a foliation and to these {\it stable motions} once the system is perturbed.
In 1954 Kolmogorov \cite{KolmogorovKam} - and later Arnol'd  \cite{ArnoldKAM} and Moser \cite{Mosertwist} in different contexts - proved that, in spite of the generic disappearence of the invariant 
tori filled by periodic orbits, already pointed out by Henri Poincar\'e,  for small perturbations of an integrable system it is still possible to find invariant Lagrangian tori corresponding to  ``{\it strongly non-resonant}'',
\ie {\it Diophantine}, rotation vectors. This result, commonly referred to as {\it KAM theorem}, from the initials of the three main pioneers, besides opening the way to a new understanding of the nature of Hamiltonian systems and their stable motions, contributed to raise new interesting questions, for instance about the destiny of the {\it stable motions} (orbits on KAM tori) that are destroyed by effect of the perturbation or about the possibility of extending these results to systems that are not close to any integrable one.\\
An answer to the first question did not take long to arrive. In 1964 V. I. Arnol'd  \cite{Arnolddiff}  constructed an example of a perturbed integrable system, in which {\it  unstable orbits} - resulting from the breaking of unperturbed KAM tori - coexist with the stable motions drawn by KAM theorem. This
striking, and somehow unexpected, phenomenon, yet quite far from being completely understood, is nowadays called {\it Arnol'd diffusion}.
This new insight led to a change of perspective and in order to make sense of the complex balance between stable and unstable motions that was looming out, new approaches needed to be exploited. Amongst these, variational methods turned out to be particularly successful. 
Mostly inspired by the so-called {\it least action principle}, a sort of widely accepted ``thriftiness'' of  Nature in all its actions, they seemed to provide the natural setting to get over the local view given by  the analytical methods and make towards a global understanding of the dynamics. 

{\it Aubry-Mather theory} represents probably one of the biggest triumphs in this direction. Developed independently by Serge Aubry \cite{AubryLeDaeron} and 
John Mather \cite{Math82} in the eighties, this novel approach to the study of the dynamics of
{\it twist diffeomorphisms of the annulus} (which correspond to Poincar\'e maps of $1$-dimensional Hamiltonian systems \cite{Mosertwist}) pointed out the existence of many {\it action-minimizing sets}, which in some sense generalize invariant rotational curves and that always exist, even after rotational curves are destroyed. Besides providing a detailed structure theory for these new sets, this powerful approach yielded to a better understanding of the destiny of invariant rotational curves and to the construction of interesting chaotic orbits as a result of their destruction \cite{Mather86, MatherForni, Matherorbitsdiffeo}.

Motivated by these achievements, John Mather \cite{Mather91, Mather93} - and later Ricardo Ma\~n\'e \cite{ManeI, ManeII}
and Albert Fathi \cite{Fathibook} in different ways - developed a generalization of this theory to higher dimensional systems.
Positive definite superlinear Lagrangians on compact manifolds, also called {\it Tonelli Lagrangians} (see Definition \ref{defTonelliLag}),
were the appropriate setting to work in.  
Under these conditions, in fact, it is possible to prove the existence of interesting invariant (action-minimizing) sets, known as {\it Mather}, {\it Aubry} and {\it Ma\~n\'e} sets, which generalize KAM tori, and which continue to exist even after KAM tori's disappearance or when it does not make sense to speak of them (for example when the system is ``far'' from any integral one).

Let us remark that these tools revealed also quite promising in the construction of {\it chaotic orbits}, such as for instance {\it connecting orbits} among the above-mentioned invariant sets \cite{Mather93, Bernardconnecting, Contrerasconnecting}. Therefore they set high hopes on the possibility of proving the generic existence of Arnold diffusion in nearly integrable Hamiltonian systems \cite{Matherannouncement}.
However, differently from the case of twist diffeomorphisms, the situation turns out to be more complicated, due to a general lack of information on the topological structure of these action-minimizing sets. 
%This lack represents one of the biggest restraints to the potentiality of such approaches.
% in view of proving the existence of chaotic motions.
%
These sets, in fact, play a twofold role. Whereas on the one hand they may provide an obstruction to the existence of ``diffusing orbits'', on the other hand their topological structure plays a fundamental role in the variational methods that have been developed for the construction of orbits with ``prescribed'' behaviors. We shall not enter further into the discussion of this problematic, but we refer the interested readers to  \cite{Bernardconnecting, Bernardpseudograph, Contrerasconnecting, Matherannouncement, Matherunpublished}.\\

In these lecture notes we shall try to to provide a brief, but hopefully comprehensive introduction to Mather's theory for Lagrangian systems and its subsequent developments by Ricardo Ma\~n\'e and  Albert Fathi. We shall consider only the autonomous case (\ie no dependence on time in the Lagrangian and Hamiltonian). This choice has been made only to make the discussion easier and to avoid some technical issues that would be otherwise involved. However, all the theory that we are going to describe can be generalized, with some ``small'' modifications, to the non-autonomous time-periodic case. Along our discussion, in order to draw the most complete picture of the theory, we shall point out and discuss such differences and the needed modifications.\\

The sections will be organized as follows:
\begin{itemize}
\item Section \ref{sec1.1}: we shall introduce Tonelli Lagrangians and Hamiltonians on compact manifolds and discuss their properties and some examples.
\item Section \ref{cartoonexample}: before entering into the description of Mather's theory, we shall discuss a cartoon example, namely, the properties of invariant probability measures and orbits on KAM tori. This will prepare the ground for understanding the ideas behind Mather's work, as well Ma\~n\'e and Fathi's ones.
\item Section \ref{sec1.2}: we shall discuss the notion of action minimizing measures and introduce the first family of invariant sets: the {\it Mather sets}.
\item Section \ref{sec1.3}: we shall discuss the notion of action minimizing orbits and introduce other two families of invariant sets: the {\it Aubry} and {\it Mane sets}.
\item Section \ref{sec1.4}: we shall discuss Fathi's {\it Weak KAM theory} and its relation to Mather and Ma\~n\'e's works.
\item Some Addenda to the single sections with some complimentary material will be provided along the way.\\

\end{itemize}

\noindent{\it Acknowledgements.}  
These lectures were  delivered at Universit\`a degli Studi di Napoli ``Federico II'' (April 2009, thematic program ``{\it New connections between dynamical systems and Hamiltonian PDEs'}') and at Universitat Polit\`ecnica de Catalunya (June 2010, summer school ``{\it Jornades d'introducci\`o als sistemes din\`amics i a les EDP's}''). I am very grateful to, respectively, Massimiliano Berti, Michela Procesi, Vittorio Coti-Zelati and  Xavier Cabr\'e, Amadeu Delshams, Maria del Mar Gonzales, Tere M. Seara, for their kind invitation. 
I would also like to thank all participants to the courses for their useful comments, stimulating suggestions and careful feedbacks, in particular  Will Merry, Joana Dos Santos and Rodrigo Trevi\~no for a careful reading of a first draft of these notes.\\
A special acknowledgement must go to John Mather, Albert Fathi and Patrick Bernard,  from whom I learnt most of the material here collected and much beyond.
I would also like to thank  all other people that have contributed, in different ways and at different times, to the realization of this project: Luigi Chierchia, Gonzalo Contreras, Rafael de la Llave, Daniel Massart, Gabriel Paternain and many others.

%%%%%%%%%%%%%%%%%%%%%%%%%%%%%%%%%%%%%%%%%%%%%%%%%%%%%%%%%%%%%%%%%%%%%%%%%%%%
%%%%%%%%%%%%%%%%%%%%%%%% CHAPTER 1 - MATHER - FATHI THEORY%%%%%%%%%%%%%%%%%%%%%%%%%%%%%%
%%%%%%%%%%%%%%%%%%%%%%%%%%%%%%%%%%%%%%%%%%%%%%%%%%%%%%%%%%%%%%%%%%%%%%%%%%%%

%\section{Mather-Fathi theory for Lagrangian systems}\label{chap1}

%%%%%%%%%%%%%%%%%%%%%%%% Section 1.1 - Tonelly Lagrangians and Optical Hamiltonians %%%%%%%%%%%%%%%%%%%

\section{Tonelli Lagrangians and Hamiltonians on compact manifolds}\label{sec1.1}

In this section we want to introduce the basic setting that we shall be considering hereafter.
Let $M$ be a compact and connected smooth manifold without boundary.
Denote by ${\rm T}M$ its tangent bundle and ${\rm T}^*M$ the cotangent one. A
point of ${\rm T}M$ will be denoted by $(x,v)$, where $x\in M$ and $v\in
{\rm T}_xM$, and a point of ${\rm T}^*M$ by $(x,p)$, where $p\in {\rm T}_x^*M$ is a
linear form on the vector space ${\rm T}_xM$. Let us fix a Riemannian
metric $g$ on it and denote by $d$ the induced metric on $M$; let
$\|\cdot\|_x$ be the norm induced by $g$ on ${\rm T}_xM$; we shall use the same notation for the norm induced on
${\rm T}_x^*M$.

We shall consider functions $L:\rT M \longrightarrow \R$ of class $C^2$, which are called {\it Lagrangians}. Associated to each Lagrangian, there is a flow on $\rT M$ called the {\it Euler-Lagrange flow}, defined as follows. Let us consider the action functional $\A_L$ from the space of continuous
piecewise $C^1$ curves $\g:[a,b]\rightarrow M$, with $a\leq b$, defined by:
$$
\A_L(\g) := \int_a^b L(\g(t),\dot{\g}(t))\,dt.
$$
Curves that extremize this functional among all curves with the same end-points are solutions of the {\it Euler-Lagrange equation}:
\beqa{ELequation}
\frac{d}{dt}\frac{\dpr L}{\dpr v} (\g(t),\dot{\g}(t)) = \frac{\dpr L}{\dpr x} (\g(t),\dot{\g}(t)) \qquad \forall\,t\in[a,b]\,.
\eeqa
Observe that this equation is equivalent to
$$
\frac{\dpr^2 L}{\dpr v^2} (\g(t),\dot{\g}(t)) \ddot{\g}(t) = \frac{\dpr L}{\dpr x} (\g(t),\dot{\g}(t)) - \frac{\dpr^2 L}{\dpr v \dpr x} (\g(t),\dot{\g}(t))\dot{\g}(t)\,,$$
therefore, if the second partial vertical derivative ${\dpr^2 L}/{\dpr v^2}(x,v)$ is non-degenerate at all points of $\rT M$, we can solve for $\ddot{\g}(t)$. This condition
$$
\det \frac{\dpr^2 L}{\dpr v^2} \neq 0
$$
is called {\it Legendre condition} and allows one to define a vector field $X_L$ on $\rT M$, such that the solutions of 
$\ddot{\g}(t)=X_L(\g(t),\dot{\g}(t))$ are precisely the curves satisfying the Euler-Lagrange equation.
This vector field $X_L$ is called the {\it Euler-Lagrange vector field} and its flow ${\Phi^L_t}$ is the {\it Euler-Lagrange flow} associated to $L$. It turns out that $\Phi^L_t$ is $C^1$ even if $L$ is only $C^2$ (see Remark \ref{ELC1}).

\begin{Def}[{\bf Tonelli Lagrangian}]\label{defTonelliLag}
A function $L:\,{\rm T}M\, \longrightarrow \,\R$ is called a {\it Tonelli
Lagrangian} if:
\begin{itemize}
\item[i)]   $L\in C^2({\rm T}M)$;
\item[ii)]  $L$ is strictly convex in the fibers, in the $C^2$ sense, \ie the second partial vertical derivative
${\dpr^2 L}/{\dpr v^2}(x,v)$ is positive definite, as a quadratic form, for all $(x,v)$;%\in~{\rm T}M$;
\item[iii)] $L$ is superlinear in each fiber, \ie
            $$\lim_{\|v\|_x\rightarrow +\infty} \frac{L(x,v)}{\|v\|_x} = + \infty.$$
            This condition is equivalent to ask that for each $A\in \R$ there exists $B(A)\in\R$ such that
            $$ L(x,v) \geq A\|v\| - B(A) \qquad \forall\,(x,v)\in \rT M\,.$$ 
\end{itemize}
\end{Def}
Observe that since the manifold is compact, then condition {\it iii)} is independent of the choice of the Riemannian metric $g$.

\begin{Rem}\label{remcompleteness}
More generally, one can consider the case of a {\it time-periodic Tonelli Lagrangian}  $L:\,{\rm T}M \times \T\, \longrightarrow \,\R$
(also called {\it non-autonomous case}), as it was originally done by John Mather \cite{Mather91}.  
In fact, as it was pointed out by J\"urgen Moser, this was the right setting to generalize Aubry and Mather's results for twist maps to higher dimensions; in fact,
every twist map can be seen as the time one map associated to the flow of a periodic Tonelli Lagrangian on the one dimensional torus (see for instance \cite{Mosertwist}).
In this case, a further condition on the Lagrangian is needed:
\begin{itemize}
\item [{\it iv)}] {\it The Euler-Lagrange flow is {\it complete}, \ie every maximal integral curve of the vector field $X_L$ has all $\R$ as its domain of definition.}
\end{itemize}
In the non-autonomous case, in fact, this condition is necessary in order to have that action-minimizing curves (or {\it Tonelli minimizers}, see section \ref{sec1.3}) satisfy the Euler-Lagrange equation. Without such an assumption Ball and Mizel \cite{BallMizel} have constructed an example of Tonelli minimizers that are not $C^1$ and therefore are not solutions of the Euler-Lagrange flow. The role of the completeness hypothesis can be explained as follows. It is possible to prove, under the above conditions, that action minimizing curves not only exist %\footnote{As pointed out by Ma\~n\'e in \cite{ManeI} [CHECK!], if the superlinearity condition is satisfied, it is not necessary to assume compactness of $M$ for proving the existence of Tonelli minimizers.} 
and are absolutely continuous, but they are $C^1$ on an open and dense  full measure subset of the interval in which they are defined. It is possible to check that they satisfy the Euler-Lagrange equation on this set, while their velocity goes to infinity on the exceptional set on which they are not $C^1$.
Asking the flow to be complete, therefore, implies that Tonelli minimizers are $C^1$ everywhere and that they are actual solutions of the Euler-Lagrange equation. \\
A sufficient condition for the completeness of the Euler-Lagrange flow, for example, can be expressed in terms of a growth condition for ${\dpr L}/{\dpr t}$:
$$
-\frac{\dpr L}{\dpr t} (x,v,t) \leq C \left( 1 + \frac{\dpr L}{\dpr v}(x,v,t)\cdot v - L(x,v,t) \right) \quad \forall\, (x,v,t)\in \rT M \times \T.
$$
\end{Rem}

{\bf Examples of Tonelli Lagrangians.}
\begin{itemize}
\item {\bf Riemannian Lagrangians.} Given a Riemannian metric $g$ on $\rT M$, the {\it Riemannian Lagrangian} on $(M,g)$ is given by the  {\it Kinetic energy}:
$$
L(x,v) = \frac{1}{2} \| v\|_x^2\,.
$$
Its Euler-Lagrange equation is the equation of the geodesics of $g$: 
$$\frac{D}{dt}\dot{x} \equiv 0\,,$$
and its Euler-Lagrange flow coincides with the geodesic flow.
\item {\bf Mechanical Lagrangians.} These Lagrangians play a key-role in the study of classical mechanics. They are given by the sum of the kinetic energy and a {\it potential} $U: M \longrightarrow \R$:
$$
L(x,v) = \frac{1}{2} \| v\|_x^2 + U(x)\,.
$$
The associated Euler-Lagrange equation is given by:
$$\frac{D}{dt}\dot{x} = \nabla U(x)\,,$$
where $\nabla U$ is the gradient of $U$ with respect to the Riemannian metric $g$, \ie
$$
d_x U \cdot v = \langle \nabla U(x), v \rangle_x \quad \forall\, (x,v)\in \rT M\,.
$$
\item {\bf Ma\~n\'e's Lagrangians.} This is a particular class of Tonelli Lagrangians, introduced by Ricardo Ma\~n\'e in \cite{Maneminmeasure} (see also \cite{FathiFigalliRifford}). If $X$ is a $C^k$ vector field on $M$, with $k\geq 2$, one can embed its flow 
$\varphi^X_t$ into the Euler-Lagrange flow associated to a certain Lagrangian, namely 
$$ L_X(x,v)= \frac{1}{2}\left\| v- X(x) \right\|_x^2\,.$$ 
It is quite easy to check that the integral curves of the vector field $X$ are solutions to the Euler-Lagrange equation. In particular, the Euler-Lagrange flow $\Phi^{L_X}_t$ restricted to ${\rm Graph}(X)=\{(x,X(x)),\;x\in M\}$ (that is clearly invariant) is conjugated to the flow of $X$ on $M$
and the conjugation is given by $\pi|{{\rm Graph}(X)}$,
where $\pi: \rT M \rightarrow M$ is the canonical projection. In other words, the following diagram commutes:
$$\xymatrix{
{{\rm Graph}(X)} \ar@{->}[d]_{\pi} \ar@{->}[r]^{\Phi^{L_X}_t}  & {{\rm Graph}(X)}\ar@{->}[d]^{\pi} 
\\  
{M} \ar@{->}[r]_{\varphi^X_t}  & {M} \\}
 $$
 that is, for every $x\in M$ and every $t\in \R$, $\Phi^{L_X}_t(x,X(x))=(\g_x^X(t),\dot{\g}^X_x(t))$, where $\g_x^X(t)=\varphi_t^X(x)$.\\ 
 %is the solution of the vector field $X$ which is equal to $x$ for $t=0$.\\
\end{itemize}

In the study of classical dynamics, it turns often very useful to consider the associated {\it Hamiltonian system}, which is defined  on the cotangent space $\rT^* M$. 
Let us describe how to define this new system and what is its relation with the Lagrangian one.\\
A standard tool in the study of convex functions is the so-called {\it Fenchel transform}, which allows one to transform functions on a vector space into functions on the dual space (see for instance \cite{Fathibook, Rockafellar} for excellent introductions to the topic).
Given a Lagrangian $L$, we can define the associated {\it Hamiltonian},
as its Fenchel transform (or {\it Fenchel-Legendre transform}): 
%a function on the cotangent bundle:
\beqano H:\; {\rm T}^*M &\longrightarrow & \R \\
(x,p) &\longmapsto & \sup_{v\in {\rm T}_xM} \{\langle p,\,v \rangle_x -
L(x,v)\}\, \eeqano where $\langle \,\cdot,\,\cdot\, \rangle_x$
denotes the canonical pairing between the tangent and cotangent
space.

If $L$ is a Tonelli Lagrangian, one can easily prove that $H$ is finite everywhere (as a consequence of the superlinearity of $L$), $C^2$, superlinear and strictly convex in each fiber (in the $C^2$ sense). Such a Hamiltonian is called {\it Tonelli} (or {\it optical}) {\it Hamiltonian}.

\begin{Def}[{\bf Tonelli Hamiltonian}]\label{opticalHamiltonian}
A function $H:\,{\rm T}^*M\longrightarrow \R$ is called a 
{\it Tonelli (or optical) Hamiltonian} if:
\begin{itemize}
\item[i)]  $H$ is of class $C^2$;
\item[ii)]  $H$ is strictly convex in each fiber in the $C^2$ 
sense, \ie the second partial 
vertical derivative ${\dpr^2 H}/{\dpr p^2}(x,p)$ is positive definite, 
as a quadratic form, for any $(x,p)\in {\rm T}^*M$;
\item[iii)] $H$ is superlinear in each fiber, \ie
$$\lim_{\|p\|_x\rightarrow +\infty} \frac{H(x,p)}{\|p\|_x} = + \infty\,.$$
\end{itemize}
\end{Def}

\vspace{20 pt}

{\bf Examples of Tonelli Hamiltonians.}\\

Let us see what are the Hamiltonians associated to the Tonelli Lagrangians that we have introduced in the previous examples.

\begin{itemize}
\item {\bf Riemannian Hamiltonians.} If $L(x,v) = \frac{1}{2} \| v\|_x^2$ is the Riemannian Lagrangian associated to a Riemannian metric $g$ on $M$, the corresponding Hamiltonian will be
$$
H(x,p) = \frac{1}{2} \| p\|_x^2,
$$
where $\|\cdot\|$ represents - in this last expression - the induced norm on the cotangent space $\rT^* M$. %Hamilton's equations in this case are simply [CONTROLLARE]
%\beqano
%\left\{
%\begin{array}{l}
%\dot{x}(t) = p(t) \\
%\dot{p}(t) = 0\,.
%\end{array}\right.
%\eeqano
%
\item {\bf Mechanical Hamiltonians.} If $L(x,v) = \frac{1}{2} \| v\|_x^2 + U(x)$ is a mechanical Lagrangian, the associated Hamiltonian is:
$$
H(x,p) = \frac{1}{2} \| p \|_x^2 - U(x),
$$
%
%HAMILTON'S EQUATIONS?
that it is sometime referred to as {\it mechanical energy}.
\item {\bf Ma\~n\'e's Hamiltonians.}  If $X$ is a $C^k$ vector field on $M$, with $k\geq 2$, and 
$ L_X(x,v)= \left\| v- X(x) \right\|_x^2$ is the associated Ma\~n\'e Lagrangian, one can check that the corresponding Hamiltonian is given by:
$$
H(x,p)= \frac{1}{2}\|p\|_x^2 + \langle p,X(x) \rangle\,.
$$
\end{itemize}

Given a Hamiltionian one can consider the associated {\it Hamiltonian flow} $\Phi^H_t$ on $\rT^* M$. In local coordinates, this flow can be expressed in terms of the so-called {\it Hamilton's equations}:
\beqa{Hamiltonequation}
\left\{
\begin{array}{l}
\dot{x}(t) = \frac{\dpr H}{\dpr p} (x(t),p(t))\\
\dot{p}(t) = - \frac{\dpr H}{\dpr x} (x(t),p(t))\,.
\end{array}\right.
\eeqa

We shall denote by $X_H(x,p):=\left(\frac{\dpr H}{\dpr p} (x,p),- \frac{\dpr H}{\dpr x} (x,p)\right)$ the {\it Hamiltonian vector field} associated to $H$. This has a more intrinsic (geometric) definition in terms of the canonical symplectic structure $\omega$ on $\rT^*M$. %(see Appendix \ref{sec3.1}). 
In fact, $X_H$ is the unique vector field that satisfies
$$
\omega \left( X_H (x,p), \cdot \right) = d_x H(\cdot)\qquad \forall (x,p)\in\rT^*M.
$$
For this reason, it is sometime called {\it symplectic gradient of} $H$.
It is easy to check from both definitions that - only in the autonomous case - the Hamiltonian is a {\it prime integral of the motion}, \ie it is constant along the solutions of these equations. \\

Now, we would like to explain what is the relation between the Euler-Lagrange flow and the Hamiltonian one. 
It follows easily from the definition of Hamiltonian (and Fenchel transform) that for each $(x,v) \in \rT M$ and $(x,p)\in \rT^* M$ the following inequality holds:
\beqa{Fenchelineq}
 \langle p,\,v \rangle_x \leq L(x,v) + H(x,p)\,.
\eeqa
This is called {\it Fenchel inequality} (or {\it Fenchel-Legendre inequality}) and plays a crucial role in the study of Lagrangian and Hamiltonian dynamics and in the variational methods that we are going to describe. In particular, equality holds if and only if $p = \dpr L/\dpr v (x,v)$. One can therefore introduce the following diffeomorphism between ${\rm T}M$ and ${\rm T}^*M$, known as  {\it Legendre transform}:
\beqa{Legendretransform} \cL:\; {\rm T}M &\longrightarrow & {\rm T}^*M\nonumber\\
(x,v) &\longmapsto & \left(x,\,\frac{\dpr L }{\dpr v}(x,v)
\right)\,. \eeqa
Moreover the following relation with the Hamiltonian holds:
$$ H \circ \cL(x,v) = \left\langle \frac{\dpr L }{\dpr v}(x,v),\,v \right\rangle_x - L(x,v)\,.$$
A crucial observation is that this diffeomorphism $\cL$ represents a conjugation between
the two flows, namely the Euler-Lagrange flow on $\rT M$ and the Hamiltonian flow on $\rT^* M$; in other words, the following diagram commutes:
$$\xymatrix{
{\rT M} \ar@{->}[d]_{\cL} \ar@{->}[r]^{\Phi^L_t}  & {\rT M}\ar@{->}[d]^{\cL} 
\\  
{\rT^* M} \ar@{->}[r]_{\Phi^H_t}  & {\rT^* M} \\}
 $$
 
 \begin{Rem}\label{ELC1}
 Since $\cL$ and the Hamiltonian flow $\Phi^H$ are both $C^1$, then it follows from the commutative diagram above that the Euler-Lagrange flow is also $C^1$.
\end{Rem}
 
Therefore one can equivalently study the Euler-Lagrange flow or the Hamiltonian flow, obtaining in both cases information on the dynamics of the system.  Each of these equivalent approaches will provide different tools and advantages, which  may prove very useful to understand the dynamical properties of the system. For instance, the tangent space is the natural setting for the classical calculus of variations and for Mather and Ma\~n\'e's approaches (sections \ref{sec1.2} and \ref{sec1.3}); on the other hand, the cotangent space is equipped with a canonical symplectic structure, %(see Appendix \ref{sec3.1}), 
which allows one to use several symplectic topological tools, coming from the study of Lagrangian graphs, Hofer's theory, Floer Homology, {\it etc} ... Moreover, a particular fruitful approach in $\rT^*M$ is the so-called {\it Hamilton-Jacobi method} (or {\it Weak KAM theory}), which is concerned with the study of {\it solutions} and {\it subsolutions} of  Hamilton-Jacobi equations. In a certain sense, this approach represents the functional analytical counterpart of the above-mentioned variational approach (section \ref{sec1.4}). In the following sections we shall provide a complete description of these methods and their implications to the study of the dynamics of the system.\\

%%%%%%%%%%%%%%%%%%%%%%%%%%%%%%%%%%%%%%%%%%%%%%%%%%%%%%%%%%%%%%%%%%%%%%%%%%%%%%%%%%%%%%%%%%%%%%

\section{A cartoon example: properties of measures and orbits on KAM tori}\label{cartoonexample}

Before entering into the details of Mather's work, we would like to discuss some properties of orbits and 
invariant probability measures that are supported on {\it KAM tori}. This will provide us with a better understanding of the ideas behind Mather's theory and allow to see in what sense these {\it action-minimizing sets} - namely, what we shall call {\it Mather}, {\it Aubry} and {\it Ma\~n\'e sets} - represent a generalization of KAM tori.

Let $H: \rT^*\T^d \longrightarrow \R$ be a Tonelli Hamiltonian and $L: \rT\T^d \longrightarrow \R$ its associated Tonelli Lagrangian and let us denote by $\Phi^H$ and $\Phi^L$ the respective flows. Observe that one can identify $\rT^*\T^d$ and $\rT \T^d$ with $\T^d\times \R^d$. First of all, let us define what we mean by KAM torus.

\begin{Def}[{\bf KAM Torus}]\label{kamtorus}\rm
$\TT\subset\TTT^d \times\RRR^d$ is a (maximal) KAM torus with 
rotation vector $\r$ if:\\
\0i) $\TT\subset\TTT^d\times\RRR^d$ is a $C^1$ Lagrangian graph, \ie $\TT=\{(x,c+du):\,x\in\T^d\}$, where $c\in \R^d$ and $u: \T^d\longrightarrow \R$;\\
\0ii) $\TT$ is invariant under the Hamiltonian flow $\Phi_t^H$ 
generated by $H$;\\
\0iii) the Hamiltonian flow on $\TT$ is conjugated to a uniform 
rotation on $\TTT^d$; \ie there exists a diffeomorphism 
$\f: \TTT^d\to \TT$
such that $\f^{-1}\circ\Phi_t^H\circ\f= R_\r^t$, $\forall t\in\RRR$,
where $R_\r^t: x \to x+\r t$ $({\rm mod}.\;\Z^d)$. In other words, the following diagram commutes:
$$\xymatrix{
{\TT}  \ar@{->}[r]^{\Phi^H_t}  & {\TT} \\  
{\T^d} \ar@{->}[u]_{\f} \ar@{->}[r]_{R_\r^t}  & {\T^d}\ar@{->}[u]^{\f}  \\}
 $$

\end{Def}

\begin{Rem}\label{remarktori}
({\it i}) First of all, let us remark that, up to some small technicalities in the proofs, most of the results that we shall discuss in this section, will continue to hold if we assume $\TT$ to be only Lipschitz rather than $C^1$.

({\it ii}) If $\r$ is rationally independent, \ie $\r\cdot k \neq 0$ for all $k\in \Z^d\setminus\{0\}$, then each orbit is dense on it (\ie the motion is said to be {\it minimal}). Therefore, it supports a unique invariant probability measure. Consider the case in which $\r$ admits some resonance, \ie there exists  $k\in \Z^d\setminus \{0\}$ such that $\r \cdot k =0$. The set of resonances of $\r$ forms a module over $\Z$. Let $r$ denote its rank. Then, $\TT$ is foliated by an $r$-dimensional family of $(d-r)$-tori, on each of which the motion is minimal.

({\it iii}) If $\r$ is rationally independent then a classical result by Michel Herman implies that this torus is automatically Lagrangian (see \cite[Proposition 3.2]{Herman}). %(see Proposition \ref{Lagrangiantori}).

({\it iv}) Since $\TT$ is Lagrangian and invariant under the Hamiltonian flow, then the Hamiltonian $H$ is constant on it, \ie $H(x,c+du)=E_c$ for some $E_c \in \R$ (it follows from the definition of the Hamiltonian flow). 
In particular, $u$ is a classical solution of Hamilton-Jacobi equation $H(x, c+dv)=k$. It is easy to check that, for a fixed $c$, all solutions have the same energy. In fact, let $u$ and $v$ be two solutions. Then, the function $u-v$ will have at least a critical point $x_0$ and at $x_0$ the two differentials must coincide: $d_{x_0}u=d_{x_0}v$. This implies that 
$H(x_0, c+d_{x_0}u)=H(x_0, c+d_{x_0}v)$. \\
Moreover, $E_c$ is the least possible $k\in\R$ such that there can exist subsolutions of the equation $H(x,c+dv)=k$,  namely $v\in C^1(\T^d)$  satisfying $H(x,c+dv)\leq k$ for all $x\in\T^d$. The proof is essentially the same as above.
In fact, since the function $u-v$ has at least a critical point $x_0$ and at $x_0$ the two differentials coincide, then  $E_c=H(x_0, c+d_{x_0}u)=H(x_0, c+d_{x_0}v) \leq k$. \\
\end{Rem}

Let us start by studying the properties of invariant probability measures that are supported on $\TT$. Let $\m^*$ be an ergodic invariant probability measure supported on $\TT$. 
It will be more convenient at this point to work in the Lagrangian setting. 
Let  us consider the correspondent invariant probability measure for the Euler-Lagrange flow $\Phi^L$, obtained using the Legendre transform: $\m=\cL^*\m^*$, where $\cL^*$ denotes the pull-back.
We would like to point out some properties of this measure and see how $\TT$ can be characterized using measures satisfying such properties.\\

\noindent {\sc \bf Properties}:
\begin{enumerate}
\item First of all, observe that $\int v d\m = \r$. To see this, let us consider the universal cover of $\T^d$, \ie $\R^d$,  and denote by:
\begin{itemize}
\item $\widetilde{\TT}:= \{(q,c+d\tilde{u}(q)),\; q\in\R^d\}$ the lift of $\TT$ to $\R^d$, where $\tilde{u}$ is a periodic extension of $u$;
\item $\widetilde{\Phi^H_t}$ the suspension of the Hamiltonian flow $\Phi^H_t$ to $\R^d\times\R^d$;
\item $\widetilde{\f}$ the lift of $\f$, \ie  $\f: \R^d\to \widetilde{\TT}$ such that  the following diagram commutes:
$$\xymatrix{
{\widetilde{\TT}}  \ar@{->}[r]^{\widetilde{\Phi^H_t}}  & {\widetilde{\TT}} \\  
{\R^d} \ar@{->}[u]_{\tilde{\f}} \ar@{->}[r]_{\widetilde{R}_\r^t}  & {\R^d}\ar@{->}[u]^{\tilde{\f}}  \\}$$
where  $\widetilde{R}_\r^t:q \to q+\r t$. Moreover, $\widetilde{\f}$ has the form
$\widetilde{\f}(q) = (\widetilde{\f}_q(q), c+d\tilde{u}(\widetilde{\f}_q(q)))$ and 
$\widetilde{\f}_q$ satisfies $\widetilde{\f}_q(q+k)=\widetilde{\f}_q(q)+k$ for each $k\in\Z^d$.\\
\end{itemize}

Now, let $(x_0,c+du(x_0))$ be a point in the support of $\m^*$ and consider the corresponding orbit $\g_{x_0}(t):= \pi \Phi^H_t(x_0,c+du(x_0))$, where $\pi: \rT^*\T^d \longrightarrow \T^d$ is the canonical projection. Denote by $\widetilde{\g_{x_0}}(t)$ its lift to $\R^d$. Then, using the ergodic theorem (remember that we are assuming that $\m^*$, and hence $\m$, is ergodic), we have that for a generic point $(x_0,c+du(x_0))$ in the support of $\m^*$ the following holds:
\beqano
\int_{\rT \T^d} v d\m &=&  \lim_{n\rightarrow +\infty} \frac{1}{n} \int_0^n {\dot{\widetilde{\g_{x_0}}}}(t)\,dt = \lim_{n\rightarrow +\infty} \frac{\widetilde{\g_{x_0}}(n)- \widetilde{\g_{x_0}}(0)}{n}\, =\\
&=& \lim_{n\rightarrow +\infty} \frac{\widetilde{\g_{x_0}}(n)}{n}\, =
\lim_{n\rightarrow +\infty} \frac{ \widetilde{\f}_q\left( \widetilde{\f}_q^{-1}\left(\widetilde{\g_{x_0}}(0)\right)+ n\r \right)  }{n}\, =\\
&=& \lim_{n\rightarrow +\infty} \left[ \frac{ \widetilde{\f}_q\left( \widetilde{\f}_q^{-1}\left(\widetilde{\g_{x_0}}(0)\right)+ \{n\r\} \right)}{n} + \frac{[n\r]}{n} \right] = \\
&=& \rho
\eeqano
where $[\cdot]$ and $\{\cdot\}$  denote respectively the integer and fractional part of the vector (component by component).

\begin{Rem}
We shall say that  $\m$ has {\it rotation vector} $\r$. See also Section \ref{sec1.2}. % and Appendix  \ref{sec3.5}
\end{Rem}

\item Let $f:\T^d \longrightarrow \R$ be a $C^1$ function. Then, $\int df(x)\cdot v d\m =0 $. In fact, take $(x_0,v_0) \in {\rm supp}\, \m$ and denote the corresponding orbit by $(x_t,v_t)=\Phi^L_t(x_0,v_0)$. Since $\m$ is invariant, then ${\Phi^L_t}_*\m=\m$ and therefore:
\beqano
\int_{\rT \T^d} df(x)\cdot v d\m &=&  \frac{1}{T} \int_0^Tdt \int_{\rT\T^d} df(x)\cdot v \, d{\Phi^L_t}_*\m = \\
&=& \frac{1}{T} \int_0^Tdt \int_{\T\T^d} df(x_t)\cdot v_t \, d \m  =\\
&=& \frac{1}{T} \int_{\rT\T^d}  d\m \int_{0}^T df(x_t)\cdot v_t \, d t  =\\
&=& \int_{\rT\T^d} \frac{f(x_T) - f(x_0)}{T}\, d\m\, \stackrel{T\rightarrow +\infty}{\longrightarrow} 0\,,
 \eeqano
 since $f$ is bounded. A measure satisfying this condition is said to be {\it closed} (see also Addendum 4.A). Observe that in the above proof we have not used anything else than the invariance of $\m$. See also Proposition \ref{closedmeasure}.\\
 
 \item Recall that $E_c$ is the {\it energy} of  $\TT$, \ie $H(x,c+du)=E_c$. We want to show that this energy value is somehow related to the average value of $L$ on the support of $\m$. Observe in fact that:
 \beqano
\qquad \qquad \int_{\rT\T^d} L(x,v)\,d\m =  \int_{\rT\T^d} \big[L(x,v) - (c+du)\cdot v \big]\,d\m + 
 \int_{\rT\T^d} (c+du)\cdot v\,d\m\,.
 \eeqano
Using the property just pointed out in (2) and observing that, because of the Fenchel-Legendre (in)equality (\ref{Fenchelineq}), along the orbits $(x_t,v_t)$ in the support of $\m$ we have  $(c+du(x_t))\cdot v_t = L(x_t,v_t)+H(x, c+du(x_t)) = L(x_t,v_t)+ E_c $, we can conclude that:
 \beqano
 \int_{\rT\T^d} L(x,v)\,d\m =
 &=& - \int_{\rT\T^d} E_c \,d\m +  \int_{\rT\T^d} c \cdot v\,d\m = \\
&=& - E_c  + c \cdot \int_{\rT\T^d} v\,d\m \;= 
\; - E_c  + c \cdot \r. 
 \eeqano
 In particular, $\int_{\rT\T^d} \big[L(x,v) - c\cdot v \big]\,d\m = -E_c$. Observe that this new Lagrangian ${L_c}(x,v)=L(x,v)-c\cdot v$ is still Tonelli and it is immediate to check that it has the same Euler-Lagrange flow as $L$.\\
 \end{enumerate} 
 
 We can now show this first  result.\\
 
 \begin{Prop}\label{propminimalactiontorus}
 If $\tilde{\m}$ is another invariant probability measure of $\Phi^L$, then 
\beqa{minimalaction}
 \int_{\rT\T^d} L_c(x,v)\, d\tilde{\m} \geq  \int_{\rT\T^d} L_c(x,v)\, d{\m}.
\eeqa
 \end{Prop}
 
\noindent Therefore, $\m$ {\it minimizes} the average value of $L_c$ (or {\it action} of $L_c$) among all invariant probability measures of $\Phi^{L_c}_t$ (or 
$\Phi^{L}_t$ since they are the same). \\ {\sc Warning}: in general, it does not minimize the action of $L\,$! 

\begin{Proof}
The proof is an easy application of Fenchel-Legendre inequality (\ref{Fenchelineq}). Observe that in the support of $\tilde{\m}$, differently from what happens on the support of $\m$ (see property (3) above), this is not an equality, \ie 
$(c+du(x))\cdot v \leq L(x,v)+H(x, c+du(x)) = L(x,v)+ E_c.$ Then:
\beqano
\int_{\rT\T^d} L_c(x,v)\, d\tilde{\m} &=& 
\int_{\rT\T^d} \big[L(x,v) - (c+du)\cdot v \big]\,d\tilde{\m} \geq \\
&\geq& - \int_{\rT\T^d} E_c \,d{\tilde \m}  
= - E_c \;=\; \int_{\rT\T^d} L_c(x,v)\, d{\m}\,.
 \eeqano
\end{Proof}

\begin{Rem}\label{Rem3.5}
({\it i}) It follows from the above proposition that 
$$
E_c = - \min \left\{\int_{\rT\T^d} L_c(x,v)\, d\tilde{\m}:\; \tilde{\m}\;\mbox{is a}\; \Phi^L\mbox{-invariant  probability measure}  \right\}\,.
$$
The minimizing action on the right-hand side is also denoted by $\a(c)$ (see the definition of Mather's $\a$-function in Section \ref{sec1.2}). 

({\it ii}) A measure $\mu$ that satisfies the inequality (\ref{minimalaction}) for all invariant probability measures $\tilde{\mu}$, is called $c$-{\it action minimizing measure} (or {\it Mather's measure with cohomology class} $c$); see Section \ref{sec1.2}.

({\it iii}) It is easy to see from the proof of Proposition \ref{propminimalactiontorus} and Fenchel-Legendre inequality, that if ${\rm supp}\,\tilde{\m}$ is not contained in $\cL^{-1}(\TT)$, then the inequality in (\ref{minimalaction}) is strict.

({\it iii}) In particular, using the previous remark and ({\it ii}) in Remark \ref{remarktori},  we obtain that 
$$
\cL^{-1}(\TT) = \bigcup \{{\rm supp}\, \m: \; \m\;\mbox{is a c-action minimizing measure}\}.
$$
The set on the right-hand side is often denoted by $\widetilde{\cM}_c$ and called {\it Mather set of cohomology class $c$}. It will be defined in Section \ref{sec1.2}.\\
\end{Rem}

Let us now observe that although $\m$ does not minimize the action of $L$ among all invariant probability measures (we have pointed out in fact that it minimizes the action of a modified Lagrangian, with the same Euler Lagrange flow as $L$), then it does minimize it if we put some extra constraints.

\begin{Prop}\label{minactiontorusrotvector}
 If $\tilde{\m}$ is another invariant probability measure of $\Phi^L$ with rotation vector $\rho$ (in the sense of property (1) above) , then 
\beqa{minactionrotvector}
 \int_{\rT\T^d} L(x,v)\, d\tilde{\m} \geq  \int_{\rT\T^d} L(x,v)\, d{\m}.
\eeqa
\end{Prop}

\begin{Proof}
The proof is the same as before, using Fenchel-Legendre inequality (\ref{Fenchelineq}), the fact that $\int v \,d\tilde{\m} = \rho$ and property (3) above. In fact: 
\beqano
\int_{\rT\T^d} L(x,v)\, d\tilde{\m} &=&  \int_{\rT\T^d} \big[L(x,v) - (c+du)\cdot v \big]\,d\tilde{\m} + 
 \int_{\rT\T^d} (c+du)\cdot v\,d\tilde{\m} \geq \\
 &\geq& -E_c + c\cdot \int_{\rT\T^d} v\,d\tilde{\m} =\\
 &=& -E_c + c\cdot \rho = \int_{\rT\T^d} L(x,v)\, d{\m}\,.
 \eeqano
\end{Proof}

\begin{Rem}\label{Remark3.7}
({\it i}) It follows from the above proposition that 
$$
-E_c  + c\cdot \r = \min \left\{\int_{\rT\T^d} L(x,v)\, d\tilde{\m}:\; \tilde{\m}\;\mbox{is a}\; \Phi^L\mbox{-inv.  prob. meas. with rot. vector}\;\r  \right\}\,.
$$
The minimizing action on the right-hand side is also denoted by $\beta(\rho)$ (see the definition of Mather's $\beta$-function in Section \ref{sec1.2}). 

({\it ii}) A measure $\mu$ that satisfies the inequality (\ref{minactionrotvector}) for all invariant probability measures $\tilde{\mu}$ with rotation vector $\rho$, is called {\it action minimizing measure} (or {\it Mather's measure}) {\it with rotation vector} $\r$; see Section \ref{sec1.2}.

({\it iii}) It is easy to see from the proof of Proposition \ref{minactiontorusrotvector} and Fenchel-Legendre inequality, that if $\tilde{\m}$ has rotation vector $\rho$, but its support is not contained in $\cL^{-1}(\TT)$, then the inequality in (\ref{minactionrotvector}) is strict.

({\it iv}) In particular, using the previous remark and  ({\it ii}) in Remark \ref{remarktori},  we obtain that 
$$
\cL^{-1}(\TT) = \bigcup \{{\rm supp}\, \m: \; \m\;\mbox{is an action minimizing measure with rotation vector}\; \r\}.
$$
The set on the right-hand side is often denoted by $\widetilde{\cM}^\r$ and called {\it Mather set of homology class $\rho$}. It will be defined in Section \ref{sec1.2}.\\
\end{Rem}

Summarizing, if $\m$ is any invariant probability measure of $\Phi^L$ supported on $\cL^{-1}(\TT)$, then:
\begin{itemize}
\item[-] $\m$ minimizes the {\it action} of $L(x,v)-c\cdot v$ amongst all invariant probability measures of $\Phi^L$;
\item[-] $\m$ minimizes the {\it action} of $L(x,v)$ amongst all invariant probability measures of $\Phi^L$ with rotation vector $\rho$.\\
\end{itemize}

\vspace{10 pt}

Let us now shift our attention to orbits on KAM tori and see what are the properties that they enjoy. Let $(x_0,c+d_{x_0}u)\in \TT$ a point on the KAM torus and consider its orbit under the Hamiltonian flow, \ie $\g(t)=\pi \Phi^H_t(x_0,c+d_{x_0}u)$, where as usual $\pi:\rT^*\T^d \longrightarrow \T^d$ denotes the canonical projection along the fiber. Let us fix any times $a < b$ and consider the corresponding {\it Lagrangian action} of this curve. Using Fenchel-Legendre (in)equality we get:
\beqano
\int_a^b L(\g(t),\dot{\g}(t))\,dt &=&  \int_a^b \Big((c+d_{\g(t)}u)\dot{\g}(t) - H(\g(t), c+d_{\g(t)}u)\Big)\, dt  = \\
&=& \int_a^b c\,\dot{\g}(t)\, dt + u(\g(b))- u(\g(a)) - E_c (b-a)\,.
\eeqano
Therefore, considering as above the action of the modified Lagrangian $L_c(x,v)$, we get:
\beqa{actionorbittorus}
\int_a^b L_c(\g(t),\dot{\g}(t))\,dt \;=\;  u(\g(b))- u(\g(a)) - E_c (b-a)\,.
\eeqa
Let us now take any other absolutely continuous curve $\xi: [a,b] \longrightarrow \T^d$ with the same endpoints as $\g$, \ie $\xi(a)=\g(a)$ and $\xi(b)=\g(b)$. Proceeding as before and using Fenchel-Legendre inequality, we obtain:
\beqano
\int_a^b L(\xi(t),\dot{\xi}(t))\,dt &\geq&  \int_a^b \Big((c+d_{\xi(t)}u)\dot{\xi}(t) - H(\xi(t), c+d_{\xi(t)}u)\Big)\, dt  = \\
&=& \int_a^b c\,\dot{\xi}(t)\, dt + u(\xi(b))- u(\xi(a)) - E_c (b-a)\,.\eeqano
Hence, using (\ref{actionorbittorus}) and the fact that $\xi(a)=\g(a)$ and $\xi(b)=\g(b)$, we can conclude that:
\beqa{ineqforcurves}
\int_a^b L_c(\g(t),\dot{\g}(t))\,dt \; \leq\;   \int_a^b L_c(\xi(t),\dot{\xi}(t))\,dt.
\eeqa
Therefore for any times $a<b$, $\g$ is the curve that minimizes the action of $L_c$ over all absolutely continuous curves $\xi: [a,b] \longrightarrow \T^d$ with $\xi(a)=\g(a)$ and $\xi(b)=\g(b)$.\\
Actually something more is true. Let us consider a curve with the same endpoints, but a different time-length, \ie $\xi: [a',b'] \longrightarrow \T^d$ with $a'<b'$ and such that $\xi(a')=\g(a)$ and $\xi(b')=\g(b)$. Proceeding as above, one obtains 
\beqano
\int_{a'}^{b'} L_c(\xi(t),\dot{\xi}(t))\,dt &\geq&  u(\xi(b'))- u(\xi(a')) - E_c (b'-a')
\eeqano
and consequently
$$
\int_a^b \Big(L_c(\g(t),\dot{\g}(t)) + E_c\Big) \,dt \; \leq \;
\int_{a'}^{b'} \Big(L_c(\xi(t),\dot{\xi}(t)) + E_c\Big)\,dt\,.
$$
Hence, for any times $a<b$, $\g$ minimizes the action of $L_c+ E_c$  amongst all absolutely continuous curves $\xi$ that connect $\g(a)$ to $\xi(b)=\g(b)$ in any given time (adding a constant does not change the Euler-Lagrange flow).\\
We have just proved the following proposition.

\begin{Prop}\label{proporbittori}
For any given $a<b$, the projection $\g$ of any orbit on $\TT$ minimizes the action of $L_c$ amongst all absolutely continuous curves that connect $\g(a)$ to $\g(b)$ in time $b-a$. Furthermore, $\g$ minimizes the action of $L_c+E_c$ amongst all absolutely continuous curves that connect $\g(a)$ to $\g(b)$ in any given time length.\\
\end{Prop}

\begin{Rem}\label{remmaneset}
({\it i}) A curve $\g: \R \longrightarrow \T^d$ such that for any $a<b$,
$\g|[a,b]$ minimizes the action of $L_c$ amongst all absolutely continuous curves that connect $\g(a)$ to $\g(b)$ in time $b-a$, is called a $c$-{\it global minimizer} of $L$.

({\it ii}) A curve $\g: \R \longrightarrow \T^d$ such that for any $a<b$,
$\g|[a,b]$ minimizes the action of $L_c$ amongst all absolutely continuous curves that connect $\g(a)$ to $\g(b)$ without any restriction on the time length, is called a $c$-{\it time free  minimizer} of $L$.

({\it iii}) Proposition \ref{proporbittori} can be restated by saying that the projection of each orbit on $\TT$ is a $c$-global minimizer of $L$ and a $c$-time free minimizer of $L+E_c$. Moreover, it follows easily from the proof and Fenchel-Legendre inequality, that if the curve $\xi$ does not lie on $\TT$, then the inequality in (\ref{ineqforcurves}) is strict.\\
In particular, we obtain that:
$$\cL^{-1}(\TT)= \bigcup\{(\g(t),\dot{\g}(t)): \; \g\;\mbox{is a c-global minimizer of}\; L \; \mbox{and}\;t\in\R\}.$$
The set on the right-hand side is often denoted by $\widetilde{\cN}_c$ and called {\it Ma\~n\'e set of cohomology class $c$}. It will be defined in Section \ref{sec1.3}.\\
\end{Rem}

Actually,  these curves on $\TT$ are more than $c$-global minimizers. As we shall see, they satisfy a more restrictive condition. In order to introduce it, let us introduce what is called the {\it Ma\~n\'e potential}. 
Let $x_1,x_2 \in \T^d$ and let us denote by $h^T_{c}(x_1,x_2)$ the minimimal action of $L_c$ along curves that connect $x_1$ to $x_2$ in time $T$. We want to consider the infimum of these quantities  for all positive times. 
Of course, without any ``correction'' this infimum might be $-\infty$. Hence, let us consider the action of $L_c+k$ for some $k\in\R$ and define:
\beqa{Manepotentialontorus}
\phi_{c,k}(x_1,x_2) := \inf_{T>0} \left( h^T_{c}(x_1,x_2) + kT\right).
\eeqa
We want to see for which values of $k$ this quantity is well-defined, \ie it is $>-\infty$.\\
Let $\g:[0,T]\longrightarrow \T^d$ be any absolutely continuous curve such that $\g(0)=x_1$ and $\g(T)=x_2$. Then, using Fenchel Legendre inequality:
\beqa{stimapotential}
\int_0^T L_c(\g(t),\dot\g(t))\,dt  + kT &\geq& \int_0^T \left( du(\g(t)\cdot\dot\g(t) -  H(\g(t), c+ d_{\g(t)}u) \right)\,dt + kT = \nonumber\\
&=& u(\g(T)) -u(\g(0)) + (k-E_c)T =\nonumber\\
&=& u(x_2) -u(x_1) + (k-E_c)T\,.
\eeqa
Taking the infimum over all curves connecting $x_1$ to $x_2$, we obtain:
$$
\phi_{c,k}(x_1,x_2) \geq u(x_2) - u(x_1) + (k-E_c)T.
$$
Therefore, if $k\geq E_c$, then $\phi_{c,k}(x_1,x_2) >-\infty$ for all $x_1,x_2\in\T^d$.
On the other hand, it is quite easy to show that for $k< E_c$,  $\phi_{c,k}(x_1,x_2) = -\infty$ for all $x_1,x_2\in \T^d$. Let us observe the following:
\begin{itemize}
\item If $k<E_c$, then $\phi_{c,k}(x,x)=-\infty$ for all $x\in \T^d$. In fact, we know that the orbit $\g(t)=\pi \Phi^H_t(x,c+du(x))$ is recurrent, \ie there exist $T_n \to +\infty$ such that $\d_n:=\|\g(T_n)-x\|\to 0$ as $n\to \infty$. Consider the curve $\xi: [0,T_n+\d_n] \to \T^d$ obtained by joining the curve $\g$ to the unit speed geodesic connecting $\g(T_n)$ to $x$. Then, proceeding as before and using Fenchel-Legendre inequality and the continuity of $u$, we obtain:
\beqano
\phi_{c,k}(x,x) &\leq&  \int_0^{T_n} L_c(\g,\dot\g)dt + \int_{T_n}^{T_n+\d_n} L_c(\x,\dot\x)dt + k(T_n+\d_n) \leq \\
&\leq&u(\g(T_n)) - u(x) + A \d_n + (k-E_c)T_n +k\d_n  \stackrel{n\to+\infty}{\longrightarrow} -\infty,
\eeqano
where $A:=\max_{\|v\|=1} |L_c(x,v)|$.
\item It is also easy to check that for any $k$, $\phi_{c,k}$ satisfies a sort of triangular inequality:
$$
\phi_{c,k}(x,y) \leq \phi_{c,k}(x,z) + \phi_{c,k}(z,y) \qquad \forall \; x,y,z\in\T^d.
$$
For the proof of this, observe that any absolutely continuous curve connecting $x$ to $z$ and any  absolutely continuous curve connecting $z$ to $y$, determine a curve from $x$ to $y$. Passing to the infimum, one gets the result.
\end{itemize}
Our claim follows easily from these two observations. In fact, let $k<E_c$  and $x_1,x_2 \in \T^d$. Then:
\beqano
\phi_{c,k}(x_1,x_2) \leq \phi_{c,k}(x_1,x_1) + \phi_{c,k}(x_1,x_2) = -\infty\,.
\eeqano

Summarizing, we have proved the following proposition.

\begin{Prop}\label{Ecmanecritical}{\rm
\beqano
E_c &=& \inf\{ k \in \R:\; \phi_{c,k}(x,y)>-\infty\; \mbox{for all}\; x,y\in \T^d\} =\\
 &=& \sup\{ k \in \R:\; \phi_{c,k}(x,y)=-\infty\; \mbox{for all}\; x,y\in \T^d\} =\\
&=& \sup \{ k \in \R:\; \phi_{c,k}(x,x)=-\infty\; \mbox{for all}\; x \in \T^d\} =\\
&=& \sup \{ k \in \R:\; \phi_{c,k}(x,x)=-\infty\; \mbox{for some}\; x \in \T^d\}\,.
\eeqano}
\end{Prop}

\noindent Observe that the second and fourth equality are simply a consequence of the triangular inequality.\\

Let us go back to our $c$-global minimizing curves. We have seen (Proposition \ref{proporbittori}) that if
$\g$ is the projection of an orbit on $\TT$, then it is a $c$-global minimizer of $L$ and furthermore it is a 
$c$- time free minimizer of $L+E_c$. The last property can be rewritten in terms of the Ma\~n\'e potential: 
$$
\int_{a}^b \left( L_c(\g(t),\dot\g(t)) + E_c\right) = \phi_{c,E_c}(\g(a),\g(b)) \qquad \mbox{for any}\; a<b.
$$

Therefore these curves seems to be related to the  Ma\~n\'e potential corresponding to $k=E_c$ (the least value for which it is defined). Observe that in this case, we obtain from (\ref{stimapotential}) that
$\phi_{c,E_c}(x,y) \geq u(y)-u(x)$  for all $x,y\in\T^d$ and consequently:
\beqa{stimainversa}
\phi_{c,E_c}(y,x) \geq u(x)-u(y) = - (u(y)-u(x)) \geq - \phi_{c,E_c}(x,y)\,,
\eeqa
which provides a lower bound  on the action needed to ``go back'' from $y$ to $x$. 
In particular, if $\x:[a,b] \to \T^d$ is an absolutely continuous curve, then:
$$
\int_a^b \left(L_c(\x,\dot\x)+ E_c \right)dt \geq \phi_{c,E_c}(\x(a),\x(b)) \geq - \phi_{c,E_c}(\x(b),\x(a))\,.
$$

Question: do there exist curves for which these inequalities are equalities. Observe, that if such a curve exists, then it must necessarily be a $c$-minimizer . 
We shall show that the answer is affirmative and characterize such curves in our case.

\begin{Prop}\label{regularminim}
Let $\g$ be the projection of an orbit on $\TT$. Then, for each $a<b$ we have:
$$
\int_a^b \left(L_c(\g,\dot\g)+ E_c \right)dt = \phi_{c,E_c}(\g(a),\g(b)) = - \phi_{c,E_c}(\g(b),\g(a))\,.
$$
\end{Prop}

\begin{Rem}\label{remregularminim}
({\it i}) A curve $\g: \R \longrightarrow \T^d$ satisfying the conditions in Proposition \ref{regularminim} is called $c$-{\it regular global minimizer} or $c$-{\it static curve} of $L$. Observe that the adjective {\it regular} (coined by John Mather) has no relation to the smoothness of the curve, since this curve will be as smooth as all other solutions of the Euler-Lagrange flow (depending on the regularity of the Lagrangian).

({\it ii}) Since a $c$-regular global minimizer is a $c$-global minimizer (see the comment above), it follows from ({\it iii}) in Remark \ref{remmaneset} that these are all and only the $c$-regular global minimizers of $L$ . Therefore:
$$\cL^{-1}(\TT)= \bigcup\{(\g(t),\dot{\g}(t)): \; \g\;\mbox{is a c- regular global minimizer of}\; L \; \mbox{and}\;t\in\R\}.$$
The set on the right-hand side is often denoted by $\widetilde{\cA}_c$ and called {\it Aubry set of cohomology class $c$}. It will be defined in Section \ref{sec1.3}.\\
\end{Rem}

\begin{Proof}[{\bf Proposition \ref{regularminim}}]
Let $\g(t):=\pi \Phi^H_t(x,c+du(x))$ and consider $a<b$. We already know from Proposition  \ref{proporbittori} that
$$
\int_a^b \left(L_c(\g,\dot\g)+ E_c \right)dt = \phi_{c,E_c}(\g(a),\g(b))\,.
$$
Since the orbit is recurrent, there exist $T_n\to +\infty$ such that 
$\d_n:=\|\g(b+T_n) - \g(a)\| \to 0$.
Consider the curve $\xi: [b,b+T_n+\d_n] \to \T^d$ obtained by joining the curve $\g|[b,b+T_n]$ to the unit speed geodesic connecting $\g(b+T_n)$ to $\g(a)$. 
Then, using Fenchel-Legendre inequality and the continuity of $u$, we obtain:
\beqano
\phi_{c,E_c}(\g(b),\g(a)) &\leq&  \int_b^{b+T_n} L_c(\g,\dot\g)dt + \int_{b+T_n}^{b+T_n+\d_n} L_c(\x,\dot\x)dt + E_c(T_n+\d_n) \leq \\
&\leq&u(\g(b+T_n)) - u(\g(b)) + (A+ E_c)\d_n  \stackrel{n\to+\infty}{\longrightarrow} u(\g(a)) - u(\g(b)),
\eeqano
where $A:=\max_{\|v\|=1} |L_c(x,v)|$.
Therefore, using (\ref{stimapotential}), we obtain:
\beqano
\phi_{c,E_c}(\g(b),\g(a)) &\leq& u(\g(a)) - u(\g(b)) = - (u(\g(b)) - u(\g(a))) =\\
&=& -\phi_{c,E_c}(\g(a),\g(b))
\eeqano
that, together with (\ref{stimainversa}), allows us to deduce that $\phi_{c,E_c}(\g(b),\g(a))  = -\phi_{c,E_c}(\g(a),\g(b))$ and conclude the proof.
\end{Proof}

  \vspace{10 pt}

%%%%%%%%%%%%%%%%%%%%%%Section 1.2 %%%%%%%%%%%%%%%%%%%%%%%%%%%%
\section{Action-minimizing measures: Mather sets}\label{sec1.2}

%AGGIUNGERE QUESTO DA QUALCHE PARTE IN QUESTA SEZIONE

%In order to generalize to more degrees of freedom Aubry and Mather's 
%variational approach to twist maps, a first important notion is that of 
%{\it minimal measure}, which replaces that of action minimizing orbit. 
%Aubry-Mather theory in higher dimension, in fact, cannot deal with such orbits, due to a 
%lack of them: a classical example due to Hedlund 
%\cite{Hedlund} shows the existence of a Riemannian metric on a three-dimensional torus, for which minimal geodesics 
%exist only in three directions. Instead, Mather proposed to look at the 
%closely related notion of action minimizing invariant probability measures. Let us try to describe this idea.

In this section, inspired by what happens for invariant measures and orbits on KAM tori,  we would like to attempt a similar approach for general Tonelli Lagrangians on compact manifolds, and prove the existence of some analogous  compact invariant subsets for the Euler-Langrange (or Hamiltonian) flow. In many senses, these sets will resemble and generalize KAM tori, even when these do not exist or it does not make sense to speak of them.

% and the following sections are meant to provide a comprehensive introduction to Mather's theory for Lagrangian systems and 
%Fathi's weak KAM theory. We shall recall most of the results that we  are going to  use, trying to give - unless where it is needed - general ideas rather  than rigourous proofs (for which we refer to \cite{Fathibook, Mather91} and references therein).

%The setting is the same introduced in Section \ref{sec1.1}. Let $M$ be a compact and connected smooth manifold without boundary of dimension $n$. We denote by $H$ a fixed Tonelli Hamiltonian and $L$ the associated Tonelli Lagrangian.

%This idea represents the keystone of the approach that we are going to describe.

\vspace{10 pt}

Let us start by studying invariant probability measures of the system and their action-minimizing properties. Then, we shall use them to define a first family of invariant sets: the {\it Mather sets} (see Remark \ref{Rem3.5} ({\it iii}) and Remark \ref{Remark3.7} ({\it iii})).\\
Let $\calM(L)$ be the space of probability measures $\mu$ on 
${\rm T}M$  that are invariant  under the Euler-Lagrange flow of $L$ and such that $\int_{{\rm T}M} L\,d\m<\infty$ (finite action). It is easy to see in the case of an autonomous Tonelli Lagrangian, that this set is non-empty. In fact,  recall that because of the conservation of the energy
$E(x,v):=H \circ \cL(x,v) = \left\langle \frac{\dpr L }{\dpr v}(x,v),\,v \right\rangle_x - L(x,v)$ along the motions, each energy level of $E$ is compact (it follows from  the superlinearity condition) and invariant under $\Phi^L_t$. It is a well-known result by Kryloff and Bogoliouboff \cite{KB} that a flow on a compact metric space has at least an invariant probability measure. 

\begin{Prop}\label{KriloffBog}
Each non-empty energy level $\cE(E):= \{E(x,v)=E\}$ contains at least one invariant probability measure of $\Phi^L_t$.
\end{Prop}

\begin{Proof}
Let $(x_0,v_0)\in\cE(E)$ and consider the curve $\g(t):=\pi\Phi^L_t(x_0,v_0)$, where $\pi:\rT M \longrightarrow M$ denotes the canonical projection. For each $T>0$, we can consider the probability measure $\m_T$ uniformly distributed on the piece of curve $(\g(t),\dot{\g}(t))$ for $t\in[0,T]$, \ie 
$$
\int f \,d\m_T = \frac{1}{T}\int_0^T f(\g(t),\dot{\g}(t))\,dt \qquad \forall\; f\in C(\cE(E)).
$$
The family $\{\m_T\}_T$ is precompact in the weak$^*$ topology; in fact the set of probability measures on $\cE(E)$ is contained in the closed unit ball in the space of signed Borel measures of $\cE(E)$, which is weak$^*$ compact by Banach-Alaouglu theorem. Therefore, I can extract a converging subsequence
$\{\m_{T_n}\}_n$, such that $\m_{T_n} \longrightarrow   \m$ in the weak$^*$ topology as ${T_n\rightarrow \infty}$. We want to prove that $\m$ is invariant, \ie for each $s\in\R$ we have ${\Phi^L_s}_*\m=\m$.
Let $f$ be any continous function on $\cE(E)$. Then :
\begin{eqnarray}
\left|\int f d{\Phi^L_s}_*\m - \int f d \m \right| \!\!&\leq&
\!\!\left|\int f d{\Phi^L_s}_*\m - \int f d{\Phi^L_s}_*\m_{T_n} \right| +
\left|\int f d{\Phi^L_s}_*\m_{T_n} - \int f d \m_{T_n} \right|\nonumber\\
&&  \!\!\!+\; \left|\int f d\m_{T_n} - \int f d \m \right| . \label{stimadifferenza}
\end{eqnarray}
From the weak$^*$ convergence $\m_{T_n}\longrightarrow \m$, we have that 
\beqano
\left|\int f d{\Phi^L_s}_*\m_{T_n} - \int f d \m \right| \longrightarrow 0 \quad \rm{as}\; T_n\rightarrow +\infty
\eeqano
and
\beqano
\left|\int f d{\Phi^L_s}_*\m - \int f d{\Phi^L_s}_*\m_{T_n} \right| = 
\left|\int f \circ {\Phi^L_s} d\m - \int f \circ {\Phi^L_s} d\m_{T_n} \right|
 \longrightarrow 0 \quad \rm{as}\; T_n\rightarrow +\infty.
\eeqano
Moreover, using the definition of $\m_{T_n}$:
\beqano
\left|\int f d{\Phi^L_s}_*\m_{T_n} - \int f d \m_{T_n} \right| &=& 
\frac{1}{T_n} \left| \int_s^{T_n+s} f(\g(t),\dot{\g}(t))\,dt - \int_0^{T_n} f(\g(t),\dot{\g}(t))\,dt   \right| =\\
&=& \frac{1}{T_n}  
\left| \int_{T_n}^{T_n+s} f(\g(t),\dot{\g}(t))\,dt - \int_0^{s} f(\g(t),\dot{\g}(t))\,dt   \right| \leq\\
&\leq& \frac{2s}{T_n} \max_{\cE(E)}|f| \longrightarrow 0 \quad \rm{as}\; T_n\rightarrow +\infty.
\eeqano
Therefore, it follows from (\ref{stimadifferenza}) that for any $f\in \cE(E)$:
$$
\left|\int f d{\Phi^L_s}_*\m - \int f d \m \right| =0.
$$
Consequently, ${\Phi^L_s}_*\m = \m$ for any $s\in\R$, \ie $\m$ is $\Phi^L$-invariant.
\end{Proof}

To each  $\m \in \calM(L)$, we may  associate  its {\it average action}
$$ A_L(\m) = \int_{{\rm T}M} L\,d\m\,. $$
%Since $L$ is bounded from below (because of the superlinear growth condition), this 
%integral exists although it might be $+\infty$. Clearly, the measures constructed above, restricting the flow to an energy level, have finite action.
Clearly, the measures constructed in Proposition \ref{KriloffBog} have finite action.
\begin{Rem}
One can show the existence of invariant probability measures with finite action, also in the case of non-autonomous time-periodic Lagrangians. As it was originally done by Mather in \cite{Mather91}, one can apply Kryloff and Bogoliouboff's result to 
a one-point compactification of ${\rm T}M$ and consider the extended Lagrangian system that leaves the point at infinity fixed.
The main step consists in showing that the measure provided by this construction has no atomic part supported at $\infty$ 
(which is a fixed point for the extended system). \\
\end{Rem}

\noindent{\bf Note:} We shall hereafter assume that $\calM(L)$ is endowed with the {\it vague topology}, \ie the weak$^*$~topology induced by the space $C^0_{\ell}$ of continuous functions $f:\rT M \longrightarrow \R$ having at most linear growth:
$$
\sup_{(x,v)\in \rT M} \frac{|f(x,v)|}{1+\|v\|} <+\infty\,.
$$
It is not difficult to check that $\calM(L)\subset \left(C^0_{\ell}\right)^*$. Moreover it is a classical result that this space, with such a topology, is metrizable. A metric, for instance, can be defined as follows. Let $\{f_n\}_n$ be a sequence of functions with compact support in 
$C^0_{\ell}$ which is dense in the topology of uniform convergence on compact sets of $\rT M$. Define a metric on $\calM(L)$ by: %(see for instance \cite[Proposition 2-3.1]{ContrerasIturriaga})
$$
d(\mu_1,\mu_2)= \left|\int |v| d\mu_1 -  \int |v| d\mu_2 \right| + \sum_n \frac{1}{2^n \|f_n\|_{\infty}} \left| \int f_n d\mu_1 - \int f_n d\mu_2 \right|.
$$
See also Addendum 4.A at the end of this section for more details on this topology.

\begin{Prop}\label{lowsemicactionmeasures}
$A_L : \calM(L) \longrightarrow \R$ is lower semicontinuous with the vague topology on $\calM(L)$.
\end{Prop}

\begin{Proof} 
Let $A_{L,K}(\mu):= \int \min\{L, K\}\,d\mu$, for $K\in \R$. Then, $A_{L,K}$ is clearly continuous, actually Lipschitz with constant $K$. In fact:
$$
|A_{L,K}(\mu) - A_{L,K}(\nu)| \leq Kd(\m,\n).
$$
Moreover, since $A_{L,K}\uparrow A_L$  as $K\to \infty$, it follows that $A_L$ is lower semicontinuous.
\end{Proof}

\begin{Rem}
In general, this functional might not be necessarily continuous.
\end{Rem}

This Proposition has an immediate, but important, consequence.

\begin{Cor}\label{existminmeas}
There exists $\m\in\calM(L)$, which minimizes $A_L$ over $\calM(L)$.
\end{Cor}

\noindent We shall call a measure $\mu\in\calM(L)$, such that $A_L(\mu)=\min_{\calM(L)}A_L$, an {\it action-minimizing measure} of $L$ (compare with Remark \ref{Rem3.5} ({\it ii})).\\

Actually, one can find many other ``interesting'' measures, besides those found by minimizing $L$. We shall see later   (see Remark \ref{remarkinclusionsmathersets}), that in fact these measures minimizing the action of $L$, correspond to special measures with ``trivial'' homology (or rotation vector).
On the other hand, as we have already pointed out in Proposition \ref{propminimalactiontorus} (and the following Warning), in order to get information on a specific KAM torus and characterize it via the invariant probability measures  supported on it, one needs necessarily to modify the Lagrangian.  Think for instance of an integrable system, for which we have a whole family of KAM tori. It would be unreasonable to expect to obtain ALL of them, just by minimizing a single Lagrangian action. One needs somehow to introduce some sort of ``weight'', which, without modifying the dynamics of the system (\ie the Euler Lagrange flow), allows one to magnify certain motions rather than others.  What we found out in Proposition \ref{propminimalactiontorus}, was that this can be easily achieved by subtracting to our Lagrangian a linear function $c\cdot v$, where $c$ represented the {\it cohomology class} of the invariant Lagrangian torus we were interested in, or, in other words, the cohomology class of its graph (that is a closed $1$-form) in $\rT^*\T^d$.

 A similar idea can be implemented for a general Tonelli Lagrangian.
Observe, in fact,  that if $\eta$ is a $1$-form on $M$, we can interpret it as a function 
on the tangent space (linear on each fiber)
\beqano 
\hat{\eta}: {\rm T}M &\longrightarrow&  \R \\
(x,v) &\longmapsto& \langle \eta(x),\, v\rangle_x
\eeqano
and consider a new Tonelli Lagrangian $L_{\eta}:= L - \hat{\eta}$. The 
associated Hamiltonian will be given by
$H_{\eta}(x,p) = H(x,\eta(x) + p)$ (it is sufficient to write down the Fenchel-Legendre transform of $L$).

\begin{Lem}\label{lemmamodifiedLagrangians}
If $\eta$ is closed, then $L$ and $L_{\eta}$ 
have the same Euler-Lagrange flow on ${\rm T} M$. 
\end{Lem}

\begin{Proof}
Observe  that, since $\eta$ is closed, the variational equations $\d \!\left[\int L\, dt \right]=0$ and $\d \!\left[\int 
(L-\hat{\eta})\, dt \right]=0$ have the same extremals for the fixed end-point problem and consequently $L$ and $L_{\eta}$ 
have the same Euler-Lagrange flows.\\
One can also check this directly, writing down the Euler-Lagrange equations associated to $L$ and $L_{\eta}$ and observe that,
since $\eta$ is closed (\ie $d\eta=0$), then
\beqano
\frac{d}{dt}\left(\frac{\dpr} {\dpr v} \hat{\eta}(\g,\dot\g)\right) = \frac{d}{dt}\left({\eta}(\g)\right) = \frac{\dpr}{\partial x}\left({\hat{\eta}}(\g,\dot\g)\right).
\eeqano
Let us show in fact that they are equal component by component.  If we denote 
by $\eta_j$  and $\gamma_j$ the $j$-th components  of $\eta$ and $\gamma$, then the $i$th component of the left-hand side of the above equation becomes:
$$
\sum_{j=1}^d \frac{\partial \eta_i}{\partial x_j }(\gamma) \cdot \dot{\gamma}_j
$$
and on the other side:
$$
\sum_{j=1}^d \frac{\partial \eta_j}{\partial x_i }(\gamma) \cdot \dot{\gamma}_j\,.
$$
Since $\eta$ is closed, then $\frac{\partial \eta_j}{\partial x_i } = \frac{\partial \eta_i}{\partial x_j }$ for each $i$ and $j$. Therefore, the two expressions are the same.
Therefore, the term coming from $\hat{\eta}$ does not give any contribution to the equations.
\end{Proof}

Although the extremals of these two variational problems are the same, this is not generally true for the orbits
``minimizing the action'' (we shall give a precise definition of  ``minimizers'' later in section \ref{sec1.3}). We have already given evidence of  this  in Proposition \ref{propminimalactiontorus}.
What one can say is that these action-minimizing objects stay the same when we change the 
Lagrangian by an exact $1$-form.

\begin{Prop}\label{closedmeasure}
If $\m \in \calM(L)$ and $\eta=df$ is an exact $1$-form, then 
$\int{\widehat{df}} d\m =0$.
\end{Prop}

\noindent The proof is essentially the same as in Section \ref{cartoonexample} (Property 2), so we omit it. See also \cite[Lemma on page 176]{Mather91}. Let us point out that a measure (not necessarily invariant) that satisfies this condition, is called {\it closed measure} (see also the Addendum at the end of this section).

Thus, for a fixed $L$, the minimizing measures will 
depend only on the de Rham cohomology class $c=[\eta] \in \rH^1(M;\R)$. 
Therefore, instead of studying the action minimizing properties of a single Lagrangian, one can consider a family of such ``modified'' Lagrangians, parameterized over $H^1(M;\R)$. \\
Hereafter, for any given $c\in H^1(M;\R)$, we shall denote by $\eta_c$ a closed $1$-form with that cohomology class.\\

\begin{Def}
Let $\eta_c$ be a closed $1$-form of cohomology class $c$. Then, if $\m \in \calM(L)$ minimizes 
$A_{L_{\eta_c}}$ over $\calM(L)$, we shall say that $\m$ is a {\it $c$-action minimizing measure} (or $c$-{\it minimal measure}, or {\it Mather's measure} with cohomology $c$).
\end{Def}

\begin{Rem}
Observe that the {\it cohomology class} of an action-minimizing invariant probability measure,  is not something intrinsic in the measure itself nor in the dynamics, but it depends on the specific choice of the  Lagrangian $L$. Changing the Lagrangian $L \longmapsto L-\eta$ by a closed $1$-form $\eta$, we shall change all the cohomology classes of its action minimizing measures by $-[\eta] \in \rH^1(M;\R)$. Compare also with Remark \ref{Remrotationvector} ({\it ii}).
\end{Rem}

One can consider the function on $\rH^1(M;\R)$, which associates to each cohomology class $c$, {\it minus} the value of the corresponding minimal action of the modified Lagrangian $L_{\eta_c}$ (the {\it minus} sign is introduced for a convention that might probably become clearer later on):
\beqa{defalfa}
\a: {\rm H}^1(M;\R) &\longrightarrow& \R \nn\\
c &\longmapsto& - \min_{\m\in\calM(L)} A_{L_{\eta_c}}(\m)\,.
\eeqa
This function $\alpha$ is well-defined (it does not depend on the choice of the representatives of the cohomology classes) and it is easy to see that it is convex. This is generally known as {\it Mather's $\alpha$-function}. As it happened for KAM tori (see Remark \ref{Rem3.5}), we shall see in section \ref{sec1.3} that the value $\a(c)$ is related to the energy level containing such $c$-action minimizing measures (they coincide). See also \cite{Carneiro}. Therefore, if we have an integrable Tonelli Hamiltonian $H(x,p)=h(p)$ it is easy to deduce that $\a(c)=h(c)$. For this and several other reasons that we shall see later on, this function is sometime called {\it effective Hamiltonian}.
It is interesting to remark that the value of this $\a$~-~function coincides 
with what is called {\it Ma\~n\'e's critical value}, which will be introduced 
later in sections \ref{sec1.3} and \ref{sec1.4}.\\

We shall denote by $\calM_c(L)$ the subset of $c$-action minimizing measures:
$$\calM_c := \calM_c(L)= \{\m \in \calM(L): \; A_{L}(\m)<+\infty \;{\rm and}\; A_{L_{\eta_c}}(\m)
=-\a (c)\}.$$

We can now define a first important family  of invariant sets (see Remark \ref{Rem3.5} ({\it iii})): 
the {\it Mather sets}. 

\begin{Def}
For a cohomology class $c \in {\rm H}^1(M;\R)$, we define 
the {\it Mather set of cohomology class} $c$  as:
\be \widetilde{\cM}_c := {\bigcup_{\m \in \calM_c} {\rm supp}\,\m} 
\subset {\rm T}M\,.\label{2.3}\ee
The projection on the base manifold $\cM_c = \pi \left(\widetilde{\cM}_c\right)
\subseteq M$ is called {\it projected Mather set} (with cohomology class $c$).
\end{Def}

\begin{Rem}\label{remarknoclosure}
({\it i}) This set is clearly non-empty and invariant. Moreover, it is also closed.
Observe that usually it is defined as the {\it closure} of the set on the right-hand side. This is indeed
the original definition given by Mather in \cite{Mather91}. However, it is not difficult to check that this set is
already closed. In fact, since the space of 
probability measures on ${\rm T}M$ is a separable metric space, one can take
a countable dense set $\{\m_n\}_{n=1}^{\infty}$ of Mather's 
measures and consider the new measure
$\tilde{\m} = \sum_{n=1}^\infty \frac{1}{2^n}\m_n$. This is still an invariant 
probability measure  and it is $c$-action minimizing (it follows from the convexity of $\a$), hence
${\rm supp}\,\tilde{\m} \subseteq \widetilde{\cM}_c$. But, from the definition of $\tilde{\m}$ one can clearly deduce that ${\rm supp}\,\tilde{\m} = 
\overline{{\bigcup_{\m \in \calM_c}}{\rm supp}\,\m} \supseteq \widetilde{\cM}_c$. This implies that ${\rm supp}\,\tilde{\m}=\widetilde{\cM}_c$.
Therefore, as the support of a single action minimizing 
measure, $\widetilde{\cM}_c $ is closed (the support of a measure, by definition, is closed).
Moreover, this remark points out that there always exists a Mather's measure
$\m_c$ of full support, \ie 
${\rm supp}\,{\m_c} = 
\widetilde{\cM}_c$.\\
\noindent({\it ii}) As we have pointed out in Remark \ref{Rem3.5}, if there is a KAM torus ${\mathcal T}$ of cohomology class $c$, then $\widetilde{\cM}_c = \cL^{-1}({\mathcal T})$. In particular, the same proof continues to hold if we replace ${\mathcal T}$ with any invariant Lagrangian graph $\L$ of cohomology class $c$, which supports an invariant measure $\m$ of full support  (\ie ${\rm supp}\,\m = \L$).\\
\end{Rem}

In \cite{Mather91} Mather proved the celebrated {\it graph theorem}:

\begin{Teo}[{\bf Mather's graph theorem}]\label{Theograph} Let $\widetilde{\cM}_c$ be defined as in 
(\ref{2.3}). 
The set $\widetilde{\cM}_c$ is compact, invariant 
under the Euler-Lagrange flow and $\pi|{\widetilde{\cM}_c}$ is an injective 
mapping of $\widetilde{\cM}_c$ into $M$, and its inverse $\pi^{-1}: \cM_c 
\longrightarrow \widetilde{\cM}_c$ is 
Lipschitz. \end{Teo}

We shall not prove this theorem here, but it will be deduced from a more general result in Section \ref{sec1.3}, namely, the graph property of the {\it Aubry set}.
Moreover, similarly to what we have seen for KAM tori (see Remark \ref{Rem3.5}), in the autonomous case this set is contained in a well-defined energy level, which  can be characterized in terms of the minimal action $\a(c)$.

\begin{Teo}[{\bf Carneiro, \cite{Carneiro}}]\label{teocarneiro}
The Mather set  $\widetilde{\cM}_c$ is contained in the energy level $\{H \circ \cL(x,v) = \a(c)\}$.
%\label{2.4}\ee
\end{Teo}

\noindent Also this result will be deduced from a similar result for the Aubry and Ma\~n\'e sets. See Proposition \ref{teocarneiropermane} in Section \ref{sec1.3}.\\

%\begin{Rem}{\rm The last statement, corresponding to (\ref{2.4}), is due to 
%Dias-Carneiro \cite{Carneiro} and holds only in the autonomous case.}
%\end{Rem}

\vspace{5 pt}
Now, we would like to shift our attention to a related problem. As we have seen in section \ref{cartoonexample}, instead of considering different minimizing problems over $\calM(L)$, obtained by modifying the Lagrangian $L$, one can alternatively try to minimize the Lagrangian $L$ putting some ``{\it constraints}'', such as, for instance, fixing the {\it rotation vector} of the measures. In order to generalize this to Tonelli Lagrangians on compact manifolds, we first need to define what we mean by rotation vector of an invariant measure.

Let $\mu\in \calM(L)$. Thanks to the superlinearity of $L$, 
the integral  $ \int_{{\rm T}M} \hat{\eta} d\m$ 
is well defined and finite for any 
closed 1-form $\eta$ on $M$.
Moreover, we have proved in Proposition \ref{closedmeasure} that if $\eta$ is exact, then such an integral is zero, \ie 
$ \int_{{\rm T}M} \hat{\eta} d\m=0$.
Therefore, one can define a linear functional: 
\beqano
{\rm H}^1(M;\R) &\longrightarrow& \R \\
c &\longmapsto& \int_{{\rm T}M} \hat{\eta} d\m\,,
\eeqano
where $\eta$ is any closed $1$-form on $M$ with cohomology class $c$. By 
duality, there 
exists $\rho (\m) \in {\rm H}_1(M;\R)$ such that
$$
\int_{{\rm T}M} \hat{\eta} \,d\m = \langle c,\rho(\m) \rangle
\qquad \forall\,c\in {\rm H}^1(M;\R)$$ 
(the bracket on the right--hand side denotes the canonical pairing between 
cohomology and 
homology). We call $\rho(\m)$ the {\it rotation vector} of $\m$ (compare with the definition given in Section \ref{cartoonexample} Property (1)). This rotation vector is 
the same as the Schwartzman's asymptotic cycle of $\mu$ (see \cite{Schwartzman} for more details).

\begin{Rem}\label{Remrotationvector}
({\it i}) It is possible to provide a more ``geometrical'' interpretation of this. Suppose for the moment that $\mu$ is ergodic. Then, it is known that a generic orbit $\g(t):=\pi \Phi^L_t(x,v)$, where $\pi:\rT M \longrightarrow M$ denotes the canonical projection, will return infinitely many often close (as close as we like) to its initial point $\g(0)=x$. We can therefore consider a sequence of times $T_n \to +\infty$ such that $d(\g(T_n),x)\to 0$ as $n\to +\infty$, and consider the closed loops $\s_n$ obtained by ``closing'' $\g|[0,T_n]$ with the shortest geodesic connecting $\g(T_n)$ to $x$. Denoting by $[\s_n]$ the homology class of this loop, one can verify \cite{Schwartzman}  that $\lim_{n\to\infty} \frac{[\s_n]}{T_n} = \rho(\mu)$, independently of the chosen sequence $\{T_n\}_n$. In other words, in the case of ergodic measures, the rotation vector tells us how on average a generic orbit winds around $\rT M$. If $\mu$ is not ergodic, $\rho(\mu)$ loses this neat geometric meaning, yet it may be interpreted as the average of  the rotation vectors of its different ergodic components. 

({\it ii}) It is clear from the discussion above that the rotation vector of an invariant measure depends only on the dynamics of the system  (\ie the Euler-Lagrange flow) and not on the chosen Lagrangian. Therefore, it does not change when we modify our Lagrangian adding a closed one form.
\end{Rem}

A natural question is whether or not for a given Tonelli Lagrangian $L$, there exist invariant probability measures for any given rotation vector. The answer turns out to be affirmative. \\

In fact,
using that the action functional $A_L: \calM(L) 
\longrightarrow \R$ is lower semicontinuous, one can prove the following \cite{Mather91}:

\begin{Prop}
{\rm(i)} The map $\rho: \calM(L) \longrightarrow \rH_1(M;\R)$ is continuous.\\
{\rm(ii)} The map $\rho: \calM(L) \longrightarrow \rH_1(M;\R)$ is affine, \ie for any $\mu,\nu \in \calM(L)$ and $a,b\geq 0$ with $a+b=1$,  $\rho(a\mu+b\nu) = a\rho(\mu)+b\rho(\nu)$.\\
{\rm(iii)}For every 
$h\in {\rm H}_1(M;\R)$ 
there exists $\m\in \calM(L)$ with  $A_L(\m) < \infty$ and $\rho(\m)=h$. In other words,
the map $\rho$ is surjective.
\end{Prop}

\begin{Proof}
(i) Let us fix a basis $h_1,\ldots, h_b$ in $\rH_1(M;\R)$, where $b$ is the first Betti number of $M$, \ie $b=\dim \rH_1(M;\R)$. If 
$\mu_n\to \mu$ and $\eta_c$ is any closed one-form on $M$ of cohomolgy class $c$, then
$$
\langle c, \rho(\mu_n) \rangle = \int \eta_c \cdot v d\mu_n \stackrel{n\rightarrow +\infty}{\longrightarrow} \int \eta_c \cdot v d\mu =  \langle c, \rho(\mu) \rangle,
$$
that is equivalent to say that:
$$
\langle c, \rho(\mu_n)-\rho(\mu) \rangle \stackrel{n\rightarrow +\infty}{\longrightarrow} 0 \quad \forall\,c \in \R^b \simeq H^1(M;\R)\,.
$$
Therefore, $\rho(\mu_n)-\rho(\mu) {\longrightarrow} 0$ as $n\to +\infty$.

(ii) The fact that the map $\rho$ is affine, is a trivial consequence of the definition of rotation vector.

(iii) Let $h \in \rH_1(M;\Z)$ be an integer homology class and choose $\gamma_h:[0,1]\longrightarrow M$ a closed loop with homology $h$. Let us consider the loop space:
$$
\Omega_h:= \{\xi: [0,1] \longrightarrow M,\; \mbox{abs. cont. loop with the same free homotopy class as }\; \g_h\}.
$$ 
One can prove \cite{ManeI}  that the Lagrangian action functional $A_L$ has a minimum on this space, which is a periodic orbit of the Euler-Lagrange flow, with period $1$ (it is essentially Tonelli's theorem). Let us consider this periodic orbit $\gamma:[0,1] \longrightarrow M$ and define an invariant probability measure $\mu_{\gamma}$ evenly distributed along this periodic orbit:
$$
\int f \, d\mu_{\gamma} := \int_0^1 f(\g(t),\dot{\g}(t))\,dt \qquad \forall f\in C^0(\rT M).
$$
It is easy to verify that such a measure is invariant, \ie ${\Phi^L_t}^*\mu_{\gamma} = \mu_{\gamma}$ and, using the definition of rotation vector and Remark \ref{Remrotationvector}, that $\rho(\mu_{\gamma})=h$.
Since the map $\rho$ is affine and $\calM(L)$ is convex, it follows that $\rho(\calM(L))$ is convex. It follows from what discussed above, that this set must contain the convex hull of $\rH_1(M;\Z)$ and therefore $\rho(\calM(L))= \rH_1(M;\R)$. 
\end{Proof}

As already pointed out in Section \ref{cartoonexample}, amongst all probability measures with a prescribed rotation vector, a peculiar role - from a dynamical systems point of view - will be played by those minimizing the average action.
Following Mather, let us consider the minimal value of the average action $A_L$ over the 
probability measures with rotation vector $h$. Observe that this minimum is actually achieved because of the lower semicontinuity of $A_L$ and the compactness of $\rho^{-1}(h)$ ($\rho$ is continuous and $L$ superlinear). Let us define
\bea
\b: {\rm H}_1(M;\R) &\longrightarrow& \R \nn\\
h &\longmapsto& \min_{\m\in\calM(L):\,\rho(\m)=h} A_L(\m)\,.\label{2.2}
\eea
This function $\beta$ is what is generally known as {\it Mather's 
$\beta$-function} and it is immediate to check that it is convex. 
As we have noticed in Remark \ref{Remark3.7}, if there is a KAM torus of cohomology class $c$ and rotation vector $\rho$, then  $\beta(\rho)= -E_c +c\cdot \rho$. Therefore, if we have an integrable Tonelli Hamiltonian $H(x,p)=h(p)$ and the associated Lagrangian $L(x,v)=\ell(v)$, it is easy to deduce that $\beta(h)=\ell(h)$. For this and several other reasons that we shall see later on, this function is sometime called {\it effective Lagrangian}.
Moreover, this function is also related to the notion of {\it stable norm} for a metric $d$ (see for instance \cite{Massartstablenorm}).\\

We can now define what we mean by action minimizing measure with a given rotation vector (compare with Remark \ref{Remark3.7}({\it ii})).

\begin{Def}
A measure $\m \in \calM(L)$ realizing the minimum in (\ref{2.2}), \ie such that $A_L(\m)~=~\b(\rho(\m))$, is called an {\it action minimizing} (or {\it minimal} or {\it Mather's}) {\it measure} with rotation vector $\rho(\m)$.
\end{Def}

\begin{Rem}
We shall see in Section \ref{sec1.3} that, differently from what happens with invariant probability measures,
it will not be always possible to find {\it action-minimizing orbits} for any given rotation vector (not even define a rotation vector for each action minimizing orbit). This is one of the main difference with the twist map case. In higher dimensions, in fact, an example due to Hedlund \cite{Hedlund} provides the existence of a Riemannian metric on a three-dimensional torus, for which minimal geodesics 
exist only in three directions.\\
 \end{Rem}

 We shall denote by $ \calM^h(L)$ the subset of action minimizing measures with rotation vector $h$:
 $$ \calM^h := \calM^h(L) = \{\m \in \calM(L): \; A_L(\m)<+\infty, \;\rho(\m)= h \; 
{\rm and}\;  A_L(\m)=\beta(h)\}.$$

 This allows us to define another important familty of invariant sets (see also Remark \ref{Remark3.7} ({\it iii})).

 \begin{Def}
 For a homology class (or rotation vector) $h\in{\rm H}_1(M;\R)$, we define the
{\it Mather set corresponding to a rotation vector} $h$ as
\be \widetilde{\cM}^h := {\bigcup_{\m \in \calM^h} {\rm supp}\,\m} 
\subset {\rm T}M\,,\label{2.5}\ee
and the projected one as $\cM^h = \pi \left(\widetilde{\cM}^h\right) 
\subseteq M$. 
 \end{Def}

 \begin{Rem}\label{Remnoclosurerotation}
({\it i}) Similarly to what we have pointed out in Remark \ref{remarknoclosure}, 
this set is also non-empty and invariant. Moreover, it is also closed. Also in this case, it is not necessary to consider the closure of this set - as it is usually done in the literature- since it is already closed. Just  take
a countable dense set $\{\m_n\}_{n=1}^{\infty}$ of Mather's 
measures with rotation vector $h$ and consider the new measure
$\tilde{\m} = \sum_{n=1}^\infty \frac{1}{2^n}\m_n$. This is still an invariant 
probability measure and its rotation vector is $h$.
Moreover, it follows from the convexity of $\beta$ that this measure is action-minimizing among all measures with rotation vector $h$. Hence,  ${\rm supp}\,\tilde{\m} \subseteq \widetilde{\cM}^h$.
Clearly ${\rm supp}\,\tilde{\m} = 
\overline{{\bigcup_{\m \in \calM^h} {\rm supp}\,\m}}
\supseteq \widetilde{\cM}^h$. 
This implies that ${\rm supp}\,\tilde{\m}=\widetilde{\cM}^h$.
Therefore, $\widetilde{\cM}^h $ is closed.
Moreover, this shows that there always exists a Mather's measure $\m^h$ of full support, \ie 
${\rm supp}\,{\m^h} = 
\widetilde{\cM}^h$.\\
({\it ii}) As we have pointed out in Remark \ref{Remark3.7}, if there is a KAM torus ${\mathcal T}$ of rotation vector $h$, then $\widetilde{\cM}^h = \cL^{-1}({\mathcal T})$. In particular, the same proof continues to hold if we replace ${\mathcal T}$ with any invariant Lagrangian graph $\L$, which supports an invariant measure $\m$ of rotation vector $\rho$ and of full support  (\ie ${\rm supp}\,\m = \L$).\\

\end{Rem}

Also for the Mather set corresponding to a rotation vector, we can prove a result  similar to Theorem \ref{Theograph}.

\begin{Teo}\label{matherrotation} 
Let $\widetilde{\cM}^h$ be defined as in (\ref{2.5}). $\widetilde{\cM}^h$
is compact, invariant under the Euler-Lagrange flow and 
$\pi|{\widetilde{\cM}^h}$ is an injective 
mapping of $\widetilde{\cM}^h$ into $M$ and its inverse $\pi^{-1}: \cM^h \longrightarrow \widetilde{\cM}^h$ is Lipschitz.
\end{Teo}

\begin{Rem} Although the graph property for $\widetilde{\cM}^h$ is not proved in 
\cite{Mather91}, it is easy to deduce it from Theorem \ref{Theograph}, using the fact that  $\widetilde{\cM}^h $ can be seen as the support of a single action-minimizing measure (Remark \ref{Remnoclosurerotation}) and that this set is included in some 
$\widetilde{\cM}_c$, for some suitable $c\in \rH^1(M;\R)$ (Proposition \ref{CaracMinim}).\\
\end{Rem}

The above discussion leads to two equivalent formulations for the minimality 
of a measure $\m$:
\begin{itemize}
\item there exists a homology class $h \in {\rm H}_1(M;\R)$, namely its 
rotation vector $\rho(\m)$, such that $\m$ minimizes $A_L$ amongst all 
measures in $\calM(L)$ with rotation vector $h$; \ie $A_L(\m)=\b (h)$.
\item There exists a cohomology class $c \in {\rm H}^1(M;\R)$, such that $\m$ minimizes $A_{L_{\eta_c}}$ 
amongst all probability measures in $\calM(L)$; \ie $A_{L_{\eta_c}}(\m)=-\a (c)$.\\
\end{itemize}

\noindent What is the relation between two these different approaches? Are they equivalent, \ie 
$\bigcup_{h \in {\rm H}_1(M;\R)} \calM^h = \bigcup_{c \in {\rm H}^1(M;\R)} 
\calM_c\,$ ?\\

In order to comprehend the relation between these two families of action-minimizing measures, we need to understand better the properties of these functions
$$\a: \rH^1(M;\R) \longrightarrow \R \quad {\rm and}\quad  \b:\rH_1(M;\R) \longrightarrow \R$$ that we have introduced above. Let us start with the following trivial remark.

\begin{Rem}\label{relalphabetainteg} As we have previously pointed out, if we have an integrable Tonelli Hamiltonian $H(x,p)=h(p)$ and the associated Lagrangian $L(x,v)=\ell(v)$, then $\a(c)=h(c)$ and $\beta(h)=\ell(h)$. In this case, the cotangent space $\rT^*\T^d$ is foliated by invariant tori ${\mathcal T}^*_c:=\T^d\times\{c\}$ and the tangent space $\rT\T^d$ by invariant tori $\widetilde{\mathcal T}^h:=\T^d\times\{h\}$. In particular, we proved that 
$$
\widetilde{\cM}_c = \cL^{-1}({\mathcal T}_c) = \widetilde{\mathcal T}^h =  \widetilde{\cM}^h,
$$
where $h$ and $c$ are such that  $h=\nabla h(c) = \nabla \a(c)$ and
$c= \nabla \ell(h) = \nabla \beta(h)$.
\end{Rem}

 We would like to prove that such a relation that links Mather sets of a certain cohomology class to Mather sets with a given rotation vector, goes beyond the specificity of this situation. Of course, one main difficulty is that in general the {\it effective Hamiltonian} $\a$ and the {\it effective Lagrangian} $\b$, although being convex and superlinear (see Proposition \ref{alphabetasuperlinear}), are not necessarily differentiable. \\
Before stating and proving the main relation between these two functions, let us recall some definitions and results from classical convex analysis (see \cite{Rockafellar}). Given a convex function $\varphi: V \longrightarrow \R\cup \{+\infty\}$ on a finite dimensional vector space $V$, one can consider a {\it dual} (or {\it conjugate}) function  defined on the dual space $V^*$, via the so-called {\it Fenchel transform}: $\f^*(p):= \sup_{v\in V} \big(p\cdot v - \f(v)\big)$.

\begin{Prop}\label{alphabetasuperlinear}
$\a$ and $\b$ are convex conjugate, \ie
$\a^* = \b$ and $\b^* = \a$. In particular, it follows that $\a$ and $\b$ have superlinear growth.
\end{Prop}

\begin{Proof}
First of all, recall that $\left(\rH_1(M;\R)\right)^*\simeq \rH^1(M;\R)$ and $\left(\rH^1(M;\R)\right)^*\simeq \rH_1(M;\R)$. Let us compute $\beta^*$:
\beqano
\beta^*(c) &=& \max_{h\in \rH_1(M;\R)} \left(\langle c,h\rangle - \beta(h)\right) =\\
&=& - \min_{h\in \rH_1(M;\R)}  \left( \beta(h) - \langle c,h\rangle \right) = \\
&=& - \min_{h\in \rH_1(M;\R)}  \left( \min_{\m \in \calM^{h}(L)} A_L(\mu) - \langle c,h\rangle \right) = \\
&=& - \min_{h\in \rH_1(M;\R)}  \min_{\m \in \calM^{h}(L)} (A_L(\mu) - \langle c,\rho(\mu)\rangle) = \\
&=& - \min_{h\in \rH_1(M;\R)}  \min_{\m \in \calM^{h}(L)} A_{L_{\eta_c}}(\mu) = \\
&=&- \min_{\m \in \calM (L)} A_{L_{\eta_c}}(\mu) = \\
&=& -\a(c)\,.
\eeqano
Similarly, one can check that $\a^*=\beta$:
\beqano
\a^*(h) &=& \max_{c\in \rH^1(M;\R)} \left(\langle c,h\rangle - \a(c)\right) =\\
&=& \max_{c\in \rH^1(M;\R)} \left(\langle c,h\rangle + \min_{\mu\in\calM(L)} A_{L_{\eta_c}}(\mu)\right) =\\
&=&  \max_{c\in \rH^1(M;\R)} \min_{\mu\in\calM(L)} \left(A_{L}(\mu) + \langle c,h-\rho(\mu)\rangle\right) =\\ 
&=&  \min_{\mu\in\calM(L)} \max_{c\in \rH^1(M;\R)}  \left(A_{L}(\mu) + \langle c,h-\rho(\mu)\rangle\right) =\\ 
&=&  \min_{\rho\in \rH_1(M;\R)} \min_{\mu\in\calM^{\rho}(L)} \max_{c\in \rH^1(M;\R)}  \left(A_{L}(\mu) + \langle c,h-\rho\rangle\right)\,,\\ 
\eeqano
where in the second last line we could exchange the order of the max and the min, using  
a general result by Rockafellar, that requires concavity in one variable, convexity in the other one, and some compactness assumption (see \cite[Section 36]{Rockafellar}).
Observe now that if $h\neq \rho$, then $\max_{c\in \rH^1(M;\R)}  \left(A_{L}(\mu) + \langle c,h-\rho\rangle\right) = +\infty$. Therefore: 
\beqano
\a^*(h) &=& \ldots \; = \min_{\rho\in \rH_1(M;\R)} \min_{\mu\in\calM^{\rho}(L)} \max_{c\in \rH^1(M;\R)}  \left(A_{L}(\mu) + \langle c,h-\rho\rangle\right)= \\
&=&  \min_{\mu\in\calM^{h}(L)} A_{L}(\mu) = \beta(h)\,.
\eeqano
The second statement of this proposition follows from a general property of convex conjugation. Let $\varphi: V\rightarrow \R$ be a convex function on a finite dimensional vector space $V$ and let $\f^*: V^* \longrightarrow \R\cup\{+\infty\}$ be its convex conjugate.
Then (see \cite{Rockafellar}): $\f^*$ is finite everywhere if and only if $\f$ has superlinear growth, \ie $\frac{\f(x)}{\|x\|} \longrightarrow +\infty$ as $\|x\|\to +\infty$.\\
\end{Proof}

 Next proposition  will allow us to clearify the relation (and duality) between the two minimizing procedures above. To state it, recall that, like any convex 
function on a
finite-dimensional space, $\b$
admits a subderivative at each point $h\in \rH_1(M;\R)$, \ie we can find $c\in 
\rH^1(M;\R)$ such that
\beqa{subdiff}
\forall h'\in \rH_1(M;\R), \quad \b(h')-\b(h)\geq \langle c,h'-h\rangle.\eeqa
As it is usually done, we shall denote by $\partial \b(h)$ the set of $c\in 
\rH^1(M;\R)$ that 
are subderivatives of $\b$ at $h$, \ie the set of $c$ which satisfy the 
inequality above.  Similarly, we shall denote by $\dpr \a(c)$ the set of subderivatives of $\a$ at $c$.\\
Fenchel's duality implies an easier characterization of subdifferentials.

\begin{Prop}
$c\in \partial \b(h)$ if and only if $\langle c,h\rangle=\a(c)+\b(h).$ 
Similarly $h\in \partial \a(c)$ if and only if $\langle c,h\rangle=\a(c)+\b(h).$ 
In particular, $c\in \partial \b(h)$ if and only if $h\in \partial \a(c)$.
\end{Prop}

\begin{Proof}
We shall prove only the first statement. The second one is analogous and the third one a trivial consequence.\\
$[\Longleftarrow]$ Using that $\beta(h') \geq \langle c,h' \rangle -\a(c)$ for each $h'\in \rH_1(M;\R)$, it follows that for each $h'\in \rH_1(M;\R)$:
\beqano
\beta(h')-\beta(h) &\geq&  \langle c,h' \rangle -\a(c) - \beta(h) = \langle c,h' \rangle - (\a(c)+ \beta(h)) =\\
&=& \langle c,h' \rangle - \langle c,h \rangle = \langle c,h'-h \rangle.
\eeqano
$[\Longrightarrow]$ It follows from Fenchel-Legendre inequality that 
$\a(c) + \beta(h') \geq \langle c,h' \rangle$ for each $h'\in \rH_1(M;\R)$. Therefore we only need to prove the reverse inequality. In fact, using (\ref{subdiff}) one can deduce that for each $h'\in \rH_1(M;\R)$:
\beqano
\a(c) +\beta(h) &\leq& \a(c) + \beta(h') - \langle c,h'-h \rangle  = 
\a(c) + \beta(h') - \langle c,h' \rangle + \langle c,h \rangle =\\
&=& \Big( \a(c) + \beta(h') - \langle c,h' \rangle\Big) + \langle c,h \rangle\,.
\eeqano
Therefore, taking the minimum over $h'$ on the left-hand side we obtain:
$$
\a(c) +\beta(h) \leq \min_{h'\in\rH_1(M;\R)} \Big( \a(c) + \beta(h') - \langle c,h' \rangle\Big) + \langle c,h \rangle = \langle c,h \rangle,
$$
where in the last equality we used that $\min_{h'\in\rH_1(M;\R)} \big( \a(c) + \beta(h') - \langle c,h' \rangle\big)=0$, as it follows easily from Fenchel-Legendre duality.
\end{Proof}

We can now prove that what observed in Remark \ref{relalphabetainteg} continues to hold in the general case

\begin{Prop}\label{CaracMinim}
Let  $\m \in \calM(L)$ be an invariant probability measure. Then:\\
{\rm (i)} $A_L(\m)=\b(\rho(\m))$ {if and only if}  there exists $c\in {\rm H}^1(M;\R)$ such that $\m$ minimizes $A_{L_{\eta_c}}$  
{\rm(}\ie $A_{L_{\eta_c}}(\mu)=-\a(c)${\rm)}.\\
{\rm (ii)} If $\m$ satisfies $A_L(\m)=\b(\rho(\m))$ and $c \in 
\rH^1(M;\R)$, then $\m$ minimizes $A_{L_{\eta_c}}$ if and only if $c\in \partial 
\b(\rho(\mu))$ {\rm(}or equivalently $\langle c,h\rangle=\a(c)+\b(\rho(\m)\rm{)}$.
\end{Prop}

\begin{Proof} We shall prove both statements at the same time.
Assume $A_L(\m_0)=\b(\rho(\m_0))$. Let $c\in \partial \b(\rho(\mu_0))$, by 
Fenchel's duality this is equivalent to
\begin{align*}
\a(c)&=\langle c,\rho(\m_0)\rangle-\b(\rho(\m_0))\\
&=\langle c,\rho(\m_0)\rangle-A_L(\m_0)\\
&=-A_{L_{\eta_c}}(\m_0).
\end{align*}
Therefore $A_{L_{\eta_c}}(\m_0)=\min_{\m \in \calM(L)} 
A_{L_{\eta_c}}(\m)$.

Assume conversely that $A_{L_{\eta_c}}(\m_0)=\min_{\m \in \calM(L)} A_{L_{\eta_c}}(\m)$, for some 
given cohomology class $c$. Then, it follows  that
$$\a(c)=-A_{L_{\eta_c}}(\m_0),$$
which can be written as 
$$ \langle c,\rho(\m_0)\rangle=\a(c)+A_L(\m_0).$$
It now suffices to use the Fenchel inequality $ \langle c,\rho(\m_0)\rangle\leq 
\a(c)+\b(\rho(\m_0))$, and the inequality $\b(\rho(\m_0))\leq A_L(\m_0)$, given 
by the definition of $\b$, to obtain the  equality 
$$ \langle c,\rho(\m_0)\rangle=\a(c)+\b(\rho(\m_0)).$$
In particular, we have $A_L(\m_0)=\b(\rho(\m_0))$.\\
 \end{Proof}

%\begin{Cor}\label{corolinclusions}
%$\widetilde{\cM^h} \cap \widetilde{\cM_c} \neq 0$ if and only if $\widetilde{\cM^h} \subseteq \widetilde{\cM_c}$.
%\end{Cor}

%\begin{Proof}
%One inclusion ($\Longleftarrow$) is obvious. Let us prove the other one.
%[$\Longrightarrow$] Suppose that $\widetilde{\cM^h} \cap \widetilde{\cM_c} \neq 0$. This intersection is invariant under the Euler-Lagrange flow (it is the intersection of two invariant sets) and compact (it is the intersection of two compact sets). Therefore, there exists an action minimizing measure with rotation vector $h$, whose support is contained in $\widetilde{\cM}_c$ and hence, it is also $c$-action minimizing.

%{\bf FINIRE QUESTA DIMOSTRAZIONE!!!!!!!}\\
%Therefore, from Kryloff-Bogoliouboff theorem (see Proposition \ref{KriloffBog}) it follows that there exists
%\end{Proof}

\begin{Rem}\label{remarkinclusionsmathersets}
({\it i}) It follows from the above proposition, that both minimizing procedures lead to the same sets of invariant probability measures:
$$
\bigcup_{h \in {\rm H}_1(M;\R)} \calM^h = \bigcup_{c \in {\rm H}^1(M;\R)} 
\calM_c\,.
$$
In other words,  minimizing over the set of invariant measures with a fixed rotation vector or minimizing - globally - the modified Lagrangian (corresponding to a certain cohomology class) are dual problems,  as the ones that often appears in linear programming and optimization.

({\it ii}) In particular, we have the following inclusions between Mather sets:
$$
 c\in \dpr \beta(h) \quad \Longleftrightarrow \quad h \in \partial \a(c) \quad \Longleftrightarrow \quad \widetilde{\cM}^h \subseteq \widetilde{\cM}_c\,.
 $$
Moreover, for any $c\in \rH^1(M;\R)$:
$$
\widetilde{\cM}_c = \bigcup_{h\in \dpr \a(c)} \widetilde{\cM}^h\,.
$$

({\it iii}) The minimum of the $\a$ function is sometime called {\it Ma\~n\'e's strict critical value}. Observe that if $\a(c_0) = \min \a(c)$, then $0\in \dpr \a(c_0)$ and $\beta(0)=-\a(c_0)$. Therefore, the measures with zero homology are contained in the least possible energy level containing Mather sets: $\widetilde{\cM}^0 \subseteq \widetilde{\cM}_{c_0}$. This inclusion might be strict, unless $\a$ is differentiable at $c_0$; in fact, there may be other action minimizing measures with non-zero rotation vectors corresponding to the other subderivatives of $\a$ at $c_0$.

{(\it iv)} Note that measures of trivial homology are not necessarily supported on orbits with trivial homology or fixed points.
For instance, one can consider the following example (see also \cite[Section 5.2]{ContrerasPaternainMacarini}). Let $M=\T^2$ equipped with the flat metric and consider a vector field $X$ with norm $1$ and such that its orbits form a Reeb foliation, \ie X has two closed orbits $\g_1$ and $\g_2$ in opposite homology classes and any other orbit asymptotically approaches $\g_1$ in forward time and $\g_2$ in backward time. As we have described in section \ref{sec1.1}, we can embed this vector field into the Euler-Lagrange vector field given by the Tonelli Lagrangian 
$L_X(x,v)= \frac{1}{2}\|v-X(x)\|^2$. Let us now consider the probability measure $\m_{\g_1}$ and $\m_{\g_2}$, uniformly distributed respectively on $\g_1$ and $\g_2$. Since these two curves have opposite homologies, then $\rho(\m_{\g_1}) = -\rho(\m_{\g_2})=: h_0 \neq 0$. Moreover, it is easy to see that 
$A_{L_X}(\m_{\g_1}) = A_{L_X}(\m_{\g_2})=0$, since the Lagrangian vanishes on ${\rm Graph}(X)$. 
Using the fact that $L_X\geq 0$ (in particular it is strictly positive outside of ${\rm Graph}(X)$) and that there are no other invariant ergodic probability measures contained in ${\rm Graph}(X)$, we can conclude that  ${\cM_0} = \g_1 \cup \g_2$ and $\a(0)=0$. 
Moreover, $\m_0:= \frac{1}{2}\m_{\g_1} + \frac{1}{2}\m_{\g_2}$ has zero homology and its support is contained in $\widetilde{\cM}_0$. Therefore (see Proposition \ref{CaracMinim} (i)), $\m_0$ is action minimizing with rotation vector $0$ and  $\widetilde{\cM}^0 \subseteq \widetilde{\cM}_0$; in particular,  $\widetilde{\cM}^0 = \widetilde{\cM}_0$. This also implies that 
$\b(0)=0$ and $\alpha(0)= \min\a(c) =0$.\\
Observe that $\a$ is not differentiable at $0$. In fact, reasoning as we have done before for the zero homology class, it is easy to see that for all $t\in[-1,1]$
$\widetilde{\cM}^{th_0} =\widetilde{\cM}_0$. It is sufficient to consider the convex combination
$\m_{\l}=\l \m_{\g_1} + (1-\l)\m_{\g_2}$ for any $\l\in [0,1]$. Therefore, $\dpr \a(0) = \{th_0, \; t\in [-1,1]\}$ and $\beta(th_0)=0$ for all $t\in[-1,1]$.\\
\end{Rem}

As we have just seen in  item ({\rm iv}) of Remark \ref{remarkinclusionsmathersets}, it may happen that the Mather sets corresponding to different homology (resp. cohomology) classes coincide or are included one into the other. This is something that, for instance, cannot happen in the integrable case: in this situation, in fact, these sets form a foliation and are disjoint. The problem in the above mentioned example, seems to be related to a lack of {\it strict convexity} of $\beta$ and $\a$. See also the discussion on the simple pendulum in Addendum 4.B: in this case the Mather sets, corresponding to a non-trivial interval of cohomology classes about $0$, coincide.\\
In the light of this, let us try to understand better what  happens when $\a$ and $\b$ are not strictly convex, \ie when we are in the presence of  ``{\it flat}'' pieces. \\
Let us first fix some notation. If $V$ is a real vector space and $v_0,v_1 \in V$, we shall denote by 
$\s(v_0,v_1)$ the segment joining $v_0$ to $v_1$, that is $\s(v_0,v_1):= \{tv_0 + (1-t)v_1:\; t \in [0,1]\}$. 
We shall say that a function $f: V\longrightarrow \R$ is {\it affine} on $\s(v_0,v_1)$, if there exists $v^*\in V^*$ (the dual of $V$), such that $f(v) = f(v_0) + \langle v^*, v-v_0 \rangle$ for each $v\in \s(v_0,v_1)$.
Moreover, we shall denote by  ${\rm Int}(\s(v_0,v_1))$ the {\it interior} of $\s(v_0,v_1)$, \ie 
${\rm Int}(\s(v_0,v_1)):= \{tv_0 + (1-t)v_1:\; t \in (0,1)\}$.\\

\begin{Prop}\label{flatness}
{\rm(i)} Let $h_0,h_1 \in \rH_1(M;\R)$. $\b$ is affine on $\s(h_0,h_1)$ if and only if for any $h\in {\rm Int}(\s(h_0,h_1))$ we have $\widetilde{\cM}^h \supseteq \widetilde{\cM}^{h_0} \cup \widetilde{\cM}^{h_1}$.\\
{\rm(i)} Let $c_0,c_1 \in \rH^1(M;\R)$. $\a$ is constant on $\s(c_0,c_1)$ if and only if for any $c \in {\rm Int}(\s(c_0,c_1))$ we have $\widetilde{\cM}_c \subseteq \widetilde{\cM}_{c_0} \cap \widetilde{\cM}_{c_1} $.\\
\end{Prop}

\begin{Rem}
The inclusions in Proposition \ref{remarkinclusionsmathersets} may not be true at the end points of $\s$. For instance,  Remark \ref{remarkinclusionsmathersets} ({\it iv}) provides an example in which the inclusion in Proposition \ref{flatness} (i) is not true at the end-points of $\s(-h_0,h_0)$.\\
\end{Rem}

\begin{Proof}
(i) [$\Longrightarrow$] Assume that $\b$ is affine on $\s(h_0,h_1)$, \ie there exists $c\in\rH^1(M;\R)$ such that $\beta(h)=\beta(h_0) + \langle c, h-h_0 \rangle$ for all $h\in\s(h_0,h_1)$.\\
For $i=0,1$, let $\m_i$ be an action minimizing measure with rotation vector $h_i$. If $t\in (0,1)$, let us consider $\m_t:= t \m_1 + (1-t)\m_0$. Clearly, $\rho(\m_t)=t h_1 + (1-t)h_0 \in {\rm Int}(\s(h_0,h_1))$, therefore $\beta(\rho(\m_t))= \beta(h_0) + t \langle c,h_1-h_0 \rangle$. Then:
\beqano
\int_{\rT M} L(x,v) d\m_t &=&  \int_{\rT M} L(x,v) d(t \m_1 + (1-t)\m_0) = \\
&=& t \int_{\rT M} L(x,v) d\m_1 + (1-t) \int_{\rT M} L(x,v) d\m_0 = \\
&=& t\beta(h_1) +(1-t)\beta(h_0) = \\
&=& t(\beta(h_0) + \langle c, h_1 -h_0 \rangle)  +(1-t)\beta(h_0) = \\
&=& \beta(h_0) + t \langle c,h_1-h_0 \rangle = \beta(\rho(\m_t))\,.
\eeqano
Therefore $\m_t$ is action minimizing with rotation vector $t h_1 + (1-t)h_0$. Since this is true for any $\m_i$ with rotation vector $h_i$ ($i=0,1$), it follows that:
$$
 \widetilde{\cM}^{h_0} \cup \widetilde{\cM}^{h_1} \subseteq \widetilde{\cM}^{th_1+ (1-t)h_0}.
 $$
 \noindent [$\Longleftarrow$] Since 
 $\widetilde{\cM}^h \supseteq \widetilde{\cM}^{h_0} \cup \widetilde{\cM}^{h_1}$ for all $h\in {\rm Int}(\s(h_0,h_1))$, using the observation in Remark \ref{remarkinclusionsmathersets} ({\it ii}) and Proposition \ref{CaracMinim}, we obtain that $\widetilde{\cM}^{h_0}$ and $\widetilde{\cM}^{h_1}$
 must be  contained in $\widetilde{\cM}_{c_0}$ for some $c_0\in \rH^1(M;\R)$. In particular, $c_0 \in \dpr \b(h_0) \cap \dpr \b(h_1)$. Let us show that $c_0 \in \dpr \beta(h)$ for all $h\in {\rm Int}(\s(h_0,h_1))$. 
 In fact, let $h=th_1+(1-t)h_0$ for some $t\in(0,1)$. Then using the convexity of $\beta$:
\beqano
 \a(c_0) +\beta(h) &=& \a(c_0) +\beta(th_1+(1-t)h_0) \leq \\
 &\leq& t \big(\a(c_0) +\beta(h_1)\big) + (1-t)\big(\a(c_0) +\beta(h_0)\big) =\\
 &\leq& t\langle c_0,h_1 \rangle + (1-t) \langle c_0, h_0 \rangle = \\
 &=& \langle c_0, th_1 + (1-t) h_0 \rangle = \langle c_0, h\rangle.
 \eeqano
 On the hand, the reverse inequality is always true (Fenchel-Legendre inequality).\\
  Now using this fact, it follows that  if $h\in \s(h_0,h_1)$ then:
 $$
 \a(c_0) +\beta(h)= \langle c_0, h \rangle  \qquad {\rm and } \qquad    \a(c_0) +\beta(h_0)= \langle c_0, h_0 \rangle
 $$
 therefore subtracting the second equality from the first one, we obtain what we wanted:
$$
 \beta(h) -   \beta(h_0) =  \langle c_0, h-h_0 \rangle  \quad \forall \; h\in \s(h_0,h_1)\,.
 $$
 (ii) [$\Longleftarrow$] It is a trivial consequence of Theorem \ref{teocarneiro}. \\
  $[\Longrightarrow]$ Let $t\in (0,1)$ and consider $c_t:=tc_1+(1-t)c_0$. We want to show that 
  for each $h\in \dpr \a(c_t)$ we have that $\langle c_1-c_0,h \rangle=0$. First observe that 
  since $c_t \in {\rm Int}(\s(c_0,c_1))$, then there exists $\d>0$ such that
  $c_t+ s (c_1-c_0) \in \s(c_0,c_1)$ for all $s\in(-\d,\d)$. Then, using the fact that $h\in \dpr \a(c_t)$ and that $\a$ is constant in $\s(c_0,c_1)$, we obtain:
  $$
  0= \a(c_t + s(c_1-c_0)) - \a(c_t) \geq s\langle c_1 -c_0, h \rangle \qquad \forall\; s\in (-\d,\d).
  $$
  But this can be true only if $\langle c_1 -c_0, h \rangle=0$.\\
  Now let us prove that,  for $i=0,1$, $\dpr \a(c_i)\supseteq \dpr \a(c_t)$  for all $t\in (0,1)$ . 
  In fact, if $h\in \dpr \a(c_t)$, then using that $\langle c_1,h\rangle = \langle c_0,h\rangle$ and that $\a$ is constant on $\s(c_0,c_1)$, we get:
  \beqano
  \a(c_i) + \beta(h) &=& \a(c_t) + \beta(h) = \langle c_t, h\rangle =\\
  &=& \langle tc_1 + (1-t)c_0, h\rangle =\\
  &=& t\langle c_1, h\rangle + (1-t)\langle c_0, h\rangle = \\
  &=& t\langle c_i, h\rangle + (1-t)\langle c_i, h\rangle = \\ 
&=&  \langle c_i, h\rangle.
  \eeqano
  This and   Remark \ref{remarkinclusionsmathersets} ({\it ii}) immediately allows us to conclude that
 $$
\widetilde{\cM}_c \subseteq \widetilde{\cM}_{c_i} \quad {\rm for}\; i=0,1 \quad \Longleftrightarrow \quad
\widetilde{\cM}_c \subseteq \widetilde{\cM}_{c_0} \cap \widetilde{\cM}_{c_1}.
 $$

  \end{Proof}

\begin{Rem}
It follows from the previous remarks and Proposition \ref{flatness}, that, in general, the action minimizing measures (and consequentely the mather sets $\widetilde{\cM}_c$ or $\widetilde{\cM}^h$) are not necessarily ergodic. Recall that an invariant probability measure is said to be {\it ergodic}, if all invariant Borel sets have measure $0$ or $1$. These measures play a special role in the study of the dynamics of the system, therefore one could ask what are the ergodic action-minimizing measures. It is a well-known result from ergodic theory, that the ergodic measures  of a flow correspond to the {\it extremal points} of the set of invariant probability measures  (see for instance \cite{KB}), where by ``extremal point'' of a convex set, we mean an element that cannot be obtained as a non-trivial convex combination of other elements of the set.
Since $\beta$ has superlinear growth, its epigraph $\{(h,t)\in \rH_1(M;\R)\times \R:\; t\geq \beta(h)\}$ has infinitely many extremal points. Let $(h,\beta(h))$ denote one of these extremal points. Then, there exists at least one ergodic action minimizing measure with rotation vector $h$. It is in fact sufficient to consider any extremal point of the set $\{\mu\in \calM^h(L):\; A_L(\m)=\beta(h)\}$: this measure will  be an extremal point of $\calM(L)$ and hence ergodic.
Moreover, as we have already recalled in Remark \ref{Remrotationvector}, for such an ergodic measure $\mu$, Birkhoff's ergodic theorem implies that $\m$-almost every trajectory of $\Phi^L$ has rotation vector $h$.\\
\end{Rem}

\vspace{10 pt}

\noindent{\bf ADDENDA}\\

\noindent{\bf {\sc 4.A - The symplectic invariance of Mather sets}}\\
\addcontentsline{toc}{subsection}{\hspace{15 pt} 4.A: The symplectic invariance of Mather sets}

In this addendum we would like to discuss some symplectic aspects of the theory that we have started to develop. In particular, as a first step, we would like to understand how the Mather sets behave under the action of {\it symplectomorphisms} (see also \cite{BernardSympl}). In order to do this, we need to move to the Hamiltonian setting, rather than the Lagrangian one, and consider the Hamiltonian  $H: \rT^*M \longrightarrow \R$, defined by Fenchel duality (see Section \ref{sec1.1}). Recall that the associated Hamiltonian flow $\Phi^H$ is conjugate to the Euler-Lagrange flow of $L$, therefore, from a dynamical systems point of view, the two systems are equivalent. We can define the associated Mather sets in the cotangent space as 
$$
\cM_c^*(H) := \cL_L(\widetilde{\cM}_c(L)),
$$
where $\cL_L: \rT M \longrightarrow \rT^*M$ is the Legendre transform associated to $L$ (see (\ref{Legendretransform})).
These sets are non-empty, compact and invariant for the Hamiltonian flow and they continue to satisfy the graph property (see Theorems \ref{Theograph} and \ref{matherrotation}).

 The main advantage of shifting our point of view, is that the cotangent space ${\rT^*M}$ can be naturally equipped with an exact symplectic form $\omega =  \sum_{i=1}^d dx_i \wedge dp_i = -d\l$, where $\lambda = \sum_{i=1}^d p_i\,dx_i\,$ is what is called the {\it tautological form} or {\it Liouville form}. \\
 %See Appendix \ref{sec3.1} for a more intrinsic definition of $\l$.\\
Let $\Psi: \rT^*M \longrightarrow \rT^*M$ be a diffeomorphism. We shall say that $\Psi$ is a {\it symplectomorphism} if it preserves the symplectic form $\omega$, \ie $\Psi^*\omega=\omega$. 
An easy class of examples is provided by {\it translations in the fibers}. Let $\eta$ be a closed $1$-form and consider
$\tau_{\eta} : \rT^*M \longrightarrow \rT^*M$, $(x,p)\longmapsto (x,p+\eta(x))$. Then, $\t$ is a symplectomorphism.\\
Observe  that since $\omega = -d\l$, one has that $d(\Psi^*\l-\l)=0$; in other words,  the $1$-form $\Psi^*\l-\l$ is closed. We shall call the {\it cohomology class } of $\Psi$, the cohomology class of $\Psi^*\l-\l$ and denote it by $[\Psi]\in \rH^1(\rT^*M;\R) \simeq \rH^1(M;\R)$, since $\rT^*M$ can be retracted to $M$ along the fibers (hereafter we shall always identify these two spaces). Obviously, going back to the previous example, $[\tau_{\eta}]=[\eta]$. In particular, $\Psi$ is said to be {\it exact} if and only if $[\Psi]=0$.

\begin{Lem}\label{symplectomorphism}
Any symplectomorphism $\Psi: \rT^*M \longrightarrow \rT^*M$ can be written as $\Psi = \Phi \circ \t_{\eta}$, where $\Phi$ is an exact symplectomorphism and $[\eta]=[\Psi]$.
\end{Lem}

\begin{Proof}
Let $\eta$ be any closed $1$-form on $M$, such that $[\eta]=[\Psi]$, and define $\tau_{\eta}$ as above.
Moreover, we define $\Phi = \Psi \circ (\tau_{\eta})^{-1} = \Psi \circ \tau_{-\eta}$. In order to conclude the proof of the Lemma, we need to check that $\Phi$ is exact. In fact,
$$
\Phi^*\l = (\Psi \circ \tau_{-\eta})^*\l = \tau_{-\eta}^*(\Psi^* \l) = \tau_{-\eta}^*(\l + \Theta) = \l + \tau_{-\eta}^*\Theta - \eta,
$$
where $[\tau_{-\eta}^*\Theta]=[\Psi]=[\eta]$.Then, $[\Phi^*\l -\l]=[\Theta - \eta]=0$ and hence it is exact.
\end{Proof}

Therefore, if we want to understand the interplay between Mather sets and the action of symplectomorphisms, it will be sufficient to analyze the behaviour of these two kinds of symplectomorphisms: {\it translations in the fibers} and {\it exact symplectomorphisms}.

\begin{Prop}\label{4c1}
Let $H:\rT^*M\longrightarrow \R$ be a Tonelli Hamiltonian and $\eta$ a closed $1$-form on $M$. Then:
$$\cM^*_c (H\circ \tau_{\eta}) = \tau_{-\eta}\left( \cM^*_{c+[\eta]}(H)\right) \qquad \forall\, c\in \rH^1(M;\R).$$
\end{Prop}

\begin{Rem}
If $H$ is a Tonelli Hamiltonian, then also $H\circ \tau_{\eta}$ is a Tonelli Hamiltonian (we are just composing it with a vertical translation in the fibers).
\end{Rem}

\begin{Proof}
Let us start by observing that the Lagrangian associated to $H\circ \tau_{\eta}$ is $L_{\eta}:=L-\hat{\eta}$ and that the associated Legendre transform $\cL_{L_{\eta}} = \tau_{-\eta} \circ \cL$ (just  derive $L_{\eta}$ with respect to $v$).
Therefore, using the definition of Mather sets for a given cohomology class, we get:
\beqano
\cM_c^*(H\circ \tau) &:=& \cL_{L_{\eta}}\left(\widetilde{\cM}_c(L-\hat{\eta})\right) = \cL_{L_{\eta}}\left(\widetilde{\cM}_{c+[\eta]}(L)\right) = \\
&=&(\tau_{-\eta}\circ \cL_L)\left( \widetilde{\cM}_{c+[\eta]}(L)\right) = 
\tau_{-\eta} \left( \cM^*_{c+[\eta]}(H)\right).
\eeqano
\end{Proof}

Let us see now what happens with exact symplectomorphisms (see also \cite{IntegTonelli}).

\begin{Prop}\label{4c2}
Let $H:\rT^*M\longrightarrow \R$ be a Tonelli Hamiltonian and
$\Phi: \rT^*M \longrightarrow \rT^*M$ an exact symplectomorphism such that $H\circ \Phi$ is still of Tonelli type. Then:
$$\cM^*_c (H\circ \Phi) = \Phi^{-1} \left( \cM^*_{c}(H)\right) \qquad \forall\, c\in \rH^1(M;\R).$$
\end{Prop}

\begin{Rem}
Observe that even if $H$ is Tonelli, $H\circ \Phi$ is not necessarily  of Tonelli type!
\end{Rem}

This proposition can be easily deduced from the following Lemma (see also \cite{IntegTonelli}).

\begin{Lem}\label{lemma4.30}
Let $H:\rT^*M\longrightarrow \R$ be a Tonelli Hamiltonian and $\Phi: \rT^*M \longrightarrow \rT^*M$ an exact symplectomorphism, such that $H':=H\circ \Phi$ is still Tonelli. Then:\\
{\rm (i)}  $\m$ is an invariant probability measure of $H$ if and only if  $\Phi^*\mu$ is an invariant probability measure of $H'=H\circ \Phi$;\\
{\rm (ii)} for any $\m$ invariant probability measure of $H$, the following holds:
\beqa{identita}
\int \left[p\frac{\dpr H}{\dpr p}(x,p) - H(x,p)\right] d \m =
\int \left[p\frac{\dpr H'}{\dpr p}(x,p) - H'(x,p)\right] d \Phi_*\m.
\eeqa
\end{Lem}

\begin{Rem}
The identity in (\ref{identita}) represents the equality between the associated Lagrangian actions (in the Hamiltonian formalism). Therefore, since the respective actions coincide on all measures and the invariant measures are in $1-1$ correspondence, there must be a correspondence between the minimizing ones: {\it $\m$ is action minimizing for $H$ if and only if $\Phi^*\mu$ is action-minimizing for $H\circ\Phi$}. From this, Proposition \ref{4c2} follows easily. 
\end{Rem}

\begin{Proof}[{\bf Lemma \ref{lemma4.30}}]
(i) The first part of the statement follows from the classical fact that $\Phi$ tranforms the associated Hamiltonian vector fields in the following way:
\beqa{trasfcampivett}
X_{H} (\Phi(x,p)) = D\Phi(x,p) X_{H\circ \Phi} (x,p) \qquad \forall\, (x,p)\in \rT^*M.
\eeqa
(ii) If we denote by $\l(x,p)$ the Liouville form $pdx$, then: %and by $X_H(x,p)$ the Hamiltonian vector field, then:
%Let us denote $\Phi(x,p)=(\Phi_1(x,p),\Phi_2(x,p))$. Using Fenchel-Legendre equality,  the coordinate representation of $\cL^{-1}(x,p)=\big(x,\frac{\dpr H}{\dpr p}(x,p)\big)$, the fact that $H$ is preserved by $\Phi$ and the exactness of $\Phi$, we get:
\beqano
\int \left[p\frac{\dpr H}{\dpr p}(x,p) - H(x,p)\right] d \m = 
\int \Big(\l(x,p)[X_H(x,p)] - H(x,p)\Big) d \m\,.
\eeqano
Therefore, using that $|\det D \Phi| =1$, the relation in (\ref{trasfcampivett})  
and that $\Phi^*\l - \l = df$ (since $\Phi$ is an exact symplectomorphism),  we obtain: 
\beqano
&& \int \left[p\frac{\dpr H}{\dpr p}(x,p) - H(x,p)\right]\;d\m = 
\int \Big(\l(x,p)[X_H(x,p)] - H(x,p)\Big)\;d \m = \\
&&  = \int \Big(\Phi^*\l(x',p')[D\Phi^{-1}(x',p')X_{H}(\Phi(x',p'))] - H(\Phi(x',p'))\Big)\;d \Phi^*\m = \\
&&  = \int \Big(\big(\l(x',p') + df(x',p')\big) [X_{H\circ \Phi}(x',p')] - H(\Phi(x',p'))\Big)\;d \Phi^*\m =\\
&& = \int \Big(\l(x',p') [X_{H'}(x',p')] - H'(x',p')\Big)\;d \Phi^*\m +  \int df(x',p')[X_{H'}(x',p')] \,d\Phi^*\m=\\
&& = \int \Big(p'\frac{\dpr H'}{\dpr p'}(x',p') - H(x',p')\Big)\;d \Phi^*\m. 
\eeqano
In the last equality we used that 
$\int df(x',p')[X_{H'}(x',p')] \,d\Phi^*\m =0$, as it follows easily from the invariance of $\Phi^*\m$. For the sake of simplifying the notation, let us denote $\nu:= d\Phi^*\m$ and assume that $\nu$ is ergodic (otherwise consider each ergodic component). Using the ergodic theorem and the compactness of the suppport of $\nu$, we obtain that for a generic point $(x_0,y_0)$ in the support of $\nu$:
\beqano
\int df(x',p')[X_{H'}(x',p')] \,d\nu &=& \lim_{N\rightarrow +\infty} \frac{1}{N} \int_0^N
df(\Phi^{H'}_t(x_0,p_0))[X_{H'}(\Phi^{H'}_t(x_0,p_0))]\,d t = \\
&=& \lim_{N\rightarrow +\infty} \frac{f(\Phi^{H'}_N(x_0,p_0)) - f(x_0,p_0)}{N} =0\,. 
\eeqano

%where for the last equality we used that $\int df(x)\cdot v\,d(\cL^*\m)=0$, as it follows easily from the invariance of the measure under the flow (see for instance \cite{Mather91}).
\end{Proof}

Finally, we can deduce the main result of this addendum.

\begin{Teo}\label{4cteofinale}
Let $H:\rT^*M\longrightarrow \R$ be a Tonelli Hamiltonian and $\Psi: \rT^*M \longrightarrow \rT^*M$ a symplectomorphisms of class $[\Psi]$, such that $H\circ \Psi$ is of Tonelli type. Then,
$$\cM^*_c (H\circ \Psi) = \Psi^{-1} \left( \cM^*_{c+[\Psi]}(H)\right) \qquad \forall\, c\in \rH^1(M;\R).$$
\end{Teo}

\begin{Proof}
Suppose that $\Psi = \Phi \circ \tau_{\eta}$, with $[\eta]=[\Psi]$. Then, using Propositions \ref{4c1} and \ref{4c2}, we get for all $c\in\rH^1(M;\R)$:
\beqano
\cM^*_c(H\circ \Psi) &=& \cM^*_c(H\circ \Phi \circ \tau_{\eta}) = \tau_{-\eta}\left( \cM^*_{c+[\eta]}(H\circ \Phi)\right) =\\
&=& \tau_{\eta}^{-1}\left( \Phi^{-1} \left(\cM^*_{c+[\eta]}(H)\right)\right) = (\Phi \circ \tau_{\eta})^{-1} \left(\cM^*_{c+[\eta]}(H)\right) =\\
&=& \Psi^{-1} \left(\cM^*_{c+[\Psi]}(H)\right)\,.
\eeqano
\end{Proof}

\begin{Cor}
Let $H:\rT^*M\longrightarrow \R$ be a Tonelli Hamiltonian and $\Psi: \rT^*M \longrightarrow \rT^*M$ a symplectomorphisms of class $[\Psi]$, such that $H\circ \Psi$ is of Tonelli type. Then:
\beqano
{\rm(i)} && \a_{H\circ \Psi} (c)  = \a_H (c+[\Psi])  \qquad \qquad \forall\, c\in \rH^1(M;\R)\\
{\rm(ii)} && \b_{H\circ \Psi} (h)  = \b_H (h) - \langle [\Psi], h \rangle \qquad\, \forall\, h\in \rH_1(M;\R).\\
\eeqano
\end{Cor}

\begin{Proof}
(i) From Theorem \ref{4cteofinale} it follows that if  $(x,p)\in \cM^*(H\circ \Psi)$, then $\Psi(x,p)\in \cM_{c+[\Psi]}^*(H)$. Therefore, using Theorem \ref{teocarneiro}, for any $(x,p)\in \cM^*(H\circ \Psi)$ we have:
\beqano
-\a_{H\circ \Psi} (c) = (H\circ \Psi)(x,p) = H(\Psi(x,p)) =  - \a_H(c+[\Psi]).
\eeqano
(ii) Observe that if $h\in \dpr \a_{H\circ\Psi}(c)$, then $h\in \dpr \a_{H}(c+[\Psi])$. In fact, for every $c'\in \rH^1(M;\R)$:
\beqano
\a_{H}(c'+[\Psi]) - \a_{H}(c'+[\Psi]) &=& \a_{H \circ \Psi}(c') - \a_{H \circ\Psi}(c) \;\geq\; \langle c'-c, h\rangle  =\\
&=& \langle (c'+[\Psi])-(c-[\Psi]), h\rangle. 
\eeqano
Let now  $h\in \dpr \a_{H\circ\Psi}(c)$. Using the fact that $\a$ and $\b$ are one the conjugate of the other, we get:
\beqano
\b_{H\circ\Psi}(h) &=& \langle c, h \rangle - \a_{H\circ\Psi}(c) = \langle c, h \rangle - \a_{H}(c+[\Psi]) = \\
&=&  \langle c + [\Psi], h \rangle  - \a_{H}(c+[\Psi]) - \langle [\Psi], h \rangle=\\
&=&  \beta_{H}(h) - \langle [\Psi], h \rangle\,.
\eeqano
\end{Proof}

We can summarize everything in the following commutative diagram.\\

$$
\xymatrix{
\ar@(l,u)^{\Phi^{H\circ\Psi}_t} \rT^*M \ar[rr]^{\Psi} \ar[rd]_{H\circ\Psi} 		&& \ar[ld]^{H}  \rT^*M \ar@(r,u)[]_{\Phi^H_t} \\ 
		&\R \\
\rH^1(M;\R) \ar[rr]^{c\, \mapsto c+[\Psi]} \ar[ru]^{\a_{H\circ\Psi}}\ar@{~>}[uu]^{c\mapsto \cM^*_c(H\circ\Psi)}&& \ar[lu]_{\a_H} \ar@{~>}[uu]_{c\mapsto \cM^*_c(H)}\rH^1(M;\R)\\ 
}
$$

\vspace{20pt}

\vspace{20 pt}

\noindent{\bf {\sc 4.B - An example: the simple pendulum I}}\\
\addcontentsline{toc}{subsection}{\hspace{15 pt} 4.B: An example: the simple pendulum I}

In this addendum we would like to describe the Mather sets, the $\a$-function and the $\beta$-function, in a specific example: the {\it simple pendulum}. This system can be  described in terms of the Lagrangian:
\beqano
L: \rT \T &\longrightarrow& \R \\
(x,v) &\longmapsto& \frac{1}{2}|v|^2 + \big(1-\cos (2\pi x)\big).
\eeqano
It is easy to check that the Euler-Lagrange equation provides exactly the equation of the pendulum:
$$
\dot{v} = 2\pi \sin (2\pi x) \qquad \Longleftrightarrow \qquad \left\{\begin{array}{l} v=\dot{x} \\ \ddot{x} - 2\pi\sin(2\pi x) = 0. \end{array}\right.
$$

\begin{figure} [h!]
\begin{center}
\includegraphics[scale=0.3]{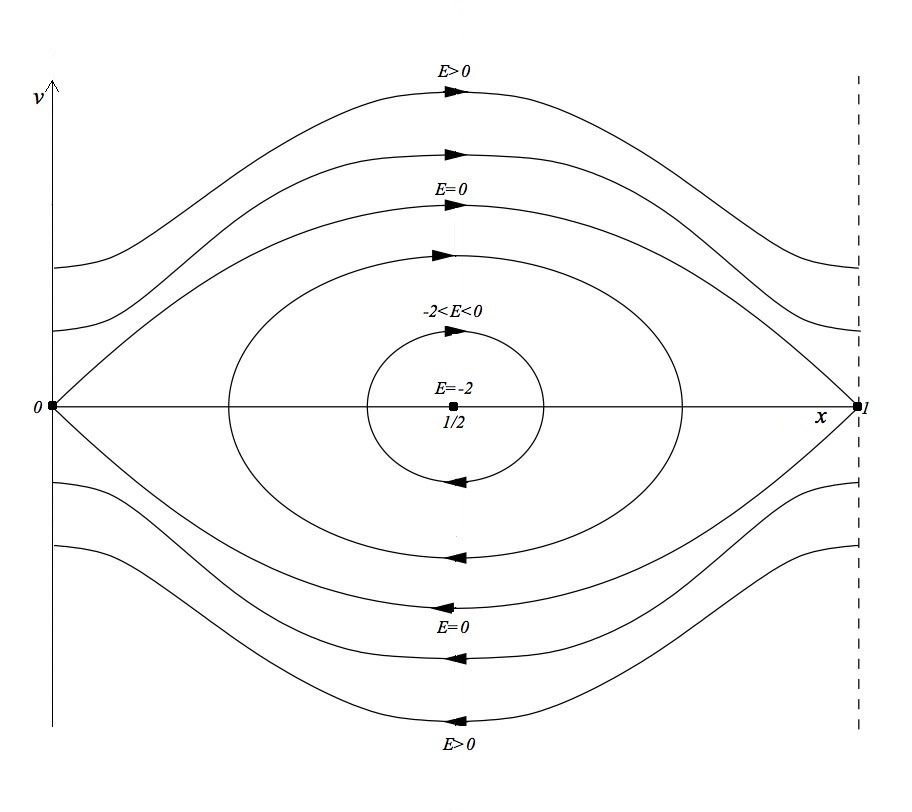}
\caption{The phase space of the simple pendulum.}
\end{center}
\end{figure}

The associated Hamiltonian (or energy) $H: \rT^* \T \longrightarrow \R$ is given by $H(x,p):= \frac{1}{2}|p|^2 - (1- \cos (2\pi x))$. 
Observe that in this case the Legendre transform $\cL_L(x,v)=(x,p)$, therefore we can easily identify the tangent and cotangent space. In the following we shall consider $\rT \T \simeq \rT^*\T \simeq \T\times \R$ and
identify $\rH^1(M;\R)\simeq \rH_1(M;\R) \simeq \R$.\\

First of all, let us study what are the invariant probability measures of this system. 
\begin{itemize}
\item Observe that $(0,0)$ and $(\frac{1}{2},0)$ are fixed points for the system (respectively {\it unstable} and {\it stable}). Therefore, the Dirac measures concentrated on each of them are invariant probability measures. Hence, we have found two first invariant measures:  $\delta_{(0,0)}$ and $\delta_{(\frac{1}{2},0)}$, both with zero rotation vector: $\rho(\delta_{(0,0)})= \rho(\delta_{(\frac{1}{2},0)})=0$. As far as their energy is concerned (\ie the energy levels in which they are contained), it is easy to check that $E(\delta_{(0,0)})=H(0,0)=0$ and $E(\delta_{(\frac{1}{2},0)})=H(\frac{1}{2},0)=-2$. Observe that these two energy levels cannot contain any other invariant probability measure.
\item If $E>0$, then the energy level $\{H(x,v)=E\}$ consists of two homotopically non-trivial periodic orbits ({\it rotation motions}):
$$
{\mathcal P}^{\pm}_E := \{(x,v): \; v=\pm \sqrt{2[(1+E)-\cos(2\pi x)]},\; \forall\, x\in\T\}.
$$
The probability measures evenly distributed along these orbits - which we shall denote $\m^{\pm}_E$ - are invariant probability measures of the system. If we denote by 
\beqa{periodorbit}
T(E):= \int_0^1\frac{1}{\sqrt{2[(1+E)-\cos(2\pi x)]}}\,dx
\eeqa
the period of such orbits, then it is easy to check that (see Remark \ref{Remrotationvector})
$\rho(\mu^{\pm}_E) = \frac{\pm1}{T(E)}$. Observe that this function $T:(0,+\infty) \longrightarrow (0,+\infty)$, which associates to a positive energy $E$ the period of the  corresponding periodic orbits ${\mathcal P}_E^{\pm}$, is continuous and strictly decreasing. Moreover, $T(E) \to \infty$ as $E\to 0$ (it is easy to see this, noticing that motions on the {\it separatrices} take an ``{infinite}'' time to connect $0$ to $1\equiv0$ mod.$1$). Therefore, 
$\rho(\mu^{\pm}_E) \to 0$ as $E\to 0$.
\item If $-2<E<0$, then the energy level $\{H(x,v)=E\}$ consists of one contractible periodic orbit ({\it libration motion}):
$$
{\mathcal P}_E := \{(x,v): \; {v^2}= {2(1+E)-2\cos(2\pi x)},\quad x\in [x_E, 1-x_E]\},
$$
where $x_E:=\frac{1}{2\pi}\arccos(1+E)$.
The probability measure evenly distributed along this orbit - which we shall denote $\m_E$ - is  an invariant probability measure of the system. 
Moreover, since this orbit is contractible, its rotation vector is zero: $\rho(\mu_E)=0$.
\end{itemize}

The measures above are the only ergodic invariant probability measures of the system. Other invariant measures can be easily obtained as convex combination of them.\\

Now we want to understand which of these are action-minimizing for some  cohomology class. 

\begin{Rem}\label{remarkpendulum}
({\it i}) Let us start by remarking that for $-2<E<0$ the support of the measure $\mu_E$ is not a graph over $\T$, therefore it cannot be action-minimizing for any cohomology class, since otherwise it would violate Mather's graph theorem (Theorems \ref{Theograph} and \ref{matherrotation}). Therefore all action-minimizing measures will be contained in energy levels corresponding to energy bigger than zero. It follows from Theorem \ref{teocarneiro}, that $\a(c) \geq 0$ for all $c\in \R$.

({\it ii}) Another interesting property of the $\a$ function (in this specific case) is that it is an {\it even function}: $\a(c)=\a(-c)$ for all $c\in\R$. This is a consequence of the particular symmetry of the system, \ie $L(x,v)=L(x,-v)$. In fact, let  us denote $\t: \T\times \R \longrightarrow \T\times \R$, $(x,v) \longmapsto (x,-v)$ and observe that if $\m$ is an invariant probability measure, then also $\t^*\mu$ is still an invariant probability measure. Moreover, $\t^* \calM(L) = \calM(L)$, where $\calM(L)$ denotes the set of all invariant probability measures of $L$.
It is now sufficient to notice that for each $\m \in \calM(L)$,
$
\int (L-c\cdot v)\,d\mu = \int (L+c\cdot v)d\tau^*\mu,
$
and hence conclude that 
$$
\a(c)= - \inf_{\calM(L)} \int (L-c\cdot v)\,d\mu = - \inf_{\calM(L)} \int (L+c\cdot v)d\tau^*\mu = \a(-c)\,.
$$

({\it iii}) It follows from the above symmetry and the convexity of $\a$, that $\min_\R\a(c) = \a(0)$.
\end{Rem}

Let us now start by studying the $0$-action minimizing measures, \ie invariant probability measures that minimize the action of $L$ without any ``correction''. Since  $L(x,v)\geq 0$ for each $(x,v)\in \T\times\R$, then $A_L(\mu)\geq 0$ for all
$\mu\in\calM(L)$. In particular, $A_L(\delta_{(0,0)})=0$, therefore $\d_{(0,0)}$ is a $0$-action minimizing measure and $\a(0)=0$. Since there are not other invariant probability measures supported in the energy level $\{H(x,v)=0\}$ (\ie on the separatrices), then we can conclude that:
$$
\widetilde{\cM}_0 = \{(0,0)\}\,.
$$
Moreover, since $\a'(0)=0$ (see Remark \ref{remarkpendulum} ({\it iii})), then it follows from Remark \ref{remarkinclusionsmathersets} that:
$$
\widetilde{\cM}^0 = \widetilde{\cM}_0= \{(0,0)\}\,.
$$
On the other hand, this could be also deduced from the fact that the only other measures with rotation vector $0$, cannot be action minimizing since they do not satisfy the graph theorem (Remark \ref{remarkpendulum}).\\

Now let us investigate what happens for other cohomology classes.
A na\"ive observation is that since the $\a$ function is superlinear and continuous, all energy levels for $E>0$ must contain some Mather set; in other words, they will be achieved for some $c$.\\
Let $E>0$ and consider the periodic orbit ${\mathcal P_E^+}$ and the invariant probability measure $\mu_E^+$ evenly distributed on it. The graph of this orbit can be seen as the graph of a closed $1$-form
$
\eta^+_E :=  \sqrt{2[(1+E)-\cos(2\pi x)]}\,dx,
$
whose cohomology class is
\beqa{cohomologyorbit}
c^+(E):=[\eta^+_E] = \int_0^1 \sqrt{2[(1+E)-\cos(2\pi x)]}\,dx,
\eeqa
which can be interpreted as the (signed) area between the curve and the positive $x$-semiaxis. This value is clearly continous and strictly increasing with respect to $E$ (for $E>0$) and, as $E\to 0$:
$$
c^+(E)  \longrightarrow \int_0^1\sqrt{2[1-\cos(2\pi x)]}\,dx = \frac{4}{\pi}\,.
$$
Therefore, it defines an invertible function $c^+: (0,+\infty) \longrightarrow (\frac{4}{\pi},+\infty)$.\\
We want to prove that $\mu^{+}_E$ is $c^+(E)$-action minimizing. The proof will be an imitation of what already seen for KAM tori in Section \ref{cartoonexample} (see Proposition \ref{propminimalactiontorus}).\\
Let us consider the Lagrangian $L_{\eta_E^+}(x,v):= L(x,v) - \eta_E^+(x)\cdot v$. Then, using Fenchel-Legendre inequality (\ref{Fenchelineq}) (on the support of $\mu^+_E$, because of  our choice of $\eta_E^+$, this is indeed an equality):
\beqano
\int L_{\eta_E^+}(x,v) d\mu_E^+ &=& \int \left(L(x,v)- {\eta_E^+}(x)\cdot v\right) d\mu_E^+ = \\
&=& \int - H(x, {\eta_E^+}(x)) d\mu_E^+ =  - E\,.
\eeqano
Now, let $\nu$ be any other invariant probability measure and apply again the same procedure as above (warning: this time Fenchel-Legendre inequality is not an equality anymore!):
\beqano
\int L_{\eta_E^+}(x,v) d\nu &=& \int \left(L(x,v)- {\eta_E^+}(x)\cdot v\right) d\nu \geq \\
&\geq& \int - H(x, {\eta_E^+}(x)) d\nu =  - E\,.
\eeqano
Therefore, we can conclude that $\mu^{+}_E$ is $c^+(E)$-action minimizing. Since it already projects over the whole $\T$, it follows from the graph theorem that it is the only one:
$$
\widetilde{\cM}_{c^+(E)}= {\mathcal P_E^+} = \{(x,v): \; v=  \sqrt{2[(1+E)-\cos(2\pi x)]},\; \forall\, x\in\T\}.
$$
Furthermore, since $\rho(\mu_E^+) = \frac{1}{T(E)}$, then:
$$
\widetilde{\cM}^{\frac{1}{T(E)}}= \widetilde{\cM}_{c^+(E)}={\mathcal P_E^+}.
$$

Similarly, one can consider the periodic orbit ${\mathcal P_E^-}$ and the invariant probability measure $\mu_E^-$ evenly distributed on it. The graph of this orbit can be seen as the graph of a closed $1$-form
$
\eta^-_E :=  - \sqrt{2[(1+E)-\cos(2\pi x)]}\,dx = - \eta^+_E,
$
whose cohomolgy class is 
$c^-(E)=-c^+(E)$. Then (see also Remark \ref{remarkpendulum} ({\it ii})): 
$$
\widetilde{\cM}_{c^-(E)}= {\mathcal P_E^-} = \{(x,v): \; v=- \sqrt{2[(1+E)-\cos(2\pi x)]},\; \forall\, x\in\T\},
$$
and
$$
\widetilde{\cM}^{-\frac{1}{T(E)}}= \widetilde{\cM}_{c^-(E)}={\mathcal P_E^-}.
$$

\vspace{10 pt}

Note that this completes the study of the Mather sets for any given rotation vector, since
$$
\rho(\mu^{\pm}_E) = \pm\frac{1}{T(E)} \stackrel{E\to+\infty}{\longrightarrow} \pm \infty
\qquad {\rm and}\qquad \rho(\mu^{\pm}_E) = \pm\frac{1}{T(E)} \stackrel{E\to 0^+}{\longrightarrow} 0\,.
$$
What remains to study is what happens for non-zero cohomology classes in $[-\frac{4}{\pi},\frac{4}{\pi}]$. The situation turns out to be quite easy. Observe that $\a(c^{\pm}(E))=E$. Thefore, 
from the continuity of  $\a$ it follows that (take the limit as $E\to 0$): $\a(\pm \frac{4}{\pi})=0$. Moreover, since 
$\a$ is convex and $\min \a(c)=\a(0)= 0$, then: $\a(c)\equiv 0$ on $[-\frac{4}{\pi},\frac{4}{\pi}]$.
Therefore, the corresponding Mather sets will lie in the zero energy level. From the above discussion, it follows that in this energy level there is a unique invariant probability measure, namely $\d_{(0,0)}$, and consequently:
$$
\widetilde{\cM}_{c}= \{(0,0)\} \qquad \mbox{for all}\; -\frac{4}{\pi} \leq c \leq \frac{4}{\pi}.
$$

Let us  summarize what we have found so far. Recall that 
in (\ref{periodorbit}) and (\ref{cohomologyorbit}) we have introduced these two functions:
$
T: (0,+\infty) \longrightarrow (0,+\infty)
$
and
$
c^+: (0,+\infty) \longrightarrow (\frac{4}{\pi},+\infty)
$
representing respectively the period and the ``cohomology'' (area below the curve) of the ``upper'' periodic orbit of energy $E$. These functions (for which we have an explicit formula in terms of $E$) are continuous and strictly monotone (respectively, decreasing and increasing). Therefore, we can define their inverses which provide the energy of the  periodic orbit with period $T$ (for all positive periods) or  the energy of the periodic orbit with cohomology class $c$  (for $|c|>\frac{4}{\pi}$). We shall denote them $E(T)$ and $E(c)$ (observe that this last quantity is exactly the $-\a(c)$). Then:
$$
\widetilde{\cM}_c= \left\{\begin{array}{lll}
\{(0,0)\} && {\rm if}\; -\frac{4}{\pi}\leq c \leq \frac{4}{\pi}\\
{\mathcal P}^+_{E(c)} && {\rm if}\; c > \frac{4}{\pi}\\
{\mathcal P}^-_{E(-c)} && {\rm if}\; c < -\frac{4}{\pi}
\end{array}\right.
$$
and
$$
\widetilde{\cM}^h= \left\{\begin{array}{lll}
\{(0,0)\} && {\rm if}\; h=0\\
{\mathcal P}^+_{E(\frac{1}{h})} && {\rm if}\; h > 0\\
{\mathcal P}^-_{E(-\frac{1}{h})} && {\rm if}\; h < 0\,.
\end{array}\right.
$$

We can provide an expression for these functions in terms of the quantities introduced above:
$$
\a(c)= \left\{\begin{array}{lll}
0 && {\rm if}\; -\frac{4}{\pi}\leq c \leq \frac{4}{\pi}\\
E(|c|) && {\rm if}\; |c| > \frac{4}{\pi}
\end{array}\right.
$$
and 

$$
\b(h)= \left\{\begin{array}{lll}
0 && {\rm if}\; h=0\\
c(E(\frac{1}{|h|}))|h| - E(\frac{1}{|h|}) && {\rm if}\; h\neq 0\,.
\end{array}\right.
$$

\begin{figure} [h!]
\begin{center}
\includegraphics[scale=0.1]{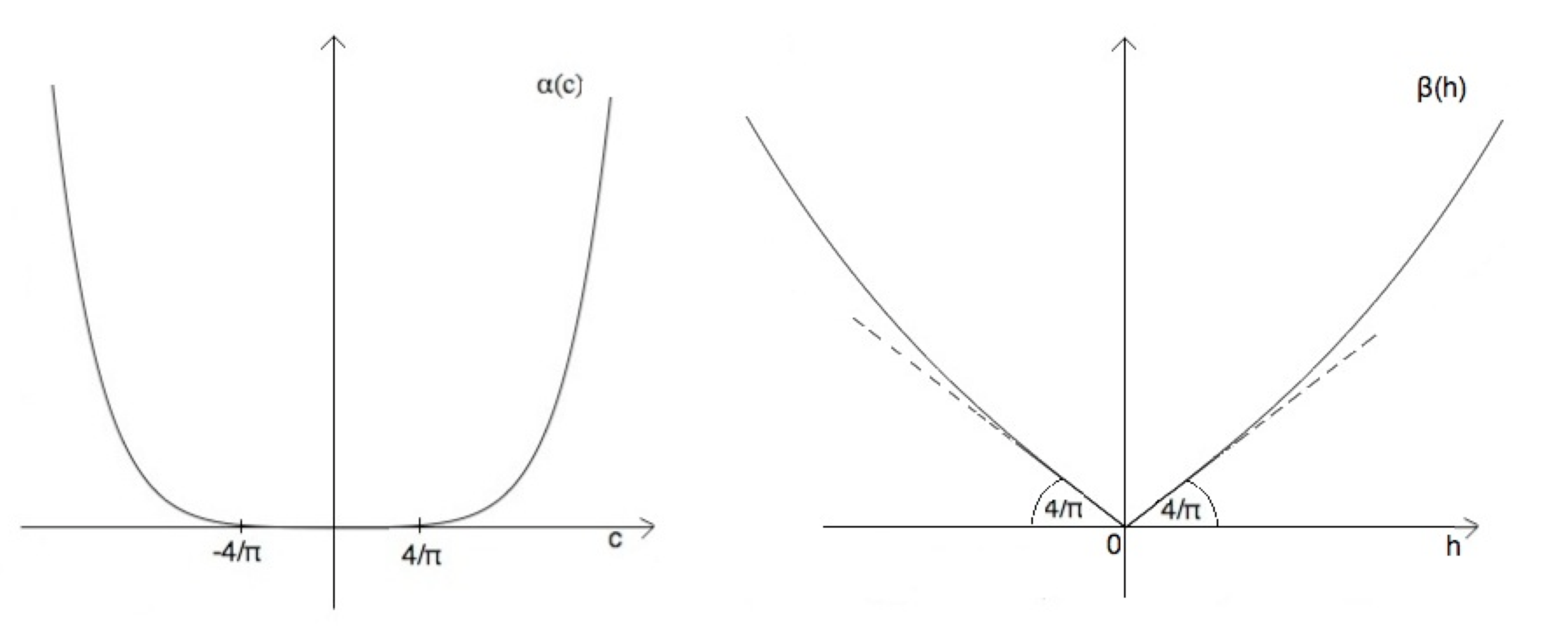}
\caption{Sketch of the graphs of the $\a$ and $\beta$ functions of the simple pendulum.}
\end{center}
\end{figure}

Observe that the $\a$ function is $C^1$. In fact, the only problem might be at $c=\pm\frac{4}{\pi}$, but also there it is differentiable, with derivative $0$. If it were not differentiable, then there would exist a subderivative $h\neq 0$ and consequently $\widetilde{\cM}^h \subseteq \widetilde{\cM}_{\pm\frac{4}{\pi}}$, which is absurd since the set on the right-hand side consists of a single point. However, $\a$ is not strictly convex, since there is a {\it flat} piece on which it is zero.\\
As far as $\beta$ is concerned, it is strictly convex (as a consequence of $\a$ being $C^1$), but it is differentiable everywhere except at the origin. At the origin, in fact, there is a corner and the set of subderivatives (\ie the slopes of tangent lines) is given by $\dpr \beta(0)=[-\frac{4}{\pi}, \frac{4}{\pi}]$ (this is related to the fact that $\a$ has a flat on this interval).\\

\vspace{20 pt}

\noindent{\bf {\sc 4.C - Holonomic measures  and generic properties}}\\
\addcontentsline{toc}{subsection}{\hspace{15 pt} 4.C: Holonomic measures and generic properties}

In this addendum we would like to stress that using the above approach the minimizing measures
are obtained through a variational principle over the set of invariant probability measures. Because of the request of  ``invariance'', this set clearly depends on the Lagrangian that one is considering. Moreover, it is somehow unnatural for variational problems to ask ``a-priori invariance''. What generally happens, in fact, is  that invariance is obtained as a byproduct of the minimization process carried out. \\
An alternative approach, slightly different under this respect, was due to Ricardo Ma\~n\'e \cite{Maneminmeasure} (see also \cite{ManeGeneric}). This deals with the bigger set of {\it holonomic measures} (or {\it closed measures}, see Remark \ref{closedequalholonomic}) and prove extremely advantageous when dealing with different Lagrangians at the same time. In this addendum we want to sketch the basic ideas behind it.

Let $C^0_{\ell}$ be the set of continuous functions $f: \rT M \rightarrow \R$ growing (fiberwise) at most linearly, \ie
$$
\|f\|_{\ell}:= \sup_{(x,v)\in \rT M} \frac{f(x,v)}{1+ \|v\|} < + \infty\,,
$$
and let ${\calM}^{\ell}_M$ be the set of probability measures on the Borel $\sigma$-algebra of $\rT M$ such that
$\int_{\rT M} \|v\|\, d\m < \infty\,,$ endowed with the unique metrizable topology given by:
$$ 
\m_n \longrightarrow\m \qquad \Longleftrightarrow \qquad \int_{\rT M} f(x,v)\,d\m_n \longrightarrow \int_{\rT M} f(x,v)\,d\m \quad \forall\,f\in C^0_{\ell}\,.
$$
Let $(C^0_{\ell})^*$ be the dual of $C^0_{\ell}$. Then ${\calM}^{\ell}_M$ can be naturally embedded in $(C^0_{\ell})^*$ and its topology
coincides with that induced by the weak$^*$ topology on $(C^0_{\ell})^*$. One can show that this topology is metrizable and a metric is, for instance:
$$
d(\m_1,\m_2)= \left| \int_{\rT M}\|v\|\,d\m_1 - \int_{\rT M}\|v\|\,d\m_2 \right| + \sum_{n} \frac{1}{2^n c_n}
\left| \int_{\rT M}\varphi_n \,d\m_1 - \int_{\rT M} \varphi_n\,d\m_2 \right|\,,
$$
where $\{\varphi_n\}_n$ is a sequence of functions with compact support on $C^0_{\ell}$, which is dense on $C^0_{\ell}$ (in the topology of uniform convergence on compact subsets of ${\rT M}$) and $c_n:=\sup_{\rT M} \left| \varphi_n (x,v)\right|$.
The space of probability measures that we shall be considering is a closed subset of ${\calM}^{\ell}_M$ (endowed with the induced topology), which is defined as follows. If $\g: [0,T] \rightarrow M$ is a closed absolutely continuous curve, let $\m_{\g}$ be such that
$$
\int_{\rT M} f(x,v) \,d\m_{\g} = \frac{1}{T} \int_0^T f(\g(t),\dot{\g}(t))\,dt \qquad \forall\,f\in C^0_{\ell}\,.
$$
Observe that $\m_{\g}\in {\calM}^{\ell}_M$ because if $\g$ is absolutely continuous then $\int |\dot{\g}(t)|\,dt < + \infty$. Let $\cC(M)$ be
the set of such $\m_{\g}$'s and $\overline{\cC(M)}$ its closure in  ${\calM}^{\ell}_M$. This set is convex and it is called the set of  {\it holonomic measures} on $M$.

\noindent One can check that the following properties are satisfied (see \cite{ManeGeneric}):

\begin{itemize}
\item[i)] $\calM (L) \subseteq \overline{\cC(M)} \subseteq {\calM}^{\ell}_M$. In particular, for every Tonelli Lagrangian $L$ on $\rT M$, all probabilities measures $\m$ that are invariant with respect to the Euler-Lagrange flow and such that $\int_{\rT M} L\, d\m < + \infty$, are contained in $\overline{\cC(M)}$.
\item[ii)] To any given probability  $\m~\in~\cC(M)$, one can  associate a rotation vector $\rho(\m)~\in~{\rH_1(M;\R)}$. This map extends continuously to a map $$\rho:~\overline{\cC(M)}~\longrightarrow~\!\!\rH_1(M;\R)$$ and this extension is surjective.
\item[iii)] For each $C\in \R$ the set $\left\{ \m\in \overline{\cC(M)}:\; A_{L}(\m)\leq C  \right\}$ is compact.
\item[iv)] If a measure $\m \in \overline{\cC(M)}$ satisfies
$$ A_L(\m) = \min \left\{ A_{L}(\n): \; \n\in \overline{\cC(M)}  \right\},$$
then $\m \in \calM (L)$ (and in particular it is invariant). Observe that the existence of probabilities attaining the minimum follows from iii).
\end{itemize}

In view of these properties, it is clear that the corresponding minimizing problem, although on a bigger space of measures, will lead to the same results as before and the same definition of Mather sets. 

\begin{Rem}\label{closedequalholonomic}
One can define the set of {\it closed measures} on $\rT M$ as:
$$
{\mathcal K}(\rT M):= \left\{\mu\in {\calM}^{\ell}_M\; \mbox{such that}
 \int df(x)\cdot v \, d\m = 0 \;\forall\,f\in C^1(M)\right\}\,.
$$
It is easy to verify that holonomic measures satisfy Proposition \ref{closedmeasure} and therefore
$\overline{\cC(M)} \subseteq {\mathcal K}(\rT M)$. It is definitely less trivial to prove that indeed these two sets coincide: $\overline{\cC(M)} \subseteq {\mathcal K}(\rT M)$.
Although this was originally noticed by John Mather, as far as I know it has never been published by the author. A proof of this result may be found in \cite{Bernardsuperposition}.\\
\end{Rem}

As we have already pointed out, this different approach is more suitable for working with different Lagrangians, for instance if one wants to study properties of family of Lagrangians or want to do some sort of perturbation theory.\\
Using these ideas, Ricardo Ma\~n\'e \cite{ManeGeneric} showed that one can prove much stronger results if, instead of considering ALL Lagrangians, one considers  {\it generic Lagrangians}.

\begin{Def}
A property $P$ is said to be {\it generic} ({\it in the sense of Ma\~n\'e}) for a Lagrangian $L$ if there exists a  {residual} set (\ie dense $G_{\d}$ set)  $\cS_L \subseteq C^2(M)$, such that  if $U\in \cS_L$ then $L+U$ satisfies the property P.
\end{Def}

\begin{Teo}[{\bf Ma\~n\'e}]
For a fixed $c\in \rH^1(M;\R)$, having a unique  $c$-action minimizing measure is a generic property in the sense of Ma\~n\'e. In other words, for any Tonelli Lagrangian $L$, there exists a residual subset 
$\cS_L \subseteq C^2(M)$, such that   for each $U\in \cS_L$, the Lagrangian $L+U$ has a unique $c$-action minimizing measure.
\end{Teo}

This result has been recently improved by Patrick Bernard and Gonzalo Contreras 
\cite{BernardContreras}:

\begin{Teo}[{\bf Bernard-Contreras}]
The following property is generic in the sense of Ma\~n\'e:
for all $c\in \rH^1(M;\R)$, there are at most $1 + \dim \rH^1(M;\R)$ ergodic  $c$-action minimizing measures.
\end{Teo}

%%%%%%%%%%%%%%%%%%%%%%%%%%%%%%%%%%%%%%%%%%%%%%%%%%%%%%%%%%%%%%%%%%%%%%%%%%%%%%%%%%%%%%%%%%%%%%%%%%%%%%%%%%%%%%%%%%%%%%%%%%%%%%%%%%%%%%%%%%%%%%%%%%%%%%%%%

%%%%%%%%%%%%%%%%%%%%%Section 1.3 %%%%%%%%%%%%%%%%%%%%%%%%%%%%%%%%%%%%%%%%%%%%%
\section{Global action minimizing curves: Aubry and Ma\~n\'e sets}\label{sec1.3}

In the previous section we have described the construction and the main properties of the Mather sets. One of the main limitations of these sets is that, being the support of invariant probability measures, they are {\it recurrent} under the flow (Poincar\'e recurrence theorem), \ie  each orbit after a sufficiently long time (and therefore infinitely many often) will return arbitrarily close to its initial point.
This property excludes many interesting invariant sets, which are somehow ``invisible'' to such a construction; for instance, think about the stable and unstable manifolds of some hyperbolic invariant set, or about heteroclinic and homoclinic orbits between invariant sets.\\
In this section we shall construct other  (possibly) ``larger'' compact invariant sets  and discuss their significance for the dynamics: the {\it Aubry sets} and the {\it Ma\~n\'e  sets}. The key idea is the same as  we have already explored in Section \ref{cartoonexample}: instead of considering action minimizing invariant probability measures, one can look at {\it action-minimizing curves} for some modified Lagrangian. We showed in Section \ref{cartoonexample} (see Remarks
\ref{remmaneset} and \ref{remregularminim}) that orbits on KAM tori could be characterized in terms of this property. In this section we shall mimic that construction in the general case of a Tonelli Lagrangians.

In the light of Lemma \ref{lemmamodifiedLagrangians} and the discussion in Section \ref{sec1.2}, let us fix a cohomology class $c\in \rH^1(M;\R)$ and choose a smooth $1$-form $\eta$ on $M$ that represents $c$.
As we have already pointed out in Section \ref{sec1.1}, there is a close relation between solutions of the Euler-Lagrange flow
and extremals of the action functional $\A_{L_{\eta}}$ for the fixed end-point problem (which are the same as the extremals of $\A_{L}$). In general, these extremals are not minima (they are local minima only if the time length is very short \cite[Section 3.6]{Fathibook}). One could wonder if such minima exist, namely if for any given end-points $x,y\in M$ and any given positive time $T$, there exists a minimizing curve connecting  $x$ to $y$ in time $T$. From what already said, this curve will correspond to an orbit for the Euler-Lagrange flow.
Under our hypothesis on the Lagrangian, the answer to this question turns out to be affirmative. This is a classical result in calculus of variations, known as {\it Tonelli Theorem}. %, that has been reproven - in the setting of Tonelli Lagrangians - by John Mather \cite{Mather91}.

\begin{Teo}[{\bf Tonelli Theorem, \cite{Mather91}}]\label{TonelliTheorem}
Let $M$ be a compact manifold and $L$ a Tonelli Lagrangian on $\rT M$.
For all $a<b\in\R $ and $x,y\in M$, there exists, in the set of absolutely continuous curves $\g:[a,b]\longrightarrow M$ such that $\g(a)=x$ and $\g(b)=y$,   a curve that minimizes the action $\A_{L_{\eta}}(\g) = \int_a^b {L_{\eta}}(\g(t),\dot{\g}(t))\,dt$.
\end{Teo}

\begin{Rem}
({\it i}) A curve minimizing $\A_{L{\eta}}(\g) = \int_a^b L_{\eta}(\g(t),\dot{\g}(t))\,dt$ subject to the fixed end-point condition $\g(a)=x$ and $\g(b)=y$, is called a $c$-{\it Tonelli minimizer}. Recall that such minimizers do only  depend on $c$ and not on the chosen representative $\eta$. In fact,  adding an exact $1$-form $df$ to $L$  will contribute with a constant term $f(y)-f(x)$, that does not play any role in selecting the minimizers (see also  Proposition \ref{closedmeasure}).\\
({\it ii}) As Ma\~n\'e pointed out in \cite{ManeI}, for these minimizers to exist it is not necessary to assume the compactness of $M$: the superlinear growth condition with respect to some complete Riemannian metric on $M$ is enough.\\
({\it iii}) A Tonelli minimizer which is $C^1$ is in fact $C^r$ (if the Lagrangian $L$ is $C^r$) and satisfies the Euler-Lagrange equation; this follows from the usual elementary arguments in calculus of variations, together with Caratheodory's remark on differentiability. In the autonomous case, Tonelli minimizers will be always $C^1$. In the non-autonomous time-periodic case (Tonelli Theorem holds also in this case \cite{Mather91}), as already remarked in remark \ref{remcompleteness}, one needs to require that the Euler-Lagrange flow is also complete.
\end{Rem}

We shall sketch here the proof of Tonelli Theorem. We refer the reader to \cite[Appendix 1]{Mather91} for more details.  A new shorter proof of this theorem can be found in \cite{Bernardsuperposition}.

\begin{Proof}[{\bf Tonelli Theorem}]
The proof of this theorem follows from the following result (by $C^{ac}$ we denote the set of absolutely continuous functions).\\

\noindent{\it Lemma.} {\it Let $K\in \R$. The set $S_K:=\{\g\in C^{ac}([a,b], M):\; A_{L_{\eta}}(\g)\leq K \}$ is compact in the $C^0$ topology.}\\

\noindent First let us see how to obtain Tonelli Theorem from this Lemma. Let $k_0:= \inf\{A_{L_{\eta}}(\g):\, \g\in C^{ac}([a,b], M)\}$. Observe that $k_0>-\infty$ since $L_{\eta}$ is bounded from below. Therefore, for any $K>k_0$, $S_{K}\neq \emptyset$ and it is compact because of the Lemma. Moreover, $S_{K} \subseteq S_{K'}$ if $K\leq K'$. Hence:
$$
\bigcap_{K>k_0} S_K \neq \emptyset
$$
and any element in this intersection is a $c$-Tonelli minimizer.\\

The proof of the above lemma consists of several steps (see \cite[Appendix 1]{Mather91} for the missing details). 
\begin{itemize}
\item[-] The first step is the observation that the family of curves in $S_k$ is absolutely equicontinuous, \ie for every $\e>0$ there exist $\d>0$ such that if $a\leq a_0 <b_0 \leq a_1 <b_1 \leq \ldots \leq a_n<b_n \leq b$ and $\sum_{i=0}^n (b_i-a_i) <\delta$, then
$\sum_{i=0}^n d(\g(a_i),\g(b_i)) <\e$. For this, one needs to use the superlinearity of $L$.
\item[-] Now, we can apply Ascoli-Arzel\`a theorem to deduce that every sequence $\{\g_n\}_n$ in $S_K$ has a convergent subsequence, with respect to the $C^0$ topology. Moreover, it follows easily from the definition of absolutely equicontinuity, that the limit of any convergent subsequence must be also absolutely equicontinuous.
\item[-] The last, and more involved, step is to show that if $\g$ is the limit of a sequence $\{\g_n\}_n$ in $S_K$, then $\g\in S_k$, namely $A_{L_{\eta}}(\g) \leq K$. See \cite[pages 199-201]{Mather91}.
\end{itemize}
These three steps conclude the proof of the lemma.
\end{Proof}

\begin{Rem}
Observe that the lemma used in the proof of Tonelli theorem is a sort of semicontinuity result for the Lagrangian action (compare with Proposition \ref{lowsemicactionmeasures}). In fact, it implies that if
$\{\g_n\}_n$ is a sequence in $C^{ac}([a,b],M)$ which converges to $\g$ in the $C^0$ topology, then $\g\in C^{ac}([a,b],M)$ and
$A_{L_{\eta}}(\g) \leq \liminf_{n\to +\infty} A_{L_{\eta}}(\g_n)$.
\end{Rem}

In the following we shall be interested in particular Tonelli minimizers that are defined for all times and whose action is minimal with respect to any given time length. We shall see that these curves present a very rich structure.

\begin{Def}[{\bf c-minimizers}]
An absolutely continuous  curve $\g: \R \longrightarrow M$ is a $c$-{\rm(}global{\rm)} minimizer for $L$, if for any given $a<b \in \R$
$$
A_{L_{\eta}}(\g\big| [a,b]) = \min A_{L_{\eta}}(\sigma)
$$
where the miminimum is taken over all $\sigma: [a,b] \rightarrow M$ such that $\s(a)=\g(a)$ and $\s(b)=\g(b)$.
\end{Def}

We have already seen in Section \ref{cartoonexample} that one can give an a-priori stronger notion of minimizer, asking that the minimum is realized amongst all curves connecting the two end-points, independently of their time length. 

\begin{Def}[{\bf c-time free minimizers}]
An absolutely continuous  curve $\g: \R \longrightarrow M$ is a $c$-time free minimizer  for $L$, if for any given $a<b \in \R$
$$
A_{L_{\eta}}(\g\big| [a,b]) = \min A_{L_{\eta}}(\sigma)
$$
where the miminimum is taken over all $\sigma: [a',b'] \rightarrow M$ such that $\s(a')=\g(a)$ and $\s(b')=\g(b)$.
\end{Def}

\begin{Rem}
({\it i}) We have proved in section \ref{cartoonexample} that orbits on KAM tori satisfy this stronger condition, modulo adding a constant to the Lagrangian (see Remark \ref{remmaneset}). In fact, it is quite easy to see that this condition is ``sensitive'' to the addition of constants to the Lagrangian (although this is something totally irrelevant for being  or not a $c$-minimizer). 
For example, suppose that  $\int_0^T L_{\eta} (\g,\dot{\g}) dt < \int_0^{T'} L_{\eta} (\s,\dot{\s}) dt$ with $T'<T$ and let $k$ be a constant such that
$k>\frac{1}{T'-T}\left( \int_0^{T'} L_{\eta} (\s,\dot{\s}) dt - \int_0^T L_{\eta} (\g,\dot{\g}) dt \right)$. Then, adding $k$ to the Lagrangian, we even reverse the inequality:
\beqano
A_{L_{\eta}+k}(\g) &=& \int_0^T L_{\eta} (\g,\dot{\g}) dt + kT =  \int_0^T L_{\eta} (\g,\dot{\g}) dt + k(T-T') + kT' \geq \\
&\geq& \int_0^T L_{\eta} (\g,\dot{\g}) dt + \int_0^{T'} L_{\eta} (\s,\dot{\s}) dt - \int_0^T L_{\eta} (\g,\dot{\g}) dt + kT' = \\
&=& \int_0^{T'} L_{\eta} (\s,\dot{\s}) dt  + kT' = A_{L_{\eta}+k}(\s)\,.
\eeqano
\noindent ({\it ii}) Obviously a $c$-time free minimizer is also a $c$-minimizer. Fathi \cite{Fathiconvsemigroup} proved that in the autonomous case, modulo adding the (unique) right constant, these two notions of minimizers indeed coincide: {\it if $\g$ be a $c$-minimizer for $L$, then, $\g$ is 
$c$-time free minimizer for $L+\a(c)$, where $\a$ is the $\a$-function associated to $L$}. 
We proved this in the special setting of Section \ref{cartoonexample} (see Proposition \ref{proporbittori}).
We shall discuss the general case in Section \ref{sec1.4}.\\
({\it iii}) The equivalence between these two notions of minimizers is not true anymore  when we consider time-periodic Tonelli Lagrangians (see  \cite{MatherFathi}). Tonelli Lagrangians for which this equivalence result holds are called {\it regular}. Patrick Bernard \cite{Bernardconnecting} showed that under suitable assumptions on the Mather set it is possible to prove that the Lagrangian is regular. For instance, if the Mather set $\tilde{\cM}_c$ is union of $1$-periodic orbits, then $L_{\eta}$ is regular.  This problem turned out to be strictly related to the convergence of the so-called {\it Lax-Oleinik semigroup} (see \cite{Fathibook} and Section \ref{sec1.4} for its definition).
\end{Rem}

Now we shall study the existence and the properties of   $c$-minimizers and $c$-time free minimizers.
 %(and therefore $c$-minimizers), trying to understand what is the ``right constant'' that we need to consider and why.
 There are two equivalent approaches that one can pursue: one is essentially due to Ricardo Ma\~n\'e  \cite{ManeI, ManeII}, while the other has been developed by John Mather \cite{Mather93}. In the following, in order to keep the analogy with the cartoon example discussed in Section \ref{cartoonexample}, we prefer following the first of these two approaches. We shall discuss Mather's approach in an Addendum 5C at the end of this section.\\
 
% In the following subsections we shall describe both of them and prove their equivalence. In subsection \ref{subsectionAubMane} we shall define the Aubry and Ma\~n\'e sets and discuss their properties.\\

%\subsection{Ma\~n\'e's approach: Ma\~n\'e potential, semi-static and static curves}

Given any $x,y\in M$ and T>0, let us denote by $C_T(x,y)$ the set of absolutely continuous curves 
$\g:[0,T]\longrightarrow M$ such that $\g(0)=x$ and $\g(T)=y$. Tonelli Theorem implies that there exists
$\g_{min} \in C_T(x,y)$ realizing the minimum, \ie $A_{\Le}(\g_{min}) = \min_{\g\in C_T(x,y)} A_{\Le}(\g)$.
Our goal here is to study the existence of $c$-time free minimizers for $L$. For, let us fix $k\in \R$ and consider the following quantity:
$$
\Fk(x,y) = \inf_{T>0} \min_{\g\in C_T(x,y)} A_{\Le+k}(\g) \in \R\cup \{-\infty\}\,.
$$

\noindent This quantity is commonly called {\it Ma\~n\'e potential} (compare with (\ref{Manepotentialontorus}) in Section \ref{cartoonexample}). First of all we would like to understand when it is finite and what are its properties. Let us introduce what is called {\it Ma\~n\'e critical value}.

\begin{Def}[{\bf Ma\~n\'e critical value}]
\beqano
c(\Le) &:=& \sup \{k\in \R:\, \exists \; \mbox{a closed curve}\;\g\; \mbox{s.t.}\; A_{\Le+k}(\g)<0\} =\\
&=& \inf \{k\in \R:\, \forall \; \mbox{closed curves}\;\g\; \mbox{s.t.}\; A_{\Le+k}(\g)\geq 0\}.
\eeqano
\end{Def}

\begin{Rem}
({\it i}) It is easy to check that $c(\Le)<\infty$. In fact, since $L$ is superlinear, there exists a sufficiently large $k$ such that $L+k\geq 0$ everywhere.\\
({\it ii}) Moreover $c(\Le)$ only depends on $c=[\eta]$ and not on the chosen representative. In fact, it is sufficient to notice that the integral of exact 1-forms along closed curves is zero. We shall see in the following that this ``critical value'' is something that we have already met before: $c(L_{\eta})=\a(c)$, where $\a$ is Mather's $\a$-function associated to $L$. In Section \ref{sec1.4} we shall also point out its relation  to  viscosity solutions and subsolutions of Hamilton-Jacobi equation and the critical value introduced by Lions, Papanicolau and Varadhan in \cite{LPV}.
\end{Rem}

\begin{Prop}[{See also \cite{ManeII}}]\label{propertymanepot}
\beqano
&(1)& \quad \forall\; k\in\R: \hspace{1cm} \Fk(x,y)\leq \Fk(x,z) + \Fk(z,y) \quad \forall\; x,y,z\in M.\\
&(2)& \quad \mbox{If}\; k<c(\Le): \quad \Fk(x,y)\equiv -\infty \hspace{3 cm} \forall\; x,y\in M.\\
&& \quad \mbox{If}\; k\geq c(\Le): \quad \Fk(x,y)\in \R \hspace{3.4 cm} \forall\, x,y\in M.\\
&(3)& \quad \mbox{If}\; k\geq c(\Le): \quad \Fk: M\times M \longrightarrow \R \mbox{ is Lipschitz}.\\
&(4)& \quad \mbox{If}\; k\geq c(\Le): \quad \Fk(x,x)\equiv 0 \hspace{3.47 cm} \forall\; x\in M.\\
&(5)& \quad \mbox{If}\; k\geq c(\Le): \quad \Fk(x,y)+\Fk(y,x)\geq 0 \hspace{1.6 cm} \forall\; x,y\in M.\\
&& \quad \mbox{If}\; k> c(\Le): \quad \Fk(x,y)+\Fk(y,x)> 0 \hspace{1.6 cm} \forall\; x\neq y\in M.\\
\eeqano
\end{Prop}

\begin{Proof}
$(1)$ First of all observe that this inequality makes sense also if $\Fk(x,y)=-\infty$ for some $x,y\in M$. 
Let $\g_1 \in C_T(x,z)$ and $\g_2\in C_{T'}(z,y)$ and consider the new curve obtained by joining them:
$\g_1*\g_2 \in C_{T+T'}(x,y)$. Since the action is linear, it follows that
$A_{\Le}(\g_1*\g_2 ) = A_{\Le}(\g_1)+A_{\Le}(\g_2 )$ and therefore
$$
\Fk(x,y) \leq A_{\Le}(\g_1)+A_{\Le}(\g_2 ).
$$
It is now sufficient to take the infimum over all possible $(\g_1, T)$  and $(\g_2,T')$ to conclude that
$$\Fk(x,y) \leq \Fk(x,z)+\Fk(z,y).$$

(2) We shall first prove that  if $k<c(\Le)$ there exists $x_0\in M$ such that $\Fk(x_0,x_0)=-\infty$. In fact, from the definition of $c(\Le)$, we know that there exists $\g:[0,T]\longrightarrow M$ closed curve with $A_{\Le}(\g) <0$. Let us denote by $\g^n$ the $n$-time iteration of $\g$, \ie $\g^n:= \g*\ldots*\g$ ($n$ times). Then,
$$
\Fk(\g(0),\g(0)) \leq A_{\Le}(\g^n) = n A_{\Le}(\g) \stackrel{n\to \infty}{\longrightarrow} -\infty\,.
$$
Choose $x_0:=\g(0)$ (or any other point on $\g$). The first claim will now follow from $(1)$. In fact, if 
$x,y\in M$, then:
\beqano
\Fk(x,y) \leq \Fk(x,x_0) + \Fk(x_0,x_0) + \Fk(x_0,y) = -\infty\,.
\eeqano
As for the second claim, if $k\geq c(\Le)$ it follows from the definition of $c(\Le)$ that all closed curves have positive action and therefore $\Fk(x,x) \geq 0$ for all $x\in M$. But, if there existed $x,y \in M$ such that $\Fk(x,y)=-\infty$, then applying $(1)$, we would get a contradiction: 
$$
\Fk(x,x) \leq \Fk(x,y) +\Fk(y,x) = -\infty\,.
$$
Therefore, $\Fk(x,y) >-\infty$ for all $x,y\in M$ if $k\geq c(L)$.

(3) Let $k\geq c(\Le)$ and let $Q:= \max_{x\in M, \|v\|=1} L(x,v)$. For any $x,y\in M$ let us consider the unit speed geodesic connecting $x$ to $y$, $\g_{x,y}:[0,d(x,y)]\longrightarrow M$. Then:
$$
\Fk(x,y) \leq A_{\Le}(\g_{x,y}) \leq (Q+k)d(x,y)\,.
$$
Using $(1)$ we  can conclude:
\beqano
\Fk(x_2,y_2)-\Fk(x_1,y_1) &=& \Fk(x_2,x_1) + \Fk(x_1,y_1) + \Fk(y_1,y_2) - \Fk(x_1,y_1) \leq\\
&\leq& (Q+k) [d(x_1,x_2)+d(y_1,y_2)].
\eeqano

(4) It follows immediately from (3).

(5) The first part is a consequence of (1) and (4). Let us prove the second part. Assume that $k>c(\Le)$ and suppose by contradiction  that there exist $x\neq y\in M$ such that $\Fk(x,y)+\Fk(y,x)=0$. Consider  
$\g_n \in C_{T_n}(x,y)$ and $\s_n\in C_{S_n}(y,x)$ such that:
$$
\lim_{n\to +\infty} A_{\Le}(\g_n) = \Fk(x,y) \qquad {\rm and} \qquad \lim_{n\to +\infty} A_{\Le}(\s_n) = \Fk(y,x).
$$
Let us prove that $\underbar{T}:=\liminf_{n\rightarrow +\infty} T_n >0$ (it might be $+\infty$). Suppose by contradiction that $\liminf_{n\rightarrow +\infty} T_n =0$ and select a subsequence $\{\g_{n_k}\}$ such that 
$T_{n_k}) \to 0$. Using the superlinearity of $L$, we know that for each $A>0$ there exists $B(A)$ such that $L_{\eta}\geq A\|v\| - B(A)$; then:
\beqano
\Fk(x,y) &=& \lim_{n_k \to \infty} \int_0^{T_{n_k}} \Le(\g_{n_k}, \dot{\g}_{n_k}) dt + k T_{n_k} \geq \\
&\geq& \lim_{n_k \to \infty} \left[ A\int_0^{T_{n_k}} \|\dot{\g}_{n_k}\| dt + T_{n_k}(k-B(A))\right]  =\\
&=&A\, d(x,y)\,.
\eeqano
From the arbitrariness of $A$ it follows that $\Fk(x,y)=+\infty$, that is a contradiction. Therefore, $\liminf_{n\rightarrow +\infty} T_n >0$.\\
Analogously one can prove that $\underbar{S}:=\liminf_{n\rightarrow +\infty} S_n >0$ (it might be $+\infty$).\\
Choose now subsequences $\{\g_{n_k}\}$ and $\{\s_{m_k}\}$ such that $T_{n_k} \to \underbar{T}>0$ and
$S_{m_k} \to \underbar{S}>0$. Using the fact that $k>c(\Le)$ and the fact that $x\neq y\in M$ are such that $\Fk(x,y)+\Fk(y,x)=0$, we obtain a contradiction  to (4) (or to the definition of $c(\Le)$):
\beqano
\Fc(x,x) &\leq& \lim_{k\to +\infty} A_{\Le+c(\Le)} (\g_{n_k}*\s_{m_k}) \leq  \\
&\leq& \lim_{k\to +\infty} A_{\Le+k} (\g_{n_k}*\s_{m_k}) + \lim_{k\to +\infty}(c(\Le)-k)(T_{n_k}+S_{n_k}) =\\
&=& \Fk(x,y) + \Fk(y,x) + \lim_{k\to +\infty}(c(\Le)-k)(T_{n_k}+S_{n_k}) =\\
&=& (c(\Le)-k)(\underbar{T}+\underbar{S}) < 0 \quad({\mbox or}\; -\infty\,). 
\eeqano
\end{Proof}

\begin{Rem}\label{remarkcritvalue}
It follows from (2) that $c(L_{\eta})$ can be equivalently defined as:
\beqano
c(\Le) &:=& \inf \{k\in \R:\, \exists\; x,y\in M\; \mbox{ s.t. }\; \Fk(x,y) > -\infty \} =\\
          &=& \sup \{k\in \R:\, \exists\; x,y\in M\; \mbox{ s.t. }\; \Fk(x,y) = -\infty \}.
\eeqano
\end{Rem}

In terms of Ma\~n\'e potential, being a $c$-time free minimizer for $L+k$ can be rewritten as:
$$
\forall\; a<b \qquad  \int_a^b \Le(\g(t),\dot{\g}(t))\,dt + (b-a) k = \Fk(\g(a),\g(b))\,.
$$
Of course, it does not make any sense to consider $k<c(L_{\eta})$.  Let us consider $k\geq c(\Le)$.

\begin{Prop}
Let $k>c(\Le)$. For all $x,y\in M$ with $x\neq y$, there exists $T>0$ and $\g\in C_T(x,y)$ such that
$A_{\Le+k}(\g) = \Fk(x,y)$.
\end{Prop}

\begin{Proof}
Let define for $T>0$ the function $f(T):=\min_{\g\in C_T(x,y)}A_{\Le + k}(\g)$. 
This function is clearly continuous (for all $k\geq c(\Le)) $ and the following properties hold.
\begin{itemize}
\item $f(T)\to +\infty$ as $t\to 0^+$ (this is true for all $k\geq c(\Le)$). In fact, let $\g_T$ be the corresponding Tonelli minimizer connecting $x$ to $y$ in time $T$. Using the superlinearity of $L$, for each $A>0$ there exists $B=B(A)$ such that $L(x,v)\geq A\|v\|-B$ for all $(x,v)$. Then:
\beqano
f(T) &=& \min_{\g\in C_T(x,y)}A_{\Le + k}(\g) = A_{\Le + k}(\g_T) \geq \\
&\geq& A \int_0^T \|\dot{\g}_T\|\,dt + (k-B)T \geq\\
&=& A\,d(x,y) + (k-B)T  \stackrel{T\to 0^+}{\longrightarrow} A\,d(x,y)\,.
\eeqano
Since $A>0$ is arbitrary and $x\neq y$ then  we can conclude that
$f(T)\stackrel{T\to0^+}{\longrightarrow} + \infty$.
\item $f(T)\to +\infty$ as $t\to +\infty$ (this is true only for  $k> c(\Le)$). In fact:
\beqano
f(T) &=& \min_{\g\in C_T(x,y)}A_{\Le + k}(\g) = \\
&=& \min_{\g\in C_T(x,y)}A_{\Le + c(\Le)}(\g) + (k-c(\Le))T \geq \\
&\geq& \Fc(x,y) + (k-c(\Le))T \stackrel{T\to+\infty}{\longrightarrow} +\infty.
\eeqano
\end{itemize}
\end{Proof}

\begin{Rem}
Hence if $k>c(\Le)$,  $c$-time free minimizers are not so special, since there are time free minimizers for $L+k$ connecting any two given points, furthermore in a finite time. In the light of this (and other ``{\it a-fortiori}'' reasons), one should probably be more interested in studying the ``critical'' case $k=c(\Le)$, \ie $c$-time free minimizers for $L+c(\Le)$, that is for the least possible value of $k$  for which they can exist.
\end{Rem}

\begin{Def}[{\bf c semi-static curves}]We say that  $\g:\R\longrightarrow M$ is a $c$ semi-static curve for $L$ if:
$$
\int_a^b L_{\eta}(\g(t),\dot{\g}(t))\,dt + c(\Le)(b-a) = \Fc(\g(a),\g(b))\qquad \forall\; a<b. 
$$
\end{Def}

\begin{Rem}\label{someremarkssemistatic}
({\it i}) If $\g$ is $c$ semi-static for $L$, then it is a $c$-time free minimizers for $L+c(\Le)$ and consequently a $c$-global minimizer for $L$. Therefore, it corresponds to a solution of the Euler-Lagrange flow of $L$.\\
({\it ii}) We shall prove in the following that for autonomous Tonelli Lagrangians, the converse is true: each $c$-global minimizer of $L$ is indeed a $c$ semi-static curve of $L+c(\Le)=L+\a(c)$, where $\a$ denotes Mather's $\a$ function associated to $L$ (see  (\ref{defalfa}) in Section \ref{sec1.2}).\\
({\it iii}) In Section \ref{cartoonexample} we proved that orbits on a KAM torus of cohomology class $c$, are $c$-time free minimizers for $L+E_c$, where $E_c$ denoted the energy of the torus. It follows from  Proposition \ref{Ecmanecritical} and the definition of Ma\~n\'e critical value (see also Remark \ref{remarkcritvalue}), that in this case $E_c=c(L_c)$. Therefore we can conclude that orbits on KAM tori are $c$ semi-static.
In particular, we can restate  Remark \ref{remmaneset} ({\it iii}) saying that if we have a KAM torus $\TT$ of cohomology class $c$, then
$$\cL^{-1}(\TT)= \bigcup\{(\g(t),\dot{\g}(t)): \; \g\;\mbox{is c semi-static for }\; L \; \mbox{and}\;t\in\R\},$$
where $\cL$ denotes the Legendre transform given by $L$.
\end{Rem}

Inspired by the last remark, we can define the following set.

\begin{Def}[{\bf Ma\~n\'e set}]
The {\it Ma\~n\'e set} (with cohomology class $c$) is:
\begin{eqnarray} 
\widetilde{\cN}_{c} &=& \bigcup \left\{(\g(t),\dot{\g}(t)):\; \text{$\g$ is a
$c$ semi-static curve and}\;t\in\R\right\} \nonumber =\\
&=& \bigcup \left\{(\g(t),\dot{\g}(t)):\; \text{$\g$ is a
$c$-global minimizer and}\;t\in\R\right\}.\label{maneset}
\end{eqnarray}
\end{Def}

\begin{Rem}
The second equality in definition follows from Remark \ref{someremarkssemistatic} and will be proved later on. 
Observe that so far we have not proved that such semi-static curves exist nor that this set is non-empty. We shall prove it later,  deducing it - amongst other properties - from analogous results for another family of sets that we are about to define: the {\it Aubry sets}.
However, if such set is non-empty, it is clearly invariant (it is union of orbits) and also closed (the proof follows the same line as Tonelli Theorem).
\end{Rem}

Let us start by recalling what happened in the case of orbits on a KAM torus. We saw in Section \ref{cartoonexample} that these orbits were not only $c$-global minimizers (or $c$ semi-static), but they also satisfied a stronger property, stated in Proposition \ref{regularminim}. Roughly speaking, the action of $L_c+E_c$ on a piece of curve between two endpoints $x$ and $y$, was not only the minimal needed  to connect $x$ to $y$,  but it was also equal to {\it minus} the minimal action to connect $y$ back to $x$. We called a curve satisfying such a condition a {\it regular minimizer}. 
Recall, in fact, that, as it follows easily from Proposition \ref{propertymanepot} (5), for each $x,y\in M$ we have $\Fc(x,y) \geq -\Fc(y,x).$\\
Let us define such curves in the general case.

\begin{Def}[{\bf c static curves}] We say that  $\g:\R\longrightarrow M$ is a $c$ static curve for $L$ (or a ``c-regular minimizer'') if:
$$
\int_a^b L_{\eta}(\g(t),\dot{\g}(t))\,dt + c(\Le)(b-a) = - \Fc(\g(b),\g(a))\qquad \forall\; a<b. 
$$
\end{Def}

\begin{Rem}\label{someremarksstatic}
({\it i}) If $\g$ is $c$ static for $L$, then it is a $c$ semi-static  (and therefore, as already observed before, it corresponds to a solution of the Euler-Lagrange flow of $L$).  It just a consequence of the fact that 
$- \Fc(\g(b),\g(a)) \leq \Fc(\g(a),\g(b$)) and therefore if $\g$ is $c$-static then:
$$
\int_a^b L_{\eta}(\g(t),\dot{\g}(t))\,dt + c(\Le)(b-a) \leq \Fc(\g(a),\g(b))\qquad \forall\; a<b. 
$$
Since $\Fc(\g(a),\g(b))$ was defined as the minimum over all connecting curves, then equality must hold.\\
({\it ii}) Observe that the adjective {\it regular} in the alternative appelation (coined by John Mather) has no relation to the smoothness of the curve, since, like all solutions of the Euler-Lagrange flow, this curve will be as smooth as the Lagrangian.\\
({\it iii}) In Section \ref{cartoonexample} (Proposition \ref{regularminim} and Remark \ref{remregularminim}) we proved that orbits on a KAM torus of cohomology class $c$, were $c$-static curves. Recall in fact that, as remarked before, $E_c$, which in that case denoted the energy of the torus, coincides with $c(L_c)$ (use Proposition \ref{Ecmanecritical} and the definition of Ma\~n\'e critical value or Remark \ref{remarkcritvalue}).  In particular, we can restate  Remark \ref{remregularminim} ({\it ii}) saying that if we have a KAM torus $\TT$ of cohomology class $c$, then
$$\cL^{-1}(\TT)= \bigcup\{(\g(t),\dot{\g}(t)): \; \g\;\mbox{is c-static for }\; L \; \mbox{and}\;t\in\R\},$$
where $\cL$ denotes the Legendre transform given by $L$.\\
({\it iv}) The names {\it semi-static} and {\it static} are probably inspired by the fact that in the case of mechanical Lagrangians (see Section \ref{sec1.1}), $0$-static curves correspond to ``some'' fixed points of the flow (namely the minima of the potential), while $0$-semi static curves may possibly include also hetero/homoclinic connections amongst these fixed points (see Addendum 5B on the pendulum).
\end{Rem}

Inspired by the last remark, we define the following set.

\begin{Def}[{\bf Aubry set}]
The {\it Aubry set} (with cohomology class $c$) is:
\begin{eqnarray} 
\widetilde{\cA}_{c} &=& \bigcup \left\{(\g(t),\dot{\g}(t)):\; \text{$\g$ is a
$c$ static curve and}\;t\in\R\right\}.\label{Aubryset}
\end{eqnarray}
The projection on the base manifold $\cA_c = \pi \left(\widetilde{\cA}_c\right)
\subseteq M$ is called {\it projected Aubry set} (with cohomology class $c$).
\end{Def}

\begin{Rem}
Observe that so far we have not proved that this set is non-empty! We shall do it now. However, if such set is non-empty, it is clearly invariant (it is union of orbits) and also closed (the proof follows the same line as Tonelli Theorem). Moreover, $\widetilde{\cA}_c$ is clearly contained in $\widetilde{\cN}_c$ (Remark \ref{someremarksstatic} ({\it i})).
\end{Rem}

We summarise in this diagram what we are going to prove in the remaing of this section.

\begin{Teo}\footnote{This (unintentional?) typographical ``coincidence'' honoring Ricardo Ma\~n\'e was first pointed out by Albert Fathi.}\label{diagramma} {\rm (1), (2), (3), (4), (5)} and {\rm (6)} in the following diagram  are true.
$$\xymatrix@R=1.5cm{
{\widetilde{\cM}_c} \ar@<0.7ex>[d]^{\pi} \ar@{}[r]|{\stackrel{(1)}{\subseteq}} & \widetilde{\cA}_c
\ar@<-0.8ex>[d]_{\pi} \ar@{}[r]|{\stackrel{(2)}{\subseteq}}  & 
\widetilde{\cN}_c \ar@{}[r]|(.2){{\stackrel{(3)}{\subseteq}}} & {\widetilde{\cE}_c:=\{E(x,v)=\a(c) {\stackrel{{(4)}}{=}}c(\Le)\}} \ar@{}[r]|(.75)
{\stackrel{}{\subseteq}} & {\rm T}M \ar@{->}[d]^{\pi}
\\  
\cM_c \ar@<0.6ex>[u]|(.4){(5)}^{\small{(\pi|\widetilde{\cM}_c)^{-1}}} \ar@{}[r]|\subseteq & 
\ar@<-0.6ex>[u]|(.4){(6)}_{\small{(\pi|\widetilde{\cA}_c)^{-1}}}\cA_c \ar@{}[r]|\subseteq  &&&  \ar@{}[l]|(.25)
\subseteq M\\}
$$
\end{Teo}

\smallskip
\begin{Rem}
({\it i}) It follows from (1), (2) and the existence of $c$-action minimizing measures (Corollary \ref{existminmeas}) that the Aubry and Ma\~n\'e sets are non-empty. Therefore, there exist semi-static and static curves.\\
({\it ii}) Inclusion $(2)$, as already observed, follows obviously from the fact that static curves are also semi-static (Remark \ref{someremarksstatic} ({\it i})).\\
({\it iii}) The above inclusions (1) and (2) may not be strict (see Remark \ref{remarkexamples} in Addendum 5B).\\
({\it iv}) Since $\widetilde{\cA}_c$ and $\widetilde{\cN}_c$ are closed, it follow from (3) that they are compact.\\
({\it v}) The proof of (3) provides a proof of Carneiro' s theorem stated in Section \ref{sec1.2} (Theorem \ref{teocarneiro}).\\
({\it vi}) Properties (5) and (6) are what are generally called Mather's graph theorem(s). Namely, the Mather set and the Aubry set are contained in a Lipschitz graph over $M$. This is probably the most important property of these sets and it has many dynamical consequences. In some sense, this is why they can be thought as generalization of KAM tori (or Lagrangian graphs).\\
({\it vii}) The graph property does not hold in general for the Ma\~n\'e set (see Addendum 5B on the pendulum).
\end{Rem}

Let us start by proving that $c$ semi-static curves have energy $c(\Le)$. As we already remarked above, the original version of this theorem (with $c(\Le)$ replaced by $\a(c)$) is due to Carneio (see Theorem \ref{teocarneiro}). The proof presented here follows an idea of Ricardo Ma\~n\'e \cite{ManeI} (see also \cite[Theorem XI]{ManeII})

\begin{Prop}[Property (3)] \label{teocarneiropermane}$\widetilde{\cN}_c\subseteq
{\widetilde{\cE}_c:=\{E(x,v)=c(\Le)\}}$, \ie $c$ semi-static curves have energy equal to $c(\Le)$.
\end{Prop}

\begin{Proof}
Let $\g:\R \longrightarrow M$ be a $c$ semi-static curve, \ie for each $T>0$ we have
$A_{\Le+c(\Le)}(\g|[0,T]) = \Fc(\g(0),\g(T))$. \\
Let us fix $T>0$ and $\l>0$, and consider a time-reparameterization of $\g$, given by
$\g_\l: [0,T/\l]\longrightarrow M$, $t\mapsto \g(\l t)$. Observe that the end-points are not changed:
$\g_\l(0)=\g(0)$ and $\g_\l(T/\l)=\g(T)$, but only the time lenght. Let us consider the action of these curves as a function of $\l>0$:
\beqano
A(\l) &:=& A_{\Le + c(\Le)}(\g_\l) = \int_0^{T/\l} \Le (\g_\l(t), \dot{\g}_\l(t))\,dt + c(\Le)T/\l =\\
&=&  \int_0^{Tl} \Le (\g(\l t), \l \dot{\g}(\l t))\,dt + c(\Le)T/\l.
\eeqano
Since  $\g$ is $c$ semi-static (and therefore it is a $c$-time free minimizer for $\Le+c(\Le)$), then
$A$ has a minimum at $\l=1$. Therefore:
\beqano
0 = A'(1) &=& -T\left[ \Le(\g(T),\dot{\g}(T))+c(\Le)T\right] + \\
&&+\; \int_0^T\left(
\frac{\dpr L}{\dpr x} (\g(t),\dot{\g}(t)) \dot{\g}(t) t + 
\frac{\dpr L}{\dpr v} (\g(t),\dot{\g}(t)) \big[\dot{\g}(t)+ t \ddot{\g}(t)\big] \right)\, dt\,.
\eeqano
Integrating by parts (observe that 
%$\frac{d}{dt}L(\g(t),\dot{\g}(t)) =  \left( \frac{\dpr L}{\dpr x} (\g(t),\dot{\g}(t)) \dot{\g}(t) + \frac{\dpr L}{\dpr v} (\g(t),\dot{\g}(t)) \ddot{\g}(t) \right) t $) 
$\frac{d}{dt}L =  \left( \frac{\dpr L}{\dpr x} \dot{\g} + 
\frac{\dpr L}{\dpr v} \ddot{\g} \right) t $)
and recalling the definition of the energy 
$E(x,v) = \frac{\dpr L}{\dpr v}(x,v)v - L(x,v)$ and the fact that is preserved along the orbit
($E(\g(t),\dot{\g}(t)) = E(\g(0),\dot{\g}(0))$ for all $t\in \R$), we obtain:
\beqano
0 &=& A'(1) = \ldots =\\
&=& -T\left[ \Le(\g(T),\dot{\g}(T))+c(\Le)\right] + \Le(\g(t),\dot{\g}(t))\Big|_{0}^T -\\
&& - \int_0^T \left(
- L (\g(t),\dot{\g}(t)) +  
\frac{\dpr L}{\dpr v} (\g(t),\dot{\g}(t)) \dot{\g}(t)\right)\, dt = \\
&=& - T c(\Le) + \int_0^T E(\g(t),\dot{\g}(t))\,dt = \\
&=&  \left[ E(\g(0),\dot{\g}(0)) - c(\Le) \right] T\,.
\eeqano
 Hence, $E(\g(0),\dot{\g}(0)) = c(\Le)$.
 \end{Proof}

Let us now prove that Ma\~n\'e critical value coincides with Mather's $\a$-function. We shall follow the proof given in \cite[Theorem II]{ManeII}

\begin{Prop}[Property (4)]\label{proofproperty4}
$c(\Le)=\a(c)$, where $c=[\eta]$.
\end{Prop}

\begin{Proof} 
Suppose that $\m$ is an invariant ergodic probability measure.
If we fix a generic point $(x,v)$ in the support of $\m$, it follows from the ergodic theorem that there exists a sequence of  times ${T_n}\rightarrow +\infty$ such that $\Phi^L_{T_n}(x,v)\to (x,v)$ as $T_n\to+\infty$ and
$$\int \Le \,d\m = \lim_{n\to +\infty} \frac{1}{T_n}\int_0^{T_n} \Le(\Phi^L_t(x,v))\,dt. $$
For the sake of simplifying the notation, let us denote $(x_t,v_t):=\Phi^L_{t}(x,v)$.
Let $B:=\max\{ |\Le(x,v)|:\; \|v\|\leq1\}$ and for each $n$ denote by $\s_{n}:[0, d(x,x_{T_n})]\longrightarrow M$ the geodesic joining $x$ to $x_{T_n}$ and by $\g_{n}:[0, {T_n}]\longrightarrow M$ the projection of the orbit (\ie $\g_n(t)=x_t$). We have:
\beqano
&& \lim_{n\to \infty} \frac{1}{T_n} A_{\Le+k}(\g_n*\s_n) \; =\\
&& \quad =\; \lim_{n\to \infty} \left(\frac{1}{T_n} A_{\Le+k}(\g_n) + \frac{1}{T_n}\int_0^{d(x,x_{T_n})} \Le(\s_n(t),\dot{\s}_n(t)) + k\,dt \right) \leq \\
&& \quad =\; \leq\; \lim_{n\to \infty} \left( \frac{1}{T_n} A_{\Le+k}(\g_n) + \frac{1}{T_n} (B+k){\rm diam}(M)\right) \leq \\
&& \quad =\; {A_{\Le}(\m)} + k\,. 
\eeqano
Therefore, if $k<-{A_{\Le}(\m)}$ then $\F_k(x,x) = -\infty$. Hence, $k\leq c(\Le)$. It follows that:
\beqano
c(\Le) &\geq& \sup\{-A_{\Le}(\m),\; \m\in \calM_{{\rm erg}}(L)\}  \geq\\
&\geq& - \inf\{A_{\Le}(\m),\; \m\in \calM(L)\}= \a(c)\,.
\eeqano
Now we want to prove the reversed inequality. Let $k<c(\Le)$ and $x,y\in M$. Since $\F_k(x,y)=-\infty$, there exists a sequence of absolutely continuous curves $\g_n: [0,T_n]\longrightarrow M$ such that
$\g_n(0)=x$, $\g_n(T_n)=y$ (\ie $\g_n\in C_{T_n}(x,y)$) and
$$
\lim_{n\rightarrow \infty} A_{\Le+k}(\g_n)=-\infty.
$$
Moreover, since $L$ is bounded from below, we have that $T_n\to +\infty$. Let now $y_n$ be a Tonelli minimizer in $C_{T_n}(x,y)$. Then, consider the invariant probability measure evenly distributed along  $\g_n$. The family of these measures is pre-compact and we can extract a subsequence converging (in the weak$^*$ topology) to an invariant probability measure $\m$. In particular,
$$
A_{\Le}(\m) + k = \lim_{n\to\infty} \frac{1}{T_n} A_{\Le+k}(\g_n).
$$
But since $T_n>0$ and $\lim_{n\to\infty} A_{\Le+k}(\g_n) = \F_k(x,y)=-\infty$, then 
$A_{\Le}(\m) + k \leq 0$. Therefore, for any $k<c(\Le)$ we can find an invariant probability measure $\m$ such that $k\leq -A_{\Le}(\m)$. Therefore:
\beqano
c(\Le) &\leq& \sup\{-A_{\Le}(\m),\; \m\in \calM_{{\rm erg}}(L)\}  \geq\\
&\leq& - \inf\{A_{\Le}(\m),\; \m\in \calM(L)\}= \a(c)\,.
\eeqano
\end{Proof}

We shall provide two other alternative proofs of the above proposition. The first one uses the fact that the Mather set is included into the Ma\~n\'e set (but observe that in our proof of this fact - Proposition \ref{inclusionmatheraubry} - we use that $\a(c)=c(L_{\eta})$). The second one is shorter but it requires the use of holonomic measures (see Addendum 4C). See also \cite[Theorem II]{ManeII} for another proof.\\

\noindent {\bf Alternative Proof I.}
Since $\widetilde{\cN}_c$ is compact (it  follows from (3)) and invariant under the Euler-Lagrange flow $\Phi^L_t$, then there exists an invariant ergodic probability measure $\m$ supported in it (Kryloff and Bogoliouboff \cite{KB}, compare also with Proposition \ref{KriloffBog}). If we fix a generic point $(x,v)$ in the support of $\m$, it follows from the ergodic theorem that there exists a sequence of  times ${T_n}\rightarrow +\infty$ such that $\Phi^L_{T_n}(x,v)\to (x,v)$ as $T_n\to+\infty$ and
$$\int \Le \,d\m = \lim_{n\to +\infty} \frac{1}{T_n}\int_0^{T_n} \Le(\Phi^L_t(x,v))\,dt. $$
Then, using the definition of $\a(c)$  and the fact that orbits in the support of this measure are semi-static, we obtain:
\beqano
-\a(c) &\leq& \int \Le \,d\m = \lim_{n\to +\infty} \frac{1}{T_n}\int_0^{T_n} \Le(\Phi^L_t(x,v))\,dt =\\
&=& \lim_{n\to +\infty} \frac{1}{T_n}\int_0^{T_n} \Le(\Phi^L_t(x,v)) + c(\Le)\,dt - c(\Le) =\\
&=& \lim_{n\to +\infty} \frac{\Fc(x, \pi(\Phi^L_{T_n}(x,v)))}{T_n}   - c(\Le)  = - c(\Le)\,,
\eeqano
where in the last equality we used that $\Fc$ is bounded on $M\times M$ (being Lipschitz on a compact manifold). Therefore $\a(c) \geq c(\Le)$.\\
For the reversed inequality,  we shall use the following version of the Ergodic theorem (see for instance 
\cite[Lemma 2.1]{ManeGeneric} for a proof).
%\cite{ContrerasIturriaga}[Corollary 3-6.5] for a proof).\\

\begin{Lem}\label{lemmateoergomane} Let $(X,d)$ be a complete metric space and $(X,\cB, \n)$ a probability space. Let $f$ be an ergodic measure preserving map and $F: X\longrightarrow \R$ a $\n$-integrable function. Then, for $\n$-almost every $x\in X$ the following property holds:
$$
\forall\, \e>0 \quad \exists\,N>0: \quad d(f^N(x),x)<\e \quad\mbox{and}\quad
\left|\sum_{j=0}^{N-1}F(f^j(x)) - N\int Fd\n \right| <\e.
$$
\end{Lem}

Let us now see how to use this lemma for our purposes. Let $\m$ be a $c$-action minimizing ergodic measure, \ie
$$
\int (\Le(x,v) +\a(c))\,d\m =0.
$$
Applying the above Lemma with $F=\Le+\a(c)$ and $X=\rT M$, we obtain that there exists a $\m$-full measure set $A$ such that if $(x,v)\in A$, then there exists a sequence $T_n\rightarrow +\infty$ such that 
$$d((x,v),\Phi^L_{T_n}(x,v))\stackrel{n\to +\infty}{\longrightarrow} 0 \quad {\rm and} \quad
\int_0^{T_n} \left(L(\Phi^L_t(x,v))+\a(c)\right)\,dt \stackrel{n\to +\infty}{\longrightarrow} 0\,.  
$$
Then:
\beqano
&& \phi_{\eta,\a(c)}(x, \pi(\Phi^L_1(x,v))) + \phi_{\eta,\a(c)}(\pi(\Phi^L_1(x,v)),x) =\\ 
&&\quad =\; \lim_{n\to\infty} \left( \phi_{\eta,\a(c)}(x, \pi(\Phi^L_1(x,v))) + \phi_{\eta,\a(c)}(\pi(\Phi^L_{1}(x,v)),\pi(\Phi^L_{T_n}(x,v))) \right) \leq \\
&&\quad \leq \; \lim_{n\to \infty} \int_0^{T_n} \left(\Le(\Phi^L_t(x,v))+\a(c)\right)\,dt  = 0.
\eeqano
It follows from the second property in Proposition \ref{propertymanepot} (5) that $\a(c)\leq c(\Le)$ and this concludes the proof. \qed \\

\noindent {\bf Alternative Proof II.} Let $\g$ be any closed curve and let $\m_{\g}$ be the probability measure evenly distributed on it (see Addendum 4C). Let $k\geq c(\Le)$. It follows from the definition of $c(\Le)$ that $A_{\Le +k} (\m_{\g}) \geq 0$ and therefore $A_{\Le} (\m_{\g}) \geq -k$. It follows from the definition of holonomic measure then, that  for any $\m$ holonomic probability measure, we have $A_{\Le} (\m) \geq -k$. Taking the infimum over all holonomic measure and using the result mentioned in iv) in Addendum 4C, we can conclude that  $- \a(c) \geq -k$. Since this holds for all $k\geq c(\Le)$, we obtain: $\a(c)\leq c(\Le)$.\\
To prove the reversed inequality, observe that if $k<c(\Le)$, then there exists $\g$ closed curve such that $A_{\Le+k}(\g)<0$. Therefore, if $\m_{\g}$ is the associated holonomic measure we obtain:
$$
-\a(c) \leq A_{\Le}(\g) < -k\,.
$$
Since this holds for all $k<c(\Le)$, we conclude that $\a(c)\geq c(\Le)$. \qed \\

The proof of Property (1) is essentially similar to the (first) proof of Proposition \ref{proofproperty4}. We shall prove this more general result, due to Ricardo Ma\~n\'e \cite{ManeI} (see also \cite[Theorem IV]{ManeII}).

\begin{Prop}[Property (1)]\label{inclusionmatheraubry}
$\m \in \calM(L)$ is $c$-action minimizing if and only if ${\rm supp}\,\m \subseteq \widetilde{\cA}_c$. In particular, $\widetilde{\cM}_c \subseteq \widetilde{\cA}_c $.
\end{Prop}

\begin{Proof}
Since $\widetilde{\cA}_c$ is closed, it is sufficent to prove the results only for the ergodic measures.

$[\Longleftarrow]$ Let $\m\in \calM(L)$ be ergodic and suppose that  ${\rm supp}\,\m \subseteq \widetilde{\cA}_c$. Applying the ergodic theorem, we know that
for a $\m$-generic point $(x,v)$ in the support of $\m$, there exists a sequence of  times ${T_n}\rightarrow +\infty$ such that $\Phi^L_{T_n}(x,v)\to (x,v)$ as $T_n\to+\infty$ and
$$\int \Le \,d\m = \lim_{n\to +\infty} \frac{1}{T_n}\int_0^{T_n} \Le(\Phi^L_t(x,v))\,dt. $$
Then, using  the fact that orbits in the support of this measure are semi-static and that $\a(c)=c(\Le)$, we obtain:
\beqano
\int \Le + \a(c) \,d\m &=& \int \Le + c(\Le )\,d\m = \\
&=& \lim_{n\to +\infty} \frac{1}{T_n}\int_0^{T_n} \Le(\Phi^L_t(x,v)) + c(\Le)\,dt =\\
&=& \lim_{n\to +\infty} \frac{- \Fc(\pi(\Phi^L_{T_n}(x,v)),x)}{T_n}=0,
\eeqano
where in the last equality we used that $\Fc$ is bounded on $M\times M$ (being Lipschitz on a compact manifold). Therefore $\int \Le \,d\m \leq -\a(c)$ and from the definition of $\a(c)$, it follows that 
$\int \Le \,d\m = -\a(c)$, \ie $\m$ is $c$-action minimizing.

$[\Longrightarrow]$ We shall use the above mentioned version of the ergodic theorem (see Lemma \ref{lemmateoergomane}).

Let $\m$ be a $c$-action minimizing ergodic measure, \ie
$$
\int (\Le(x,v) +\a(c))\,d\m =0.
$$
Applying the above Lemma with $F=\Le+\a(c)$ and $X=\rT M$, we obtain that there exists a $\m$-full measure set $A$ such that if $(x,v)\in A$, then there exists a sequence $T_n\rightarrow +\infty$ such that 
$$d((x,v),\Phi^L_{T_n}(x,v))\stackrel{n\to +\infty}{\longrightarrow} 0 \quad {\rm and} \quad
\int_0^{T_n} \left(L(\Phi^L_t(x,v))+\a(c)\right)\,dt \stackrel{n\to +\infty}{\longrightarrow} 0\,.  
$$
Then, let $a>0$ and recall that $\a(c)=c(\Le)$:
\beqano
&& \Fc(x, \pi(\Phi^L_a(x,v))) + \Fc(\pi(\Phi^L_a(x,v)),x) =\\ 
&&\quad =\; \lim_{n\to\infty} \left( \Fc(x, \pi(\Phi^L_a(x,v))) + \Fc(\pi(\Phi^L_{a}(x,v)),\pi(\Phi^L_{T_n}(x,v))) \right) \leq \\
&&\quad \leq \; \lim_{n\to \infty} \int_0^{T_n} \left(\Le(\Phi^L_t(x,v))+c(\Le)\right)\,dt  = \\
&&\quad \leq \; \lim_{n\to \infty} \int_0^{T_n} \left(\Le(\Phi^L_t(x,v))+\a(c)\right)\,dt  = 0.
\eeqano
Recalling property (5) in Proposition  \ref{propertymanepot}, we can conclude that for any $a>0$:
$$\Fc(x, \pi(\Phi^L_a(x,v))) + \Fc(\pi(\Phi^L_a(x,v)),x) = 0$$
and therefore the orbit through the point $(x,v)$ is $c$-static. Since the points for which this reasoning can be applied are dense in the support of $\m$ and $\widetilde{\cA}_c$ is closed, then we prove the claim: ${\rm supp}\,\m \subseteq \widetilde{\cA}_c$.
\end{Proof}

\begin{Rem}
Looking at the proof of Proposition \ref{inclusionmatheraubry}, it is quite easy to see that we actually proved that: if $\m \in \calM(L)$ is such that ${\rm supp}\,\m \subseteq \widetilde{\cN}_c$, then $\m$ is
$c$-action minimizing (check that the same proof still works in this case). In fact, one can prove this stronger version:
\begin{Prop}\label{inclusionmathermane}
$\m \in \calM(L)$ is $c$-action minimizing if and only if ${\rm supp}\,\m \subseteq \widetilde{\cN}_c$.
\end{Prop}
We shall deduce such a proposition from the fact (to be proved later in Section \ref{sec1.4}, Proposition \ref{propnonwandering}) that the ``{\it non-wandering set}'' of the Ma\~n\'e' set is contained in the Aubry set (see also \cite[Theorem V.c]{ManeII}).
\end{Rem}

It is also quite easy to check that orbits in the Ma\~n\'e (resp. Aubry) set are asymptotic to the Mather set.

\begin{Prop}
If $\g:\R\longrightarrow M$ is a $c$ semi-static curve, then $$\liminf_{t\rightarrow \pm\infty} d(\g(t), \cM_c) = 0.$$
\end{Prop}

\begin{Proof}
Let $T\geq 1$ and consider the probability measure $\m_T$ evenly distributed along the piece of curve $\{(\g(t), \dot \g(t)):\; t\in [0,T]\}$ (for a definition, for instance, the reader may check the proof of Proposition \ref{KriloffBog}). The Lagrangian actions of these measures  are equi-bounded (we use here that the orbit is semi-static):
\beqano
A_{\Le}(\m_T) &=& \frac{1}{T}\int_0^T \Le(\g(t),\dot\g(t))\,dt = \frac{\Fc(\g(0),\g(T))}{T} - c(\Le) \leq \\
&\leq& \max_{M\times M} \Fc(x,y) - c(\Le) < \infty\,.
\eeqano
Therefore, this family of measures is pre-compact with respect to the weak$^*$ topology. Let us consider any converging subsequence $\m_{T_{k}}\to \m$, with $T_k\to +\infty$. Then, $\m$ is invariant (see, again, the proof of Proposition \ref{KriloffBog}) and:
\beqano
\int \Le d\m &=& \lim_{k\to\infty} \int \Le d\m_{k} = \lim_{k\to\infty} \frac{\Fc(\g(0),\g(T))}{T} - c(\Le) = \\
&=&  - c(\Le) = -\a(c).
\eeqano
Therefore, $\m$ is $c$-action minimizing.
\end{Proof}

Finally, we prove the most important result of this theory: the graph property of the Mather and Aubry sets, respectively, Property (5) and (6) in the diagram (see also Theorem \ref{Theograph}). Because of the inclusion proved in Proposition \ref{inclusionmatheraubry}, it is sufficient to prove it for the Aubry set.

\begin{Teo}[{\bf Mather's graph theorem}, \cite{Mather91}]\label{graphtheoremAubry} {\rm (Property (6))}
 $\pi|{\widetilde{\cA}_c}$ is an injective mapping of $\widetilde{\cA}_c$ into $M$, and its inverse $(\pi |\widetilde{\cA}_c)^{-1}: \cA_c  \longrightarrow \widetilde{\cA}_c$ is Lipschitz.
\end{Teo}

The proof of this theorem will be based on the following ``crossing'' Lemma, proved by Mather in \cite{Mather91} (to which we refer the reader for a complete proof).

\begin{Lem}[{\bf Mather's crossing lemma, \cite{Mather91}}]\label{crossinglemma}
Let $K>0$. There exist $\e,\d,\vartheta >0$, $C>0$ such that if $\a,\b:[-\e,\e]\longrightarrow M$ are solutions of the Euler-Lagrange equation with $\|(\a(0),\dot\a(0))\|\leq K$, $\|(\b(0),\dot\b(0))\|\leq K$  and
$$
d(\a(0),\b(0))\leq \delta \quad \mbox{and} \quad d\left((\a(0),\dot\a(0)),(\b(0),\dot\b(0))\right) > C d(\a(0),\b(0)),
$$
then there exist $C^1$ curves $a,b:[-\e,\e]\longrightarrow M$ with end-points 
$a(-\e)=\a(-\e),\; a(\e)=\b(\e)$ and $b(-\e)=\b(-\e),\; b(\e)=\a(\e)$ such that:
$$
A_{\Le}(\a)+A_{\Le}(\b) - A_{\Le}(a)-A_{\Le}(b) \geq \vartheta\, d\left((\a(0),\dot\a(0)),(\b(0),\dot\b(0))\right)^2 >0\,.
$$
\end{Lem}

\begin{center}
\setlength{\unitlength}{1 cm}
\begin{picture}(6,4)(-2,-2)
\qbezier(0,0)(0.8853,0.8853)
(2,0.9640)
\qbezier(0,0)(-0.8853,-0.8853)
(-2,-0.9640)
\qbezier(-2,0.9640)(-0.8853,0.8853)
(-0.05,0.05)
\qbezier(0.05,-0.05)(0.8853,-0.8853)
(2,-0.9640)
\qbezier[50](-2,0.9640)(0,0.6)
(2,0.9640)
\qbezier[50](-2,-0.9640)(0,-0.6)
(2,-0.9640)
\put(0,0.85){{\it a}}
\put(0,-1){{\it b}}
\put(-1,0.4){$\small{\a}$}
\put(-1,-0.5){$\small{\b}$}
\put(-2,0.9640){\circle*{0.1}}
\put(-2.5,1.2){$\small{\a(-\e)}$}
\put(2,0.9640){\circle*{0.1}}
\put(1.7,1.2){$\small{\b(\e)}$}
\put(-2,-0.9640){\circle*{0.1}}
\put(-2.5,-1.3){$\small{\b(-\e)}$}
\put(2,-0.9640){\circle*{0.1}}
\put(1.7,-1.3){$\small{\a(\e)}$}
\end{picture}
\end{center}

\begin{Proof}[{\bf Theorem \ref{graphtheoremAubry}}] We shall first use Lemma \ref{crossinglemma} to  prove the Lipschitz property and then the rest will be just a consequence of this property. Let choose $K:=\max_{\widetilde{\cA}_c}{\|(x,v)\|}$ (this is finite since $\widetilde{\cA}_c$ is compact) and let  $\e,\d,\vartheta, C$ be as in Lemma \ref{crossinglemma}. Then, we shall prove that: 
{\it  if $(x_1,v_1), (x_2,v_2) \in \widetilde{\cA}_c$ are such that $d(x_1,x_2)\leq \delta$, then
$d\left((x_1,v_1), (x_2,v_2)\right) \leq C d(x_1,x_2).$}\\
Suppose by contradiction that  $d\left((x_1,v_1), (x_2,v_2)\right) > C d(x_1,x_2)$ and consider the flow lines through these points, namely $\a(t):=\Phi^L_t(x_1,v_1)$ and $\b(t):=\Phi^L_t(x_2,v_2)$. They satisfy the hypothesis of Lemma \ref{crossinglemma} (with our choice of $K$) and hence we can deduce the existence of two other curves $a,b:[-\e,\e]\longrightarrow M$ with end-points 
$a(-\e)=\a(-\e),\; a(\e)=\b(\e)$ and $b(-\e)=\b(-\e),\; b(\e)=\a(\e)$, such that:
$$
A_{\Le}(a)+A_{\Le}(b) < A_{\Le}(\a)+A_{\Le}(\b).
$$
But then:
\beqano
&& \Fc(a(-\e),a(\e)) + \Fc(b(-\e),b(\e)) \; \leq\\
&& \quad\leq\;  A_{\Le}(a) + A_{\Le}(b) < A_{\Le}(\a) + A_{\Le}(\b) = \\
&& \quad=\; - \Fc(\a(\e),\a(-\e)) - \Fc(\b(\e),\b(-\e))\,.
\eeqano
where in the last equality we used that $\a$ and $\b$ are $c$-semistatic.
The above inequality and the triangle inequality for $\Fc$ (see Proposition \ref{propertymanepot} (1)) lead to a contradiction to Proposition \ref{propertymanepot} (5):
\beqano
&&\Fc(\a(-\e),\b(\e)) \;=\; \Fc(a(-\e),a(\e)) <\\
&&\quad<\; - \left(\Fc(b(-\e),b(\e)) + \Fc(\a(\e),\a(-\e)) + \Fc(\b(\e),\b(-\e)) \right) =\\
&&\quad=\; - \left(\Fc(\b(-\e),\a(\e)) + \Fc(\a(\e),\a(-\e)) + \Fc(\b(\e),\b(-\e)) \right) \leq \\
&&\quad \leq\; - \left(\Fc(\b(-\e),\a(-\e)) + \Fc(\b(\e),\b(-\e)) \right) \leq \\
&&\quad \leq\; - \Fc(\b(\e), \a(-\e))\,.
\eeqano
Therefore the inverse of the projection is locally Lipschitz and this concludes the proof.
\end{Proof}

\begin{Rem}
({\it i}) Actually it follows from the proof (in the choice of $K$) that the graphs of the Aubry sets (or Mather sets) corresponding to compact sets of cohomology classes, are equi-Lipschitz.\\
({\it ii}) An alternative proof of the graph property will be presented in Section \ref{sec1.4}, following Fathi's weak KAM theory.
\end{Rem}

One can show several other properties of these sets. For instance, we have remarked that the Mather sets, being the support of invariant probability measures, are recurrent under the flow. This is not true anymore for the Aubry and Ma\~n\'e sets, but something can still be said. Let us first recall the definition of 
{\it $\e$-pseudo orbit}. Given a (compact) metric space $X$ and a flow $\varphi$ on it, we say that there exists an $\e$-pseudo orbit between two points $x,y\in X$, if we can find  $\{x_n\}_{n=0}^{k_\e} \subset X$ and positive times $t_1,\ldots, t_{k_\e} > 0$ such that $x_0=x$, $x_{k_\e}=y$ and ${dist}\!\left(\varphi_{t_{i+1}}\left(x_i\right), x_{i+1} \right) \leq \e$ for all $i=0,\ldots,k_{\e}$.

\begin{Prop}\label{chainproperties}
{\rm (i)} $\Phi^L\big|\widetilde{\cN}_c$ is chain transitive, \ie  for each $\e>0$ and for all $(x,v),(y,w)\in \widetilde{\cN}_c$, there exists an $\e$-pseudo-orbit for the flow $\Phi^L$ connecting them. \\
{\rm (ii)} $\Phi^L\big|\widetilde{\cA}_c$ is chain recurrent, \ie  for each $\e>0$ and for all $(x,v)\in \widetilde{\cA}_c$, there exists an $\e$-pseudo-orbit for the flow $\Phi^L$ connecting $(x,y)$ to itself. 
\end{Prop}

\noindent The proof of this result can be found for instance in \cite[Theorem V]{ManeII}\\

As a consequence of the chain-transitivity it follows that the Ma\~n\'e set must be connected (the Aubry set in general not).

\begin{Cor}\label{Maneconnected}
The Ma\~n\'e set is connected
\end{Cor}

%%%%%%%%%%%%%%%%%%%%%%%%%%%%%%%%%%%%%%%%%%%%%%%%%
%%%%%%%%%%%%%%%%%%%%%%%%%%%%%%%%%%%%%%%%%%%%%%%%5
\vspace{20pt}

\noindent{\bf ADDENDA}\\

\noindent{\bf {\sc 5.A - Some topological and symplectic properties of these sets}}\\
\addcontentsline{toc}{subsection}{\hspace{15 pt} 5.A - Some topological and symplectic properties of these sets}

In this Addendum we want to discuss (without any proof) some topological and symplectic properties of the Aubry and Ma\~n\'e sets, similar to what we have alread seen and proved for the Mather sets.

In Section \ref{sec1.3}, Proposition \ref{flatness}, we had related the intersection of Mather sets corresponding to different cohomology classes, to the ``flatness'' of the $\a$ function. The same result holds for the Aubry set and has been proved by Daniel Massart in \cite[Proposition 6]{MassartActionfunctional}. However, the proof in this case is less straightforward and more involved.\\

\begin{Prop}[{\bf Massart}, \cite{MassartActionfunctional}]
Let $c\in \rH^1(M;\R)$ and denote by $F_c$ maximal face of the epigraph of $\a$  containing $c$ in its interior. \\
{\rm (i)} If a cohomology class $c_1$ belongs $F_c$, then $\cA_c\subseteq \cA_{c_1}$. In particular, if $c_1$ belongs to the interior of $F_c$, then they coincide, \ie $\cA_c = \cA_{c_1}$.\\
{\rm (ii)} Conversely, if two cohomology classes $c$ and $c_1$ are such that $\widetilde{\cA}_c \cap \widetilde{\cA}_{c_1} \neq \emptyset$, then for each $\l\in [0,1]$ we have
$\a(c)=\a(\l c + (1-\l)c_1)$, \ie the epigraph of $\a$ has a face containing $c$ and $c_1$.
\end{Prop} 

In particular,  Massart proved that it is possible to relate the dimension of a ``face'' of the epigragph of the $\a$-function to the topological complexity  of the Aubry sets corresponding to  cohomologies in that face (see \cite[Theorem 1]{MassartActionfunctional}).  More precisely, for any sufficiently small $\e>0$, let us define $C_c(\e)$ be the set of integer homology classes which are represented by a piecewise $C^1$ closed curve made with arcs contained in $\cA_c$ except for a remainder of total length less than $\e$. Let $C_c:=\bigcap_{\e>0} C_c(\e)$. Let $V_c$ be the space spanned in $\rH_1(M;\R)$ by $C_c$. Note that $V_c$ is an integer subspace of $H_1(M;\R)$, that is it has a basis of integer elements (images in $\rH_1(M;\R)$ of elements in $\rH_1(M;\Z)$).\\

We denote by:\\
- $F_c$ the maximal face (flat piece) of the epigraph of $\a$, containing $c$ in its interior;\\
- Vect $F_c$ the underlying vector space of the affine subspace generated by $F_c$ in $H^1(M;R)$;\\
- $V_c^{\perp}$ the vector space of cohomology classes of $C^1$ $1$- forms that vanish on $V_c$;\\
- $G_c$ the vector space of cohomology classes of  $C^1$ $1$-forms that vanish in $\rT_xM$ for each $x\in \cA_c$;\\
- $E_c$ the space of cohomology classes of $1$-forms of class $C^1$, the supports of which are disjoint from $A_c$.\\

\begin{Teo}[{\bf Massart, \cite{MassartActionfunctional}}]
$E_c \subseteq {\rm Vect}\,F_c \subseteq G_c \subseteq V_c^{\perp}$.
\end{Teo}

\vspace{10 pt}

Moreover, as we have already proved for the Mather sets (see Addendum 4A, Proposition \ref{4c2}), also these sets are symplectic invariant. The same holds also for these other sets, but the proof in this case is definitely less trivial and requires a more subtle study of these action-minimizing orbits. It can be deduced, for instance, as a special case of \cite[Theorem 1.10]{BernardSympl} (which also applies to the non-autonomous case). 
Let us denote by $\cA^*_c(H)$ and $\cN^*_c(H)$ the Aubry and Ma\~n\'e sets associated to a Tonelli Hamiltonian $H$ (in the sense of the Legendre transform of the corresponding ones for the associated Lagrangian). Then:

\begin{Teo}[{\bf Bernard, \cite{BernardSympl}}]
Let $L:\rT M \longrightarrow \R$ be a Tonelli Lagrangian and $H:\rT^* M \longrightarrow \R$ the associated Hamiltonian. If $\Phi: \rT^*M \longrightarrow \rT^*M$ is an exact symplectomorphism, then
$$
\cA^*_c(H\circ \Phi) = \Phi^{-1}\left(\cA^*_c(H)\right) \qquad \mbox{and}\qquad
\cN^*_c(H\circ \Phi) = \Phi^{-1}\left(\cN^*_c(H)\right).
$$
\end{Teo}

This result can be easily extended to non-exact symplectomorphisms, using 
Lemma \ref{symplectomorphism} (as we have already done in the case of Proposition \ref{4c1}).

\begin{Teo}
Let $L:\rT M \longrightarrow \R$ be a Tonelli Lagrangian and $H:\rT^* M \longrightarrow \R$ the associated Hamiltonian. If $\Psi: \rT^*M \longrightarrow \rT^*M$ is a symplectomorphism of  class $[\Psi]$, then
$$
\cA^*_c(H\circ \Psi) = \Psi^{-1}\left(\cA^*_{c+[\Psi]}(H)\right) \qquad \mbox{and}\qquad
\cN^*_c(H\circ \Psi) = \Psi^{-1}\left(\cN^*_{c+[\Psi]}(H)\right).
$$
\end{Teo}

\begin{Rem}
More geometric proofs of this result can be obtained using weak KAM theory, for instance  \cite[Lemma 1]{IntegTonelli} (it does not apply to all symplectomorphism, but only to those  that are Hamiltonianly isotopic to the identity) and \cite{BernSantos}. See also \cite{PaternainSiburg}.
\end{Rem}

\vspace{20 pt}
\noindent{\bf {\sc 5.B - An example: the simple pendulum II}}\\
\addcontentsline{toc}{subsection}{\hspace{15 pt} 5.B - An example: the simple pendulum II}

In Addendum 4.B we discussed the Mather sets, the $\a$-function and the $\b$-function, in the case of the special case of a simple pendulum, described by the Lagrangian:
\beqano
L: \rT \T &\longrightarrow& \R \\
(x,v) &\longmapsto& \frac{1}{2}|v|^2 + \big(1-\cos (2\pi x)\big).
\eeqano
In particular we proved that:
\begin{itemize}
\item  {For all} $-\frac{4}{\pi} \leq c \leq \frac{4}{\pi}$, $\widetilde{\cM}_{c}= \{(0,0)\}$;
\item  if  $c>\frac{4}{\pi}$, 
 $$ \widetilde{\cM}_{\pm c} = \{(x,v): \; v=  \pm \sqrt{2[(1+\a(c))-\cos(2\pi x)]},\; \forall\, x\in\T\}.$$
\end{itemize}

We want to see which are the Ma\~n\'e and Aubry sets in this case. Let us start recalling the following fact, that we have proved, in a slightly different form, in Section \ref{cartoonexample}.

\begin{Prop}\label{minorbLaggraphs}
Let $\Lambda$ be a $c$-invariant (Lipschitz) Lagrangian graph in $\rT^*M$. Then, the projection on $\rT M$ of each orbit on $\L$ is $c$ semi-static.
\end{Prop}

See for instance \cite{torikam}. The proof is essentially the same as Proposition \ref{proporbittori} (see also Remarks \ref{remmaneset} ({\it ii}) and  \ref{someremarkssemistatic} ({\it iii})). The proof extends to the case of Lipschitz $c$-Lagrangian graphs (\ie Lipschitz sections that are  locally the graph of closed $1$-forms of cohomology class $c$). Observe in fact that  also in the Lipschitz case, the Hamiltonian keeps constant on invariant Lagrangian graphs (see for instance \cite{IntegTonelli}).\\

Hence, it follows from Proposition \ref{minorbLaggraphs} that for $c>\frac{4}{\pi}$, 
 $$ \widetilde{\cN}_{\pm c} \supseteq \{(x,v): \; v=  \pm \sqrt{2[(1+\a(\pm c))-\cos(2\pi x)]},\; \forall\, x\in\T\} = \widetilde{\cM}_{\pm c}.$$
 In fact, it is easy to check that   $\{(x,v): \; v=  \pm \sqrt{2[(1+\a(c))-\cos(2\pi x)]},\; \forall\, x\in\T\}$ is the graph of a closed $1$-form of cohomology class $\pm c$ (the cohomology is just the signed area enclosed between this graph and the $x$-axis). See also the discussion in Addendum 4.B.\\
 Moreover, $\widetilde{\cN}_{\pm c}$ must be equal to $\widetilde{\cM}_{\pm c}$, since it is connected (Corollary \ref{Maneconnected}) and it must be contained in the energy level corresponding to the value $\a(\pm c)$.
 Therefore, recalling the inclusions in Theorem \ref{diagramma}, we can conclude that 
  $$\widetilde{\cN}_{\pm c} = \widetilde{\cA}_{\pm c} = \widetilde{\cM}_{\pm c} \quad \mbox{for all}\; c>\frac{4}{\pi}.$$
 
Let us see what happens for $|c|\leq \frac{4}{\pi}$. Observe that they all correspond to the same energy level, namely the one of the separatrices ($\a(c)=0$ in this case). In this energy level there are exactly three orbits:
\begin{itemize}
\item[-] The fixed point $(0,0)$, \ie $\g_0(t)\equiv 0$.
\item[-] The upper separatrix $\g_+$; for instance let us choose the parametrization given by 
$$\g_+(t)=\pi \left(\Phi^L_t({1}/{2},2)\right).$$
\item[-] The lower separatrix $\g^-$; for instance let us choose the parametrization given by 
$$\g_-(t)=\pi\left(\Phi^L_t({1}/{2},-2)\right).$$
Observe that because of the symmetry of $L$ and the chosen parametrizations, we have that 
$\g_+(t) = \g_{-}(-t)$ for all $t\in \R$.
\end{itemize}

First of all, let us show that for $|c|<\frac{4}{\pi}$, neither $\g_+$ or $\g_-$ can be $c$ semi-static. In fact, let us consider  a $|c|$-closed $1$-form $\eta_c$ whose graph is contained in the region between the separatrices . This is possible since $c$ is less than $4/\pi$ and the cohomology represents the signed area of the region between the curve and the $x$-axis. Moreover, since $|c|$ is strictly less that $4/\pi$, there will be a positive measure subset of $\T$ on which $H(x,\eta_c(x))< 0$.\\
If $\g_+$ were $c$-semistatic (similarly for $\g_-$), then using that $\g_+$ is asymptotic in the past and in the future to $0$, that $\a(c)=0$ and that the Ma\~n\'e potential $\phi_{\eta_c, 0}$ is Lipschitz continuous (Proposition \ref{propertymanepot} (3)), we obtain:
\beqano
\phi_{\eta_c, 0}(0,0) &=& \lim_{T\to +\infty} \phi_{\eta_c, 0}(\g_+(-T),\g_+(T)) = \\
&=& \lim_{T\to +\infty} \int_{-T}^T L_{\eta_c}(\g_+(t),\dot{\g}_+(t))\,dt = \\
&\geq & \lim_{T\to +\infty} \int_{-T}^T \big( \eta_c(\g_+(t))\cdot \dot{\g}_+(t) -  H(\g_+(t), {\eta_c}(\g_+(t)))\big)\,dt =\\
&= &  - \lim_{T\to +\infty} \int_{-T}^T  H(\g_+(t), {\eta_c}(\g_+(t)))\,dt >0 
\eeqano
where the third inequality comes from Fenchel-Legendre inequality and the last one from the fact that there exists a positive measure set of $\T$ in which $H(\g_+(t), {\eta_c}(\g_+(t)))$ is strictly negative. But this leads to a contradiction, since the action of the constant path $\g_0$ is zero: $A_{L_{\eta_c}}(\g_0)=0$.\\
We have just proved that (use also Theorem \ref{diagramma}):
$$\widetilde{\cN}_{c} = \widetilde{\cA}_{c} = \widetilde{\cM}_{c} = \{(0,0)\} \quad \mbox{for all}\; |c|<\frac{4}{\pi}.$$

Finally, let us consider the case $c_{\pm}:=\pm \frac{4}{\pi}$. As above, it follows from Proposition \ref{minorbLaggraphs} (taking the graphs of the separatrices as invariant Lipschitz Lagrangian graphs) that: 
 $$ \widetilde{\cN}_{c_{\pm}} = \{(x,v): \; v=  \pm \sqrt{2[1-\cos(2\pi x)]},\; \forall\, x\in\T\} \supset\widetilde{\cM}_{c_{\pm}}.$$
We want to show that also in this case $\widetilde{\cA}_{c_{\pm}} = \widetilde{\cN}_{c_{\pm}}$ (but observe that this time they contain  the Mather set properly).
Let us denote by $\eta_{\pm}$ the closed $1$-forms given by the graphs of (respectively) the upper and lower separatrix. We have already pointed out that $[\eta_{\pm}]=c_{\pm}$. The key observation is the following lemma.

\begin{Lem} For every $x,y\in \T$ we have
$
\phi_{\eta_{+}, 0}(x,y) = - \phi_{\eta_{-}, 0}(y,x).
$
\end{Lem}

Observe that it follows immediately from this Lemma that the upper (resp. lower) separatrix is not only 
$c_{+}$ semi-static (resp. $c_{-}$ semi-static), but it is indeed $c_{+}$ static (resp. $c_{-}$ static). Therefore, 
$
\widetilde{\cN}_{c_{\pm}} = \widetilde{\cA}_{c_{\pm}}.
$

\begin{Proof}
The above equality is always true if $x=y$ (see Proposition \ref{propertymanepot} (4)). 
Suppose without any loss of generality that $x<y$ and let $S<T$ such that
$\g_{+}(S)=x = \g_{-}(-S)$ and $\g_{+}(T)=y = \g_{-}(-T)$ (this is true only for the parametrization that we chose above). Now, using that these curves are semi-static (for their respective cohomologies), the symmetry of $L$, the relation between $\g_+$ and $\g_-$ and the fact that $\eta_+ = - \eta_-$, we obtain:
\beqano
\phi_{\eta_{+}, 0}(x,y) &=& \int_S^T \big(L(\g_+(t),\dot{\g}_+(t)) - \eta_+(\g_+(t))\cdot \dot{\g}_+(t) \big)\,dt = \\
&=& \int_S^T \big(L(\g_-(-t), -\dot{\g}_-(-t)) + \eta_-(\g_-(t))\cdot (-\dot{\g}_-(-t)) \big)\,dt = \\
&=& - \int_{-T}^{-S} \big(L(\g_-(s), \dot{\g}_-(s)) - \eta_-(\g_-(s))\cdot \dot{\g}_-(s) \big)\,dt = \\
&=& - \phi_{\eta_{-}, 0}(y,x)\,.
\eeqano
\end{Proof}

Summarising:
\begin{itemize}
\item For $c > \frac{4}{\pi}$:
$$\widetilde{\cN}_{\pm c} = \widetilde{\cA}_{\pm c} = \widetilde{\cM}_{\pm c} = 
 \pm \sqrt{2[(1+\a(c))-\cos(2\pi x)]},\; \forall\, x\in\T\}\,.$$
\item For $|c| < \frac{4}{\pi}$:
$$\widetilde{\cN}_{c} = \widetilde{\cA}_{c} = \widetilde{\cM}_{c} = \{(0,0)\}\,.$$
  \item For $c_{\pm}=\pm \frac{4}{\pi}$:
   $$ \widetilde{\cN}_{c_{\pm}} = \widetilde{\cA}_{c_{\pm}}= \{(x,v): \; v=  \pm \sqrt{2[1-\cos(2\pi x)]},\; \forall\, x\in\T\} \supset \widetilde{\cM}_{c_{\pm}} = \{(0,0)\}.$$
\end{itemize}
\vspace{10pt}

\begin{Rem}\label{remarkexamples}
({\it i}) For $c_{\pm}=\pm \frac{4}{\pi}$ we have examples in which the Mather set is strictly contained in the Aubry and Ma\~n\'e sets.\\
\noindent({\it ii}) In all the above examples, the Aubry set coincides with the Ma\~n\'e set. However it is possible to find examples in which this is not anymore true.  For instance, consider the double covering $2\T$ of $\T$ and let $\rho: 2\T \longrightarrow \T$ be the covering map. Let us lift the Lagrangian of the simple pendulum to a Lagrangian $\tilde{L}$ on $\rT(2\T)$, given by $\tilde{L}(\tilde{x},v) = L(\rho(\tilde{x}),v)$. Observe that one can equivalently consider the system on $\rT\T$ given by 
$L(x,v)=\frac{1}{2}|v|^2 - \big(\cos (4\pi x) - 1\big)$. One can verify that (we specify the dependence on the Lagrangian and consider the projected sets):
\begin{itemize}
\item ${\cM}_{0}(\tilde{L}) = \rho^{-1}\left( {\cM}_{0}(L)\right) = \rho^{-1}\left( \{0\}\right) = \{0,1/2\}$.
\item ${\cA}_{0}(\tilde{L}) = \rho^{-1}\left( {\cA}_{0}(L)\right) = \rho^{-1}\left( \{0\}\right) = \{0,1/2\}$.
\end{itemize}
This result is true in general when we consider a finite covering (see \cite[Lemma 2.3]{Contrerasconnecting}).
However, the same is not true anymore for the Ma\~n\'e set. In fact, one can check that the four separatrices connecting the two minimizing fixed points are all $0$-semistatic and therefore:
$${\cN}_{0}(\tilde{L}) = 2\T \supset \{0, \frac{1}{2}\} = \rho^{-1}\left( {\cN}_{0}(L)\right). 
$$
In fact, it is easy to verify that the lifted system has two exact invariant Lagrangian graphs (zero-cohomology = zero area)  made by the lifts of the separatrices (combined so to enclose zero area). Therefore, the lifts of the separatrices are $0$ semi-static (since they are contained on $0$-Lagrangian graphs), but they are not $0$-static  (see figure 3).

\begin{figure}
\begin{center}
\setlength{\unitlength}{1 cm}
\begin{picture}(6,5)(-2,-2) 
\put(-4, -2){\vector(0, 1){5}}
\qbezier(-4,0)(0,0)(4,0)
\put(-4,0){\circle*{0.15}}
\put(0,0){\circle*{0.15}}
\put(4,0){\circle*{0.15}}
\put(4.2,-0.2){$2\T$}
\qbezier(-4,0)(-2,4)
(0,0)
\qbezier(0,0)(2,-4)(4,0)
\put(-1.9, 2){\vector(1, 0){0}}
\put(2, -2){\vector(-1, 0){0}}
\qbezier[70](-4,0)(-2,-4)
(0,0)
\qbezier[70](0,0)(2,4)(4,0)
\put(-1.9, -2){\vector(-1, 0){0}}
\put(2, 2){\vector(1, 0){0}}
\put(-0.7,1.2){$\L_1$}
\put(3.3,1.2){$\L_2$}
\end{picture}
\vspace{10 pt}
\caption{The two invariant Lagrangian graphs $\L_1$ (continuous line) and $\L_2$ (dashed line) for the ``lifted'' pendulum on $2\T$. }
\end{center}
\end{figure}
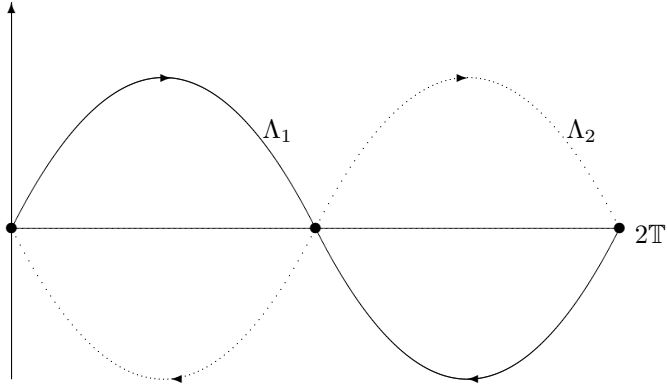

\end{Rem}

%%%%%%%%%%%%%%%%%%%%%%%%%%%%%%%%%%%%%%%%%%%%%%%%
%%%%%%%%%%%%%%%%%%%%%%%%%%%%%%%%%%%%%%%%%%%%%%%%
%%%%%%%%%%%%%%%%%%%%%%%%%%%%%%%%%%%%%%%%%%%%%%%%
\vspace{20 pt}
\noindent{\bf {\sc 5.C - Mather's approach: Peierls' Barrier}}\\
\addcontentsline{toc}{subsection}{\hspace{15 pt} 5.C - Mather's approach: Peierls' Barrier}

In this addendum we want to describe Mather' original approach to the study of action-minimizing curves. As far as the definition of the Ma\~n\'e  set is concerned, it is pretty much the same as the one we have described before. The main difference is in the definition of the Aubry set. In fact, Mather's original definition considered what he called   $c$-{\it regular-minimizers}. We shall see that they indeed coincide with $c$-static curves.

The main ingredient in Mather's approach is the notion of {\it Peierls Barrier}, introduced in 
\cite{Mather93}.\footnote{The function that we are defining here is a actually a slight generalization of $h^{\infty}_c$ defined in \cite{Mather93}. Pay attention that throughout this article, the sign of the $\a$ function is wrong: wherever there is $\a(c)$, it should be substituted by $-\a(c)$.}

For $t>0$ and $x,\,y \in M$, let us consider :
\beqa{ht}
 h_{\eta,t}(x,y) = \min \int_0^t L_{\eta}(\g(s),\dot{\g}(s))\,ds\,,
\eeqa
where the minimum is taken over all piecewise $C^1$ paths $\g: [0,t]\longrightarrow M$, such that $\g(0)=x$ and $\g(t)=y$. This minimum is achieved because of Tonelli theorem (Theorem \ref{TonelliTheorem}).
We define the {\it Peierls barrier} as:
\beqa{defPeierlsbarrier} h_{\eta}(x,y) = \liminf _{t \rightarrow +\infty} (h_{\eta,t}(x,y) + \a(c)t)
\,.\eeqa

\begin{Rem}
({\it i}) Observe that $h_{\eta}$ does not depend only on the cohomology class $c$, but 
also on the choice of the representative $\eta$; namely, if $\eta' = \eta + df$, then 
$h_{\eta'}(x,y) = h_{\eta}(x,y) + f(y) - f(x)$. Anyhow, this dependence will not be harmful for what we are going to do in the following (it will not change the set of action-minimizing curves).\\
({\it ii}) The main difference between Peierls barrier and Ma\~n\'e potential is that in this case we  consider curves defined over longer and longer time intervals. In particular, it is easy to check that:
$$
h_{\eta} (x,y) \geq \Fac(x,y) \quad \forall\; x,y\in\, M.
$$
({\it iii}) This function $h_{\eta}$ is a generalization of Peierls Barrier introduced by Aubry \cite{AubryLeDaeron} and Mather \cite{Mather85,Mather87,Matherorbitsdiffeo}  in their study of twist maps. In some sense we are comparing, in the limit, the action of Tonelli minimizers of time length $T$ with the corresponding average $c$-minimimal action $-\a(c)T$. Remember, in fact, that $-\a(c)$ is the ``average action'' of a $c$-minimal measure.\\
({\it iv}) Albert Fathi \cite{Fathibook} showed that - in the autonomous case - this $\liminf$ can be replaced with a $\lim$. This is not generally true in the non-autonomous time-periodic case (see for instance \cite{MatherFathi} for some counterexamples); Tonelli Lagrangians for which this convergence result holds are called {\it regular}. Patrick Bernard \cite{Bernardconnecting} showed that under suitable assumptions on the Mather set it is possible to prove that the Lagrangian is regular. For instance, if the Mather set $\tilde{\cM}_c$ is union of $1$-periodic orbits, then $L_{\eta}$ is regular.  This problem turned out to be strictly related to the convergence of the so-called {\it Lax-Oleinik semigroup} (see \cite{Fathibook} for its definition).
\end{Rem}

Analogously to Proposition \ref{propertymanepot}, one can prove the following.

\begin{Prop} \label{Peierls}
The values of the map $h_{\eta}$ are finite. Moreover, the following
properties hold:
\begin{itemize}
\item[{i)}]  for each $x,\,y,\,z \in M$ and $t>0$ 
$$h_{\eta}(x,y) \leq h_{\eta}(x,z) + h_{\eta,t}(z,y) + \a(c)t; $$
\item[ii)] for each $x,y,z\in M$,  $h_{\eta}(x,y) \leq h_{\eta} (x,z) + h_{\eta}(z,y)$.
\item[iii)] for each $x,y,z\in M$,  $h_{\eta}(x,y) \leq h_{\eta} (x,z) + \Fac(z,y)$.
\item[{iv)}]  for each $x,\,y,\,z \in M$ and $t>0$ 
$$h_{\eta}(x,y) \leq h_{\eta}(x,z) + h_{\eta,t}(z,y) + \a(c)t; $$
\item[{v)}]  for each $x\in M$, $h_{\eta}(x,x)\geq 0$;
\item[{vi)}]  for each $x,\,y \in M$, $h_{\eta}(x,y) + h_{\eta}(y,x) 
\geq 0$.
%\item[{iv)}]   for each $x\in \cM_c$, $h_{\eta}(x,x) = 0$; 
\end{itemize}
\end{Prop}

%\begin{Proof}
%The fact that the values of this map are finite, will follow from (i).... DA FARE!!!\\
%\end{Proof}

It is interesting to consider the following symmetrization:
\beqa{pseudodistance} \d_c:\; M\times M &\longrightarrow & \R \nonumber\\
(x,y) &\longmapsto & h_{\eta}(x,y) + h_{\eta}(y,x). \eeqa 
Observe that this function does now depend only on the cohomology class $c$ and moreover
it is non-negative, symmetric and satisfies the triangle inequality.\\
An interesting property of  $\d_c$ is the following (see \cite[Section 8]{Mather93}). If $d$ denotes the distance induced on $M$ by the Riemannian metric $g$, then there exists $C>0$ such that for each $x,y\in M$ we have 
$$\d_c(x,y) \leq C d(x,y)^2.$$ 

\begin{Rem}
The same estimate continues to be true for the non-autonomous time-periodic case. In this case we have that
$$\d_c((x,\t_0),(y,\t_1)) \leq C [d(x,y) + \|\t_1-\t_0\|]^2$$
for each $(x,\t_0),(y,\t_1)\in \cA_c$, where 
$$
\|\t_1 -\t_0\| = \inf \left\{ |t_1-t_0|:\; t_i\in\R,\; t_i \equiv \t_i\, ({\rm mod.}\ 1),\; i=0,1 \right\}.
$$
\end{Rem}

Let us see now some relation between this Peierls barrier (or equivalently $\d_c$) and $c$-action minimizing curves.  Let $\g:\R \longrightarrow M$ be a $c$-minimizer and consider $x_{\a},x'_{\a}$ in the $\a$-limit set\footnote{
 Recall that a point $z$ is in the $\a$-limit set of $\g$, if there exists a sequence $t_{n} \rightarrow -\infty$  such that $\g(t_n)\rightarrow z$.}
 of $\g$ and $x_{\om},x'_{\om}$ in the $\om$-limit set\footnote{
 Recall that a point $z$ is in the $\omega$-limit set of $\g$, if there exists a sequence $t_{n} \rightarrow +\infty$  such that $\g(t_n)\rightarrow z$.} of $\g$.
John Mather in \cite[Section 6]{Mather93} proved that $\d_c(x_{\a},x'_{\a}) = \d_c(x_{\om},x'_{\om}) = 0$. 
In general, it is not true that $\d_c(x_{\a},x_{\om})$=0; what one can prove is that this value does not depend on the particular $x_{\a}$ and $x_{\om}$, \ie $\d_c(x_{\a},x_{\om})~=~\d_c(x'_{\a},x'_{\om})$: it is a property of the limit sets rather than of their elements. 
Nevertheless, there will exist particular $c$-minimizers for which this value is equal to $0$ and these will be the $c$-minimizers that we want to single out.

\begin{Def}[{\bf c-regular minimizers}]\label{alphaomegalimit}
A $c$-minimizer $\g:\R \longrightarrow M$ is called a $c$-{\it regular minimizer}, if $\d_c(x_{\a},x_{\om})=0$ for each $x_{\a}$ in the $\a$-limit set of $\g$ and $x_{\om}$ in the $\om$-limit set of $\g$.
\end{Def}

Mather defined   the {\it Aubry set} as the union of the support of all these $c$-regular minimizers.

\begin{Def}[{\bf Aubry set}]
The {\it Aubry set} (with cohomology class $c$) is:
\beqano 
\widetilde{\cA}_{c}= \bigcup \left\{(\g(t),\dot{\g}(t)):\; \text{$\g$ is a
$c$-regular minimizer and}\;t\in\R\right\}.
\eeqano
\end{Def}

it turns out that this set coincides exactly with the one that we have defined in Section \ref{sec1.3}. In fact, one ca prove that:

\begin{Prop}
$\g$ is a $c$-regular minimizer of $L$ if and only if $\g$ is a $c$-static curve of $L$.
\end{Prop}

A proof of this can be found in \cite[Proposition 9.2.5]{Fathibook}. The essential ingredient is that $\Fac(x,y)=h_{\eta}(x,y)$ if $x,y \in \cA_c$ (see also Remark \ref{peierlcoincidemane}).\\

Moreover, one can also provide another alternative definition of the (projected) Aubry set:

\begin{Prop}[{See \cite[Proposition 5.3.8]{Fathibook}}]\label{othercharacAubry}
The following properties are equivalent.
\begin{itemize}
\item[i) ]$x \in \cA_c$; 
\item[ii)]$h_{\eta}(x,x)=0$;
\item[iii)] there exists a sequence of absolutely continuous curves $\g_n:[0,t_n]\rightarrow M$ such that:
\begin{itemize}
\item[-] for each $n$, we have $\g_n(0)=\g_n(t_n)=x$;
\item[-] the sequence $t_n \rightarrow +\infty$, as $n\rightarrow +\infty$;
\item[-] as $n\rightarrow +\infty$, $\int_0^{t_n} L_{\eta}(\g_n(s),\dot{\g}_n(s))\,ds + \a(c)t_n \rightarrow 0$.\\
\end{itemize}
\end{itemize}
\end{Prop}

\begin{Rem}
({\it i}) Therefore, the Aubry set consists of points that are contained in loops with period as long as we want and action as close as we want to the minimal average one.\\
({\it ii}) Moreover, it follows from $ii)$ in Proposition \ref{othercharacAubry} that $\d_c$ is  a
pseudometric on the projected Aubry set 
$$\cA_{c} =\{x\in M : \; \d_c(x,x)=0 \}\,.$$
({\it iii}) One can easily construct a metric space out of ($\cA_c,\d_c$). We call {\it quotient 
Aubry set} %, or {\it Mather quotient},  
the metric space $(\cbA_{c},\, \bd_c)$
obtained by identifying two points in $\cA_{c}$, if their $\d_c$-pseudodistance
is zero. This set plays quite an interesting role in the study of the dynamics; see for example
\cite{Mather03, Mather04, quotientAubry, FathiFigalliRifford} for more details.
\\
\end{Rem}

%%%%%%%%%%%%%%%%%%%%%%%%%%%%%%%%%%%%%%%%%%%%%%%%%%%
%%%%%%%%%%%%%%%%%%%%%%%%%%%%%%%%%%%%%%%%%%%%%%%%%%

%%%%%%%%%%%%%%%%%%%%%%Section 1.4  %%%%%%%%%%%%%%%%%%%%%%%%%%%%%%%%%%%%%%%%%%%%%

\section{Weak KAM theory}\label{sec1.4}

Another interesting approach to the study of these invariant sets
is provided by the so-called {\it weak KAM theory}, which represents the functional analytical counterpart of the variational methods discussed in the previous sections.
In section \ref{cartoonexample} in fact, we pointed out the relation between KAM tori (or more generally, invariant Lagrangian graphs) and classical solutions and subsolutions of Hamilton-Jacobi equation
(see Remark \ref{remarktori} ({\it iv})). 
This approach that we are going to describe, will be based on studying   ``weak'' (non-classical)
solutions of Hamilton-Jacobi equation and some special class of subsolutions ({\it critical subsolutions}). From a more geometrical point of view, this can be interpreted as the study of particular Lagrangian graphs (not necessarily invariant) and their non-removable intersection (see also \cite{PaternainSiburg}). This point of view makes this approach  particularly interesting, since it relates the dynamics of the system to the geometry of the space and might potentially open the way to a ``symplectic''  definition of Aubry-Mather theory (see also Sections 4A, 5A and \cite{BernardSympl,  IntegTonelli, BernSantos}).\\
In this section we want to provide a brief presentation of this theory, omitting most  of the proofs, for which we refer the reader to  the excellent - and self-contained - presentation \cite{Fathibook}.\\

The main object of investigation is represented by Hamilton-Jacobi (H-J) equation: 
$$ H_{\eta}(x, d_xu) = H(x, \eta(x) + d_xu) = k\,,  $$
where $\eta$ is a closed $1$-form on $M$ with a certain cohomology class $c$. Observe that considering H-J equations for different
$1$-forms corresponding to different cohomology classes, is equivalent to Mather's idea of changing Lagrangian (see section \ref{sec1.2}).

From now on, we shall consider $L$ to be a Tonelli Lagrangian on a compact manifold $M$ and $H$ its associated Hamiltonian. Let us fix $\eta$ to be a closed $1$-form on $M$ with cohomology class $c$, and, as before, denote by $L_{\eta}$ and $H_{\eta}$ the modified Lagrangian and Hamiltonian.
In classical mechanics, one is interested in studying solutions of this equation, \ie $C^1$ functions $u:M\rightarrow \R$ such that
$H_{\eta}(x, d_xu)=k$. It is immediate to check that for any given cohomology class there exists at most one value of $k$ for which these $C^1$ solutions may exist. In fact, it is enough to observe that if $u$ and $v$ are two $C^1$ functions on a compact manifold, there will exist a point $x_0$ at which their differentials coincide (take any critical point of $u-v$).
We shall see (Theorem \ref{weakKAMtheorem}) that this value of $k$ for which solutions may exist, coincides with $\a(c)$ or Ma\~n\'e critical value $c(L_{\eta})$ (defined in sections \ref{sec1.2} and \ref{sec1.3}).

\begin{Rem}
The existence of such solutions has significant implications to the dynamics of the system and it is, consequently, quite rare. In particular, they correspond to Lagrangian graphs, which are invariant under the Hamiltonian flow $\Phi^H_t$ (Hamilton-Jacobi theorem). For instance, in the case of $M=\T^d$ and nearly-integrable systems these solutions correspond to KAM tori (this might give an idea of their rareness).
\end{Rem}

One of the main results of weak KAM theory is that,  in the case of Tonelli Hamiltonians, a weaker kind of solutions do always  exist. In the following we are going to define these generalized solutions and their relation with the dynamics of the system. It is important to point out that one of the main ingredient in the proof of all these results is provided by Fenchel inequality (cf. \ref{Fenchelineq} in Section \ref{sec1.1}).\\

Let us start by generalizing the concept of subsolution. In the $C^1$-case it is easy to check - using Fenchel inequality - that the following property holds (the proof is essentially the same as for (\ref{ineqforcurves}) in Section \ref{cartoonexample}).
%In order to define these generalized solutions, let us recall some properties of classical solutions and subsolutions.

\begin{Prop}
Let $u: M \rightarrow \R$ be $C^1$; $u$ satisfies $H(x,\eta(x) + d_xu) \leq k$ for all $x\in M$ if and only if for all $a<b$ and $\g:[a,b]\rightarrow M$
$$
u(\g(b))-u(\g(a)) \leq \int_a^b L_{\eta}(\g(t),\dot{\g}(t))\,dt + k(b-a).
$$
\end{Prop}

This last inequality provides the ground for to a definition of  subsolution in the $C^0$-case.

\begin{Def}[{\bf Dominated functions}]\label{defdomination}
Let $u:M\rightarrow \R$ be a continuous function; $u$ is dominated by $L_{\eta}+k$, and we shall write $u\prec L_{\eta}+k$, if  
for all $a<b$ and  $\g:[a,b]\rightarrow M$
\beqa{domination}
u(\g(b))-u(\g(a)) \leq \int_a^b L_{\eta}(\g(t),\dot{\g}(t))\,dt + k(b-a).
\eeqa
\end{Def}

One can check that if $u \prec L_{\eta}+k$ then $u$ is Lipschitz and its Lipschitz constant can be bounded by a constant $C(k)$ independent of $u$; in fact, it is sufficient to apply the definition of dominated function with  the speed-one geodesic connecting any two points $x$ and $y$ and consider the maximum of $L$ over the unit tangent ball (see \cite[Proposition 4.2.1 (iii)]{Fathibook}).  In particular, all dominated functions for values of $k$ in a compact set are equiLipschitz. On the other hand, it is easy to check that each Lipschitz function is dominated by $L_{\eta}+k$, for a suitable $k$ depending on its Lipschitz constant: this shows that dominated functions exist.\\
Dominated functions generalize subsolutions of H-J to the continuous case. In fact:

\begin{Prop}[see {\cite[Theorem 4.25]{Fathibook}}]
\label{propsubsolutions}
If $u\prec L_{\eta}+k$ and $d_{x}u$ exists, then $H(x, \eta(x)+ d_xu) \leq k$. Moreover, if  $u:M\rightarrow \R$ is Lipschitz and 
$H(x, \eta(x)+ d_xu) \leq k$ a.e., then $u\prec L_{\eta}+k$.
\end{Prop}

\begin{Rem}
Using the fact that any Lipschitz function is differentiable almost everywhere (Rademacher theorem), one could equivalently define subsolutions in the following way:
{\it a locally Lipschitz function $u: M \longrightarrow \R$ is a {\it
subsolution} of $H_{\eta}(x,d_xu)=k$, with $k\in \R$, if $H_{\eta}(x,d_xu)\leq 
k$ for almost every $x\in M$.}
\end{Rem}

\begin{Rem}\label{manecriticalsec1.4}
One interesting question is: for which values of $k$ do there exist functions dominated by $L_{\eta}+k$ (or equivalently subsolutions of $H(x,\eta(x)+d_xu)=k$)? It is possible to show that there exists a value $k_c \in \R$ such that 
$H(x,\eta + d_xu)=k$ does not admit any subsolution for $k<k_c$, while it has subsolutions for $k\geq k_c$, see \cite{LPV, Fathibook}. In particular, if $k>k_c$ there exist $C^{\infty}$ subsolutions. It turns out that the constant $k_c$ coincides with $\a(c)$ and the {\it Ma\~n\'e's critical value} (where $c=[\eta]$). See \cite{ContrerasLaggraph,Fathibook}. 
\end{Rem}

Functions corresponding to this ``critical domination'' play an important role, since they encode significant information about the dynamics of the system.

\begin{Def}[{\bf Critical subsolutions}] 
A function $u \prec L_{\eta}+\a(c)$ is said to be {\it critically dominated}. Equivalently, we shall also call it an $\eta$-{\it critical subsolution}, since
$H (x, \eta(x)+ d_xu)\leq \a(c)$ for almost every $x\in M$.\label{etacriticalsub}
\end{Def}

\begin{Rem}\label{charactalpha}
The above observation provides a further definition of $\a(c)$:
$$
\a(c) = \inf_{u\in C^{\infty}(M)} \max_{x\in M} H(x, \eta(x) + d_xu)\,.
$$
This theorem has been proven in \cite[Theorem A]{ContrerasLaggraph} (see also in \cite{Fathibook}).
This infimum is not a minimum, but it becomes a minimum over the set of Lipschitz functions on $M$ (also over the smaller set of $C^{1,1}$ functions, see the addendum at the end of this section and \cite{FathiSiconolfi, BernardC11}).
This characterization has the following geometric interpretation. If we consider the space ${\rm T}^*M$ equipped 
with the canonical symplectic form, 
the graph of the differential of a $C^1$ $\eta$-critical subsolution 
(plus the $1$-form $\eta$) is nothing else than a $c$-Lagrangian graph (\ie a Lagrangian graph with 
cohomology class $c$). Therefore Ma\~n\'e $c$-critical energy level $\cE^*_{c}=\{(x,p)\in {\rT^*M:\; H(x,p)}=\a(c)\}$ corresponds to a 
$(2d-1)$-dimensional hypersurface, such that the region it bounds is convex in each fiber
and does not contain in its interior any $c$-Lagrangian graph, while any of its neighborhoods does. 
\end{Rem}

Analogously to what we have already seen for subsolutions, it would be interesting to investigate if there existed an equivalent characterization of classical solutions of Hamilton-Jacobi, that does not involve the regularity of the solution. This would allow us to define ``{\it weak}'' solutions, that hopefully are not so rare as the classical ones. \\

Let us now recall some properties of classical solutions, which will allow us to provide a ``weaker'' definition of solution (the proof of this proposition is essentially the same as for 
(\ref{ineqforcurves}) in Section \ref{cartoonexample}).

\begin{Prop}[{see \cite[Theorem 4.1.10]{Fathibook}}]
Let $u: M\rightarrow \R$ a $C^1$ function and $k\in \R$. The following conditions are equivalent:
\begin{enumerate}
\item $u$ is solution of $H(x,\eta(x)+d_xu)=k$;
\item $u\prec L_{\eta}+k$ and for each $x\in M$ there exists $\g_x:(-\infty,+ \infty)\rightarrow M$ such that 
$\g_x(0)=x$ and for any $[a,b]$:
$$u(\g_x(b))-u(\g_x(a)) = \int_a^b L_{\eta}(\g_x(t),\dot{\g}_x(t))\,dt + k(b-a).$$
\item $u\prec L_{\eta}+k$ and for each $x\in M$ there exists $\g_x:(-\infty,0]\rightarrow M$ such that $\g_x(0)=x$ and for any $a<b\leq 0$:
$$u(\g_x(b))-u(\g_x(a)) = \int_a^b L_{\eta}(\g_x(t),\dot{\g}_x(t))\,dt + k(b-a).$$
\item $u\prec L_{\eta}+k$ and for each $x\in M$ there exists $\g_x:[0,+\infty)\rightarrow M$ such that $\g_x(0)=x$ and for any $0\geq a<b$:
$$u(\g_x(b))-u(\g_x(a)) = \int_a^b L_{\eta}(\g_x(t),\dot{\g}_x(t))\,dt + k(b-a).$$
\end{enumerate}
\end{Prop}

Inspired by this fact, let us consider the curves for which equality in (\ref{domination}) holds.

\begin{Def}[{\bf Calibrated curves}]
Let $u\prec L_{\eta}+k$. A curve $\g: I \rightarrow M$ is $(u,L_{\eta},k)$-calibrated  if for any $[a,b]\subseteq I$
$$u(\g(b))-u(\g(a)) = \int_a^b L_{\eta}(\g(t),\dot{\g}(t))\,dt + k(b-a).$$
\end{Def}

\noindent These curves are very special curves and it turns out that they are orbits of the Euler-Lagrange flow. In fact:

\begin{Prop}\label{calibminim}
If $u\prec L_{\eta}+k$ and $\g: [a,b]\rightarrow M$ is $(u, L_{\eta}, k)$-calibrated, then $\g$ is a $c$-Tonelli minimizer, \ie
$$
\int_a^b L_{\eta}(\g(t),\dot{\g}(t))\,dt \leq  \int_{a}^{b} L_{\eta}(\s(t),\dot{\s}(t))\,dt
$$
for any $\s: [a,b]\rightarrow M$ such that $\s(a)=\g(a)$ and $\s(b)=\g(b)$. Most of all, this implies that $\g$ is a solution of the Euler-Lagrange flow and therefore it is $C^{r}$ (if $L$ is $C^r$). 
\end{Prop}

\noindent The proof of this result is the same as the one of Proposition \ref{proporbittori} (see also \cite[Proposition 4.3.2 and Corollary 4.3.3]{Fathibook}).

Moreover, the following differentiability result holds.

\begin{Prop}[{see \cite[Theorem 4.3.8]{Fathibook}}]
Let $u\prec L_{\eta}+k$ and $\g:[a,b]\rightarrow M$ be $(u,L_{\eta},k)$-calibrated.
\begin{itemize}
\item[i)] If $d_{\g(t)}u$ exists for some $t\in[a,b]$, then $H(\g(t), \eta(\g(t)) + d_{\g(t)}u)=k$ and $d_{\g(t)}u = \dfrac{\dpr L}{\dpr v}(\g(t),\dot{\g}(t))$.
\item[ii)] If $t\in (a,b)$, then $d_{\g(t)}u$ exists.
\end{itemize}
\end{Prop}

\begin{Rem}
Calibrated curves are ``Lagrangian gradient lines'' of ${\rm grad}_Lu$ (where ${\rm grad}_Lu$ is a multivalued vector field given by the equation
$d_x u = \dfrac{\dpr L}{\dpr v}\left( x, {\rm grad}_Lu\right)$. Therefore, there is only one possibility for calibrated curves, at each point of differentiability of $u$.
\end{Rem}

This suggests the following definitions.

\begin{Def}[{\bf weak KAM solutions}]
Let $u\prec L_{\eta}+ k$.
\begin{itemize}
\item $u$ is a weak KAM solution of negative type (or backward Weak KAM solution) if for each $x\in M$ there exists $\g_x:(-\infty,0]\rightarrow M$ such that $\g_x(0)=x$ and $\g_x$ is $(u,L_{\eta},k)$-calibrated;
\item $u$ is a weak KAM solution of positive type (or forward Weak KAM solution) if for each $x\in M$ there exists $\g_x:[0,+\infty)\rightarrow M$ such that $\g_x(0)=x$ and $\g_x$ is $(u,L_{\eta},k)$-calibrated.
\end{itemize} 
\end{Def}

\begin{Rem}\label{symmetricalLag}
Observe that any weak KAM solution of negative type $u_-$ (resp. of positive type $u_+$) for a given Lagrangian $L$, can be seen as a weak KAM solution of positive type (resp. of negative type) for the {\it symmetrical} Lagrangian $\tilde{L}(x,v):=L(x,-v)$.
\end{Rem}

Let us denote with $\cS^-_{\eta}$ the set of Weak KAM solutions of negative type and $\cS^+_{\eta}$ the ones of positive types. Albert Fathi \cite{FathitheoremeKAM, Fathibook} proved that these sets are always non-empty.

\begin{Teo}[{\bf Weak KAM theorem}]\label{weakKAMtheorem}
There is only one value of $k$ for which weak KAM solutions of positive or negative type of $H(x,\eta(x) + d_xu)= k$ exist. This value coincides with $\a(c)$, where $\a:\rH^1(M;\R)\rightarrow \R$ is Mather's $\a$-function. In particular, for any $u\prec L_{\eta}+\a(c)$ there exist a weak KAM solution of negative type $u_-$ and a weak KAM solution of positive type $u_+$, such that $u_-=u=u_+$ on the projected Aubry set $\cA_c$.
\end{Teo}

Therefore, for any given weak KAM solution of negative type $u_-$ (resp. of positive type $u_+$), there exists a weak KAM solution of positive type of positive type $u_+$ (resp. of negative type $u_-$) such that $u_-=u_+$ on the projected Aubry set $\cA_c$. In particular:

\begin{Prop}[{see \cite[Theorem 4.12.6]{Fathibook}}]
The projected Mather set $\cM_c$ is the uniqueness set for weak KAM solutions of the same type. Namely, if $u_-,v_-$ are weak KAM solutions of negative type (resp. $u_+,v_+$ are weak KAM solutions of positive type) and $u_-=v_-$ on $\cM_c$ (resp. $u_+=v_+$ on $\cM_c$), then they coincide everywhere on $M$.
\end{Prop}

Two solutions $u_-$ and $u+$ that coincide on the (projected) Mather set are said to be {\it conjugate}. We shall denote by $(u_-,u_+)$ a couple of conjugate subsolutions. \\

Let us try to understand the dynamical meaning of such solutions. Albert Fathi \cite{Fathibook}  - using these generalized solutions - proved a {\it weak} version of Hamilton-Jacobi theorem, showing the relation between these weak solutions and the dynamics of the associated Hamiltonian system. We shall state it for weak KAM solution of negative type, but - using remark \ref{symmetricalLag} - one can deduce an analogous statement for weak KAM solutions of positive type.

\begin{Teo}[{{\bf Weak Hamilton-Jacobi Theorem}, \cite[Theorem 4.13.2]{Fathibook}}]\label{HJtheo}
Let $u_-:M \rightarrow \R$ be a weak KAM solution of negative type and consider
$$
{\rm Graph}(\eta + du_-):= \{(x,\eta(x)+d_xu_-), \; {\rm where}\; d_xu_-\; {\rm exists}\}.
$$
Then:
\begin{itemize}
\item[i)] $\overline{{\rm Graph}(\eta + du_-)}$ is compact and is contained in the energy level $\cE^*_{c}=\{(x,p)\in {\rT^*M:\; H(x,p)}=\a(c)\}$;
\item[ii)] $\Phi^H_{-t}\left(\overline{{\rm Graph}(\eta + du_-)} \right) \subseteq {\rm Graph}(\eta + du_-)$ for each $t>0$;
\item[iii)] $M = \pi \left( \overline{{\rm Graph}(\eta + du_-)}\right)$, where $\pi: \rT^*M\rightarrow M$ is the canonical projection.
\end{itemize}
Moreover, let us define:
$$
\cI^*(u_-) := \bigcap_{t\geq 0} \Phi^H_{-t}\left(\overline{{\rm Graph}(\eta + du_-)}\right).
$$
$\cI^*(u_-)$ is non-empty, compact and invariant under $\Phi^H_t$. Furthermore, its ``unstable set'' contains $\overline{{\rm Graph}(\eta + du_-)}$; \ie
$$
\overline{{\rm Graph}(\eta + du_-)} \subseteq W^u\left( \cI^*(u_-)\right) := \left\{ (y,p):\; {\rm dist}\!\left( \Phi^H_{-t}(y,p), \cI^*(u_-)\right) \stackrel{t\rightarrow +\infty}{\longrightarrow} 0   \right\}.
$$
\end{Teo}
\vspace{15 pt}

There is a relation between these invariant sets $\cI^*(u_-)$ (or $\cI^*(u_+)$) and the Aubry set $\widetilde{\cA}_c$ (recall that $\cI^*(u_-), \cI^*(u_+) \subset \rT^* M$, while $\widetilde{\cA}_c \subset \rT M$).\\

\begin{Teo}
\beqa{aubrysetsolution}
\cA^*_c := \cL \left( \widetilde \cA_c \right) &=& \bigcap_{u_{-} \in \cS_{\eta}^-} \cI^*(u_-) = 
\bigcap_{u_{-} \in \cS_{\eta}^-}{\rm Graph}(\eta + du_-) = \nonumber\\
&=& \bigcap_{u_{+} \in \cS_{\eta}^+} \cI^*(u_+) = \bigcap_{u_{+} \in \cS_{\eta}^+}{\rm Graph}(\eta + du_+) =\nonumber\\
&=& \bigcap_{(u_-,u_{+})} \{x\in M:\; u_-(x)=u_+(x)\}, 
\eeqa
where $\cL : \rT M \rightarrow \rT^* M$ denotes the Legendre transform of $L$ and $(u_-,u_+)$ are  conjugates solutions.\\
\end{Teo}

Since it is easier to work with subsolutions rather than weak solutions, we want to discuss now how $\eta$-critical subsolutions, although they contain less dynamical information than weak KAM solutions, can be used to characterize Aubry and Ma\~n\'e sets in a similar way.

Consider $u\prec L_{\eta}+\a(c)$. For $t\geq 0$ define
$$
\widetilde{\cI}_t(u) := \left\{(x,v)\in \rT M:\; \g_{(x,v)}(s):=\pi \Phi^L_s((x,v))\;{\rm is}\; (u,L_{\eta},\a(c))-\mbox{calibr. on}\;(-\infty,t] \right\}.
$$
We shall call the {\it Aubry set of} $u$: $\widetilde{\cI}(u) := \bigcap_{t\geq 0} \widetilde{\cI}_t(u)$, that can be also  defined as
$$
\widetilde{\cI}(u) := 
 \left\{(x,v)\in \rT M:\; \g_{(x,v)}(s):=\pi \Phi^L_s((x,v))\;{\rm is}\; (u,L_{\eta},\a(c))~-~\mbox{calibr. on}\,\,\R \right\}.
$$

These sets $\widetilde{\cI}(u)$ are non-empty, compact and invariant. Moreover, here are some properties of these sets (compare with theorem \ref{HJtheo}).

\begin{Prop}\label{propertyaubry} Let $u\prec L_{\eta}+\a(c)$.
\begin{enumerate}
\item $\widetilde{\cI}_t(u)$ is compact;
\item $\widetilde{\cI}_{t'}(u)\subseteq \widetilde{\cI}_t(u) \subseteq \widetilde{\cI}_0(u)$ for all $t'\geq t\geq 0$;
\item $\cL \left(\widetilde{\cI}_0(u)\right)$ is contained in the energy level $\cE^*_c$ corresponding to $\a(c)$;
\item $\cL \left(\widetilde{\cI}_t(u)\right) \subseteq {\rm Graph}(\eta + du)$ for all $t>0$;
\item $\cL \left(\widetilde{\cI}_0(u)\right) \subseteq \overline{{\rm Graph}(\eta + du)}$  (observe that for weak solutions these two sets coincide);
\item $\Phi^L_{-t} \left(\widetilde{\cI}_0(u)\right) = \widetilde{\cI}_t(u)$ for all $t>0$;
\item $\overline{\bigcup_{t>0} \widetilde{\cI}_t(u)} = \widetilde{\cI}_0(u)$;
\item $ \widetilde{\cI}(u) = \bigcap_{t\geq 0} \Phi^L_{-t}\left(\widetilde{\cI}_0(u)\right)$;
\item $ \widetilde{\cI}_0(u) \subseteq W^u\left( \widetilde{\cI}(u)\right)$.
\item $\pi : \widetilde{\cI}(u) \longrightarrow \pi(\widetilde{\cI(u)})$ is a bi-Lipschitz homeomorphism {\rm[}Graph Theorem{\rm]}. The same is true for $\widetilde{\cI}_t(u)$ for each $t>0$.
\end{enumerate}
\end{Prop}

\begin{Teo}[{\bf Fathi}]\label{theoaubryset}
The Aubry and Ma\~n\'e sets defined in (\ref{Aubryset}) and (\ref{maneset}) can be equivalently defined in the following ways:
\beqano
\widetilde{\cA}_c &=& \bigcap_{u\prec L_{\eta}+\a(c)} \widetilde{\cI}(u) \quad =\quad \bigcap_{u\prec L_{\eta}+\a(c)} \cL^{-1} \left( {\rm Graph}(\eta + du) \right) \\
%                     &=& \bigcap_{u_{-} \in \cS^-} \cL^{-1}\left(\cI^*(u_-)\right) =  \bigcap_{u_{-} \in \cS^-}\cL^{-1} \left({\rm Graph}(\eta + du_-)\right) = \\
%		   &=& \bigcap_{u_{+} \in \cS^+} \cL^{-1}\left(\cI^*(u_+)\right) = \bigcap_{u_{+} \in \cS^+}\cL^{-1}\left({\rm Graph}(\eta + du_+)\right).\\
		   \\
\widetilde{\cN}_c &=& \bigcup_{u\prec L_{\eta}+\a(c)} \widetilde{\cI}(u).
\eeqano
Moreover, there exists $u_{\infty} \prec L_{\eta}+\a(c)$ such that $\widetilde{\cA}_c = \widetilde{\cI}(u_{\infty}) \cap \cL^{-1}\left({\cE}^*_c\right)$. 
\end{Teo}

The proof of this theorem follows from the results in \cite[Chapter 9]{Fathibook}.
For the last statement is sufficient to observe that the set of critically dominated functions is a separable subset of $C(M)$. Let $\{u_n\}$ be a countable dense family of such functions and define $u_{\infty}$ as a convex combination of their normalization (with respect to a fixed point $x_0\in M$), \eg $u_{\infty}(x) = \sum_{n=0}^{\infty} \frac{1}{2^n}\left( u_n(x) - u_n(x_0)\right)$. 

\begin{Rem}
Using this characterization, the graph property of the Aubry set (Theorem \ref{graphtheoremAubry}) follows easily from property (10) in Proposition \ref{propertyaubry}. Moreover, the non-emptiness of $\widetilde{\cN}_c$ is a result of the non-emptiness of $\widetilde{\cI}(u)$. As far as the non-emptiness of $\widetilde{\cA}_c$ is concerned, one can deduce it from this characterization and proposition \ref{propnonwandering}.
\end{Rem}

From theorem \ref{theoaubryset} one can also deduce another interesting property of critically dominated functions:  their differentiability on the projected Aubry set (recall that {a-priori} these functions are only Lipschitz, so they are differentiable almost everywhere).

\begin{Prop}[{see \cite[Theorem 4.3.8 (i)]{Fathibook}}]\label{propdifferential}
Let $u\prec L_{\eta}+\a(c)$. For each $x\in \cA_c$, $u$ is differentiable at $x$ and $d_x u$ does not depend on $u$; namely,
$d_x u = \dfrac{\dpr L}{\dpr v}(x, \pi_{|_{\widetilde{\cA}_c}}^{-1} (x))$.
\end{Prop}

In addition to the Aubry set and Ma\~n\'e set, one can also recover the definition of Ma\~n\'e potential $\Fac$ (see Section \ref{sec1.3}) and Peierls barrier $h_{\eta}$ (see (\ref{defPeierlsbarrier})) in terms of these solutions and subsolutions.

Let us start by observing that, from definition \ref{defdomination}, if $u\prec L_{\eta}+\a(c)$ then for each $x,y\in M$ and $t>0$ we have that
$u(y)-u(x) \leq h_{\eta,t}(x,y) + \a(c)t$, where $h_{\eta,t}$ is defined as
\beqano
 h_{\eta,t}(x,y) = \min \int_0^t L_{\eta}(\g(s),\dot{\g}(s))\,ds\,,
\eeqano
where the minimum is taken over all piecewise $C^1$ paths $\g: [0,t]\longrightarrow M$, such that $\g(0)=x$ and $\g(t)=y$. This minimum is achieved because of Tonelli theorem (Theorem \ref{TonelliTheorem}).

This implies that for each $u\prec L_{\eta}+\a(c)$ and for each $x,y \in M$,  $u(y)-u(x) \leq h_{\eta}(x,y)$, \ie
\beqano
\Fac(x,y) \geq \sup_{u\prec L_{\eta}+\a(c)} \left[u(y) - u(x)\right] \qquad \forall\, x,y\in M\,.
\eeqano 

One can actually show that they are equal.

\begin{Prop}[{see \cite[Corollary 9.1.3]{Fathibook}}]
For each $x,y\in M$, we have the equality
$$
\Fac(x,y) = \sup_{u\prec L_{\eta}+\a(c)} \left[u(y) - u(x)\right].
$$
\end{Prop}

\begin{Rem}
The quantity on the right-hand side is also called ``{\it viscosity semi-distance}'' (see \cite[Section 8.4]{Fathibook}).
\end{Rem}

As far as Peierls barrier is concerned, let us observe that similarly to what happens for Ma\~n\'e potential, also in this case we have that
for each $u\prec L_{\eta}+\a(c)$ and for each $x,y \in M$,  $u(y)-u(x) \leq h_{\eta}(x,y)$, \ie
\beqano
h_{\eta}(x,y) \geq \sup_{u\prec L_{\eta}+\a(c)} \left[u(y) - u(x)\right] \qquad \forall\, x,y\in M\,.
\eeqano 

Moreover, if $u_- \in \cS_{\eta}^-$ and $u_+\in\cS_{\eta}^+$ are conjugate solutions, the same result holds: $u_-(y)-u_+(x) \leq h_{\eta}(x,y)$ and consequently
\beqano
h_{\eta}(x,y) \geq \sup_{(u_-,u_+)} \left[u_-(y) - u_+(x)\right] \qquad \forall\, x,y\in M\,,
\eeqano 
where $(u_-,u_+)$ denotes conjugate weak KAM solutions. In addition to this, it is possible to show that the above inequality is actually an equality. In fact:

\begin{Prop}[{see \cite[Th\'eor\`eme 7]{FathisolutionsKAM}}]
For $x\in M$ let us define the function $h_{\eta}^x: M \rightarrow \R$ (resp. $h_{\eta,x}: M \rightarrow \R$) by $h_{\eta}^x(y)=h_{\eta}(x,y)$
(resp. $h_{\eta,x}(y)= h_{\eta}(y,x))$. For each $x\in M$, the function $h^x_{\eta}$ (resp. $-h_{\eta,x}$) is a weak KAM solution of negative (resp. positive) type. Moreover, its conjugate function $u^x_+ \in \cS_{\eta}^+$  (resp. $u^x_- \in \cS_{\eta}^-$) vanishes at $x$.
\end{Prop}

Therefore:

\begin{Cor}
For each $x,y\in M$, we have the equality
$$
h_{\eta}(x,y) = \sup_{(u_-,u_+)} \left[u_-(y) - u_+(x)\right],
$$
where the supremum is taken over pairs of conjugate solutions. Moreover, for any given $x,y\in M$ this supremum is actually attained.
\end{Cor}

Observe that, since for any $u\prec L_{\eta}+\a(c)$ there exists a weak KAM solution of negative type $u_-$ and a weak KAM solution of positive type $u_+$, such that $u_-~=~u~=~u_+$ on the projected Aubry set $\cA_c$ (see Theorem \ref{weakKAMtheorem}), one can get the following representations for Peierls barrier $h_{\eta}$ on the projected Aubry set:
\beqa{hsubsol}
h_{\eta}(x,y) &=& \sup_{u\prec L_{\eta}+\a(c)} \left[u(y) - u(x)\right]
\eeqa
for all  $x,y\in \cA_c$.  This supremum is actually attained for any fixed $x\in \cA_c$.

\begin{Rem}\label{peierlcoincidemane}
In particular, (\ref{hsubsol}) shows that Peierls barrier and Ma\~n\'e potential coincide on the projected Aubry set.
\end{Rem}

%Moreover, a similar representation also holds for Mather's pseudodistance $\d_c$  (see (\ref{pseudodistance})).
%In fact, from the definition of $\d_c(x,y)$ we immediately get: 
%\beqa{deltasubsol}
%\d_c(x,y) &=& h_{\eta}(x,y) + h_{\eta}(y,x) = \nonumber\\
%&=& \sup_{u\prec L_{\eta}+\a(c)} (u(y)-u(x)) +
%\sup_{v \prec L_{\eta}+\a(c)} (v(x)-v(y)) = \nonumber\\
%&=& \sup_{u, v \prec L_{\eta}+\a(c)} [(u(y)-v(y)) - (u(x)-v(x))] 
%\eeqa 
%for all $x,y\in \cA_c$. This supremum is also attained for any fixed $x,y\in\cA_c$.\\

\vspace{20pt}

%%%%%%%%%%%%%%%%%%%%%%%%%%%%%%%%%%%%%%%%%%%%%%%%%%%%%%%%%%%%%%%%%%%%%%%
{\bf ADDENDUM}\\
\vspace{10 pt}

\noindent{\bf {\sc 6.A - Regularity of critical subsolutions}}\\
\addcontentsline{toc}{subsection}{\hspace{15 pt} 6.A - Regularity of critical subsolutions} 

In this addendum we want to say more about the regularity of critically dominated functions or $\eta$-critical subsolutions.

We have remarked above in this section, that for $k<\a(c)$ there do not exist functions dominated by $L_{\eta}+k$, while for $k\geq \a(c)$ they do exist. Moreover, if $k>\a(c)$ these functions can be chosen to be $C^{\infty}$ (see also characterization of $\a(c)$ in remark \ref{charactalpha}). The critical case has totally different features. As a counterpart of their relation with the dynamics of the system, critical dominated functions have very rigid structural properties, that become an obstacle when someone tries to make them smoother. For instance, as we have recalled in proposition \ref{propdifferential}, if $u\prec L_{\eta} + \a(c)$ then its differential $d_x u$ exists on $\cA_c$ and it is prescribed over there. This means that although it is quite easy to make these functions smoother (\eg $C^{\infty}$) out of the projected Aubry set, it is impossible to modify them on this set.

Nevertheless, Albert Fathi and Antonio Siconolfi \cite{FathiSiconolfi} managed to prove that $C^1$ $\eta$-critical subsolutions do exist and are dense, in the following sense:

\begin{Teo}[{\bf Fathi, Siconolfi}] \label{FathiSiconolfi}
Let $u\prec L_{\eta}+\a(c)$. For each $\e>0$, there exists a $C^1$ function $\tu: M \longrightarrow \R$
such that:
\begin{itemize}
\item[{i)}] $\tu \prec L_{\eta}+\a(c)$;
\item[{ii)}]   $\tu(x)=u(x)$ on $\cA_{c}$;
\item[{iii)}]  $|\tu(x)-u(x)|<\e$ on $M\setminus \cA_{c}$.
\end{itemize}
Moreover, one can choose $\tu$ so that it is a {\it strict $\eta$-critical subsolution}, \ie we have $H_{\eta}(x,d_x\tu)<~a(c)$ on $M\setminus \cA_{c}$.
\end{Teo}

This result has been extended by Patrick Bernard \cite{BernardC11}, who
showed that every $\eta$-critical subsolution coincides, on the Aubry set, 
with a $C^{1,1}$ $\eta$-critical subsolution.

\begin{Teo}[{\bf Bernard, \cite{BernardC11}}]
Let $H$ be a Tonel li Hamiltonian. If the Hamilton-Jacobi equation  has a subsolution, then it has a $C^{1,1}$ subsolution. Moreover, the set of $C^{1,1}$ subsolutions is dense  for the uniform topology in the set of subsolutions. 
\end{Teo}

\begin{Rem}\label{remregsubsolutions}
 In general $C^{1,1}$ is  the best regularity that one  can expect: it is easy in fact to construct examples in which $C^2$ $\eta$-critical subsolutions do not exist. For example, consider the case in which the Aubry set projects over all the manifold $M$ and it is not a $C^1$ graph (\eg on $M=\T$ take $L(x,v)=\frac{1}{2}\|v\|^2 + \sin^2 (\pi x)$ and $\eta= \frac{2}{\pi} dx$). In this case there is only one critical subsolution (up to constants), that is an actual solution: its differential is Lipschitz but not $C^{1}$.\\
It is therefore clear that the structure of the Aubry set plays a crucial role. Patrick Bernard \cite{BernardSmooth} proved that if the Aubry set is a union of finitely many hyperbolic periodic orbits or hyperbolic fixed points, then smoother subsolutions can be constructed. In particular, if the Hamiltonian is $C^k$, then these subsolutions will be $C^k$ too. The proof of this result is heavily based on the hyperbolic structure of the Aubry set and the result is deduced from the regularity of its local stable and unstable manifolds.\\
See also \cite{Fathismooth} for a survey on this problem and some results concerning Denjoy-type obstructions for the existence of regular critical subsolutions on $\T^2$.
\end{Rem}

Using the density of $C^1$  critically dominated functions, one can emprove some of the results in (\ref{aubrysetsolution}), Theorem \ref{theoaubryset} and (\ref{hsubsol}).
Let us  denote by $\cS^1_{\eta}$ the set of $C^1$ $\eta$-critical subsolutions and 
$\cS^{1,1}_{\eta}$ the set of $C^{1,1}$ $\eta$-critical subsolutions. Then:
\beqa{2.8}   \widetilde{\cA}_c &:=& \bigcap_{u\in \cS^1_{\eta}} \widetilde{\cI}(u) = \bigcap_{u\in \cS^{1,1}_{\eta}} \widetilde{\cI}(u)=\nonumber\\
&=& \bigcap_{u\in \cS^1_{\eta}} \cL^{-1}\left({\rm Graph}(\eta + du)\right) = 
\bigcap_{u\in \cS^{1,1}_{\eta}} \cL^{-1}\left({\rm Graph}(\eta + du)\right).
\eeqa
In particular, there exists a $C^{{1,1}}$ $\eta$-critical subsolution $\tilde{u}$
such that:
\beqa{2.9}
{ \widetilde{\cA}}_c &=& \cL^{-1}\left({\rm Graph}(\eta + d\tilde{u}) \cap \cE^*_c\right)=\nonumber\\
&=& \widetilde{\cI}(\tilde{u}) \cap \cL^{-1}\left(\cE^*_c\right).
\eeqa
Moreover, for any $x,\,y \in \cA_{c}$:
\begin{eqnarray}
h_{\eta}(x,y) &=& \sup_{u \in \cS^{1}_{\eta}} \left\{u(y)-u(x)\right\}  = \sup_{u \in \cS^{1,1}_{\eta}} \left\{u(y)-u(x)\right\}\label{hsubsolC1}
%\d_c(x,y) &=& \sup_{u,v\in \cS^{1}_{\eta}} \{(u-v)(y)-(u-v)(x)\} = \nonumber\\
%&=& \sup_{u,v\in \cS^{1,1}_{\eta}} \{(u-v)(y)-(u-v)(x)\},\label{deltasubsolC1}
\end{eqnarray}
where all the above suprema are maxima for any fixed $x,y \in \cA_c$.\\

%%%%%%%%%%%%%%%%%%%%% ADDENDUM 6B

\vspace{10 pt}

\noindent{\bf {\sc 6.B - Non-wandering set of the Ma\~n\'e set}}\\
\addcontentsline{toc}{subsection}{\hspace{15 pt} 6.B - Non-wandering set of the Ma\~n\'e set}

In this addendum  we want now to use this approach and the above-mentioned results to show that the non-wandering set of the Euler-Lagrange flow restricted to the Ma\~n\'e set is contained in the Aubry set (we mentioned this in Proposition \ref{inclusionmathermane}).  %We have already recalled this result for sketching a proof of the second equivalence in theorem \ref{emmeinenne}.
Let us first recall the definition of 
{\it non-wandering point} for a flow $\Phi_t: X \longrightarrow X$.

\vspace{7 pt}

\begin{Def}
A point $x\in X$ is called {\it non-wandering} if for each neighborhood 
$\mathcal{U}$ and each positive integer $n$, there exists $t>n$ such that 
$f^t(\mathcal{U})\cap \mathcal{U} \neq \emptyset$.
\end{Def}

We shall denote the set of {\it non-wandering} points for $\Phi_t$ by 
$\Omega(\Phi_t)$. %Note that, if $\m$ is an invariant measure, then ${\rm supp}\
%\m\subseteq \Omega(\Phi_t)$. In fact, by the ergodic decomposition theorem,
%every point $x\in{\rm supp}\ \m$ is in the support of an ergodic invariant 
%measure $\m_1$: therefore, $x$ is non-wandering, by the ergodicity of $\m_1$. 
%Given this remark, (\ref{app.1}) is a simple consequence of the following 
%Proposition.

\vspace{7pt}

\begin{Prop}\label{propnonwandering} If $M$ is a compact manifold and $L$ 
a Tonelli Lagrangian on ${\rm T} M$, then  
$\Omega \left( \Phi^L_t\big|\widetilde{\cN}_c\right) \subseteq 
\widetilde{\cA}_c$ 
for each $c\in \rH^1(M;\R)$.\\
\end{Prop}

\vspace{5 pt}

\noindent{\bf{Remarks.}}\\
1) Proposition \ref{propnonwandering} also shows that the Aubry set is 
non-empty. In fact, any 
continuous flow on a compact space possesses non-wandering points.\\
2) Since every point in the support of an invariant 
measure $\m$ is non-wandering, then this also shows that $\widetilde{\cM}_c \subseteq 
\widetilde{\cA}_c$.

%if $x\in{\rm supp}\ \m\subseteq\widetilde{\cN}_c$,
%then $x\in \Omega \left( \Phi^L_t\big|\widetilde{\cN}_c
%\right)$ and, by the above proposition, $x\in \widetilde{\cA}_c$. This proves
%(\ref{app.1}).\\

\vspace{10 pt}
\begin{Proof}
%Without any loss of generality we can assume that $c=0$ and $\a(0)=0$.
Let $(x,v) \in \Omega \left( \Phi^L_t\big|\widetilde{\cN}_0\right)$. By the 
definition of non-wandering point, there exist a sequence $(x_k,v_k) \in 
\widetilde{\cN_c}$ and $t_k \rightarrow +\infty$, such that $(x_k,v_k) 
\rightarrow (x,v)$  and $\Phi^L_{t_k}(x_k,v_k)~\rightarrow~\!\!(x,v)$ as 
$k\rightarrow +\infty$.
From (\ref{2.8}),  for each $(x_k,v_k)$ there exists a $\eta$ critical 
subsolution $u_k$, such that the curve $\g_k (t) = \pi\left( \Phi^L_{t} (x_k,v_k)\right)$ is 
$(u_k,L_{\eta},\a(c))$~-~calibrated. Moreover, up to extracting a subsequence, we can 
assume that, on any compact interval, $\g_k$ converge in the $C^1$-topology
to $\g(t)=\pi\left(\Phi^L_t(x,v)\right)$.

Pick now any critical subsolution $u$. If we show that $\g$ is 
$(u,L_{\eta},\a(c))$~-~calibrated, using (\ref{2.8}) we can conclude that $(x,v)\in 
\widetilde{\cA}_c$. First of all, observe that, by the continuity of $u$,
$$u(\g_k(t_k)) - u(x_k)\stackrel{k\rightarrow \infty}{\longrightarrow} 0\;.$$
Using that $\eta$-critical 
subsolutions are equi-Lipschitz (as remarked after definition \ref{defdomination}), we can also conclude that
$$u_k(\g_k(t_k)) - u_k(x_k)\stackrel{k\rightarrow \infty}{\longrightarrow} 0$$
and, therefore, 
\be
\int_0^{t_k} L_{\eta}(\g_k(s),\dot{\g}_k(s)) + \a(c)\,ds = 
u_k(\g_k(t_k)) - u_k(x_k) \stackrel{k\rightarrow \infty}{\longrightarrow} 0\;. 
\ee
Let $0\leq a\leq b$ and choose $t_k \geq b$. Observe now that
\beqano
u(\g_k(b)) - u(\g_k(a)) &=& u(\g_k(t_k)) - u(x_k) - \left[u(\g_k(t_k)) -  u(\g_k(b)) \right] -\\
&& \; -\; \left[u(\g_k(a))- u(x_k) \right] \geq\\
&\geq& u(\g_k(t_k)) - u(x_k) -  \int_b^{t_k} L_{\eta}(\g_k(s),\dot{\g}_k(s)) + \a(c)\,ds - \\
&& \;-\; \int_0^{a} L_{\eta}(\g_k(s),\dot{\g}_k(s))+\a(c)\,ds = \\
&=& u(\g_k(t_k)) - u(x_k) +  \int_a^b L_{\eta}(\g_k(s),\dot{\g}_k(s)) + \a(c)\,ds - \\
&&\;-\;  \int_0^{t_k} L_{\eta}(\g_k(s),\dot{\g}_k(s))+\a(c)\,ds\,;
\eeqano

taking the limit as $k\rightarrow \infty$ on both sides,
one can conclude:
$$
u(\g(b)) - u(\g(a)) \geq \int_a^b L_{\eta}(\g(s),\dot{\g}(s))+\a(c)\,ds
$$
and therefore, from the fact that $u\prec L_{\eta}+\a(c)$, it follows the equality. This shows 
that $\g$ is $(u,L_{\eta},\a(c))$-calibrated on $[0,\infty)$. To show that it is indeed 
calibrated on all $\R$,  one can make a symmetric argument, letting $(y_k,w_k)= 
\Phi^L_{t_k}(x_k,v_k)$ play the role of $(x_k,v_k)$ in the previous argument.
In fact, one has $(y_k,w_k) \rightarrow (x,v)$  and $\Phi^L_{-t_k}(y_k,w_k) 
\rightarrow (x,v)$ as $k\rightarrow +\infty$ and the very same argument works.

\end{Proof}

%%%%%%%%%%%%%%%%%%%%%%%%%%%%%%%%%%%%%%%%%%%%%%%%%%%%%%%%%%%%%%%%%%%%%%%%%%%%%%%%%%%%%%%%%%%%%%%%%%%%
%\nocite{*}
%\addcontentsline{toc}{chapter}{Bibliography}
\bibliographystyle{plain}
\bibliography{biblioLectureNotesMatherTheory}
%{}

%%%%%%%%%%%%%%%%%%%%%%%%%%%%%%%%%%%%%%%%%%%%%%%%%%%%%%%%%%%%%%%%%%%%%

\end{document}